\documentclass[a4paper,11pt]{article}
\usepackage[margin=1.2in]{geometry}
\pdfoutput=1
\usepackage{mathtools}

\usepackage[utf8]{inputenc}
\usepackage{amsfonts,amssymb, bm}
\usepackage{graphicx,color,epstopdf}
\usepackage{hyperref}
\allowdisplaybreaks
\newtheorem{proof}{Proof}
\newtheorem{mylemma}{Lemma}[section]
\newtheorem{mytheorem}{Theorem}[section]
\newtheorem{mycorollary}{Corollary}[section]
\newtheorem{myremark}{Remark}[section]

\newtheorem{myprop}{Proposition}[section]

\def\XXint#1#2#3{{\setbox0=\hbox{$#1{#2#3}{\int}$}
    \vcenter{\hbox{$#2#3$}}\kern-.5\wd0}}

\def\tg {\tilde{g}}
\def\tL {L_c}
\def\tS {S_c}
\def\tx {\tilde{x}}
\def\ty {\tilde{y}}
\def\tG {\tilde{G}}

\def\tu {\tilde{u}}

\def\dd {''}

\def\wtd#1{{\widetilde{#1}}}

\def\cS {{\cal S}}
\def\cK {{\cal K}}

\def\bi{{\bf i}}

\newcommand{\tcr}{\textcolor{black}}

\begin{document}
\title{PML and high-accuracy boundary integral equation solver for wave
  scattering by a locally
  defected periodic surface} \author{Xiuchen Yu$^1$, Guanghui Hu$^2$, Wangtao
  Lu$^3$ and Andreas Rathsfeld$^4$}
\footnotetext[1]{School of Mathematical Sciences, Zhejiang University, Hangzhou
  310027, China. Email: yuxiuchen@zju.edu.cn.}
\footnotetext[2]{School of Mathematical Sciences, Nankai University, Tianjin
  300071, China. Email: ghhu@nankai.edu.cn}

\footnotetext[3]{School of
  Mathematical Sciences, Zhejiang University, Hangzhou 310027, China. Email:
  wangtaolu@zju.edu.cn. This author is partially
  supported by NSF of Zhejiang Province for Distinguished Young Scholars
  (LR21A010001).}
\footnotetext[4]{Weierstrass Institute, Mohrenstr. 39,
10117 Berlin,
Germany. Email: rathsfeld@wias-berlin.de}

\maketitle
\begin{abstract}
  This paper studies the perfectly-matched-layer (PML) method for wave
  scattering in a half space of homogeneous medium bounded by a two-dimensional,
  perfectly conducting, and locally defected periodic surface, and develops a
  high-accuracy boundary-integral-equation (BIE) solver. Along the vertical
  direction, we place a PML to truncate the unbounded domain onto a strip and
  prove that the PML solution converges linearly to the true solution in the
  physical subregion of the strip with the PML thickness. Laterally, we divide
  the unbounded strip into three regions: a region containing the defect and two
  semi-waveguide regions, separated by two vertical line segments. In \tcr{both
  semi-waveguides}, we prove the well-posedness of an associated scattering
  problem so as to well define a Neumann-to-Dirichlet (NtD) operator on the
  associated vertical segment. The two NtD operators, serving as exact lateral
  boundary conditions, reformulate the unbounded strip problem as a boundary
  value problem onto the defected region. Due to the periodicity of the
  semi-waveguides, \tcr{both NtD operators turn out }
to be closely related to a
  Neumann-marching operator, governed by a nonlinear \tcr{Riccati }
  equation. It is
  proved that the Neumann-marching operators are contracting, so that the PML
  solution decays exponentially fast along \tcr{both lateral directions}. The
  consequences \tcr{culminate in two opposite aspects. }
  Negatively, the PML solution \tcr{cannot
  exponentially converge} to the true solution in the whole physical region of
  the strip. Positively, from a numerical perspective, the \tcr{Riccati } 
  equations can
  now be efficiently solved by a recursive doubling procedure and a
  high-accuracy PML-based BIE method so that the boundary value problem on the
  defected region can be solved efficiently and accurately. Numerical
  experiments demonstrate that the PML solution converges exponentially fast to
  the true solution \tcr{in any compact subdomain }
  of the strip.
\end{abstract}
\section{Introduction}
Due to its nearly reflectionless absorption of outgoing waves, perfectly matched
layer or PML, since its invention by Berenger in 1994 \cite{ber94}, has become a
primary truncation technique in a broad class of unbounded wave scattering
problems \cite{chew95,mon03,fic11}, ranging from quantum mechanics, acoustics,
electromagnetism (optics), to seismology. Mathematically, a PML can be
equivalently understood as a complexified transformation of a coordinate
\cite{chewee94}. A wave outgoing along the coordinate is then analytically
continued in the complex plane and becomes exponentially decaying in the PML.
However, it is such a double-edged feature that makes PML be placed only in the
direction where the medium structure is invariant so as to guarantee the
validity of analytic continuation. Consequently, PML loses its prominence for
some complicated structures, such as periodic structures \cite{joh08}. Motivated
by this, this paper studies wave scattering in a half space of homogeneous
medium bounded by a two-dimensional, perfectly conducting, and locally defected
periodic surface, and investigate the potential of PML in designing an accurate
boundary integral equation (BIE) solver for the scattering problem.

Let a cylindrical wave due to a line source, or a downgoing plane wave be
specified above the defected surface. Then, a primary question is to understand
clearly how the scattered wave radiates at infinity. Intrinsically, PML is
highly related to the well-known Sommerfeld radiation condition (SRC), which,
arguably, is an alternative way of saying ``{\it wave is purely outgoing at
  infinity}''. However, SRC is considered to be no longer valid for
characterizing the scattered wave even when the surface is flat \cite{arehoh05}.
Instead, upward propagation radiation condition (UPRC), a.k.a angular spectrum
representation condition \cite{desmar96} is commonly used, and can well pose the
present problem or even more general rough surface scattering problems
\cite{chaels10,chamon05, charoszha99}. Milder than SRC, UPRC only requires that
the scattered wave contain no downgoing waves on top of a straight line above
the surface, allowing waves incoming horizontally from infinity.

If the surface has no defects, the total wave field for the plane-wave incidence
is quasi-periodic so that the original scattering problem can be formulated in a
single unit cell, bounded laterally but unbounded vertically. According to UPRC,
the scattered wave at infinity can then be expressed in terms of upgoing Bloch
waves, so that a transparent boundary condition or PML of a local/nonlocal
boundary condition can be successfully used to terminate the unit cell
vertically; readers are referred to \cite{baodobcox95,chewu03,lulu12josaa,
  zhowu18} and the references therein, for related numerical methods as well as
theories of exponential convergence due to a PML truncation. But, if the
incident wave is nonquasi-periodic, e.g., the cylindrical wave, or if the
surface is locally defected, much fewer numerical methods or theories have been
developed as it is no longer straightforward to laterally terminate the
scattering domain. Existing laterally truncating techniques include recursive
doubling procedure (RDP) \cite{yualu07,ehrhanzhe09}, Floquet-Bloch mode expansion \cite{helpre98,hulu09,leczha17},
and \tcr{Riccati}
-equation based exact boundary condition \cite{jollifli06}.

In a recent work \cite{hulurat20}, we proved that the total field for
the cylindrical incidence, a.k.a the Green function, satisfies the standard SRC
on top of a straight line above the surface. Based on this, we further
revealed that for the \tcr{plane-wave} incidence, the perturbed part of the total field due
to the defect satisfies the SRC as well. Consequently, this suggests to use a
PML to terminate the vertical variable so as to truncate the unbounded domain to
a strip, bounded vertically but unbounded laterally. In fact, such a natural
setup of PML had already been adopted in the literature \cite{yualu07,chamon09,sunzhe09}, without a
rigorous justification of the outgoing behavior, though. It is worthwhile to
mention that \tcr{Chandler-Wilde} and Monk in \cite{chamon09} rigorously proved that under a
Neumann-condition PML, the PML solution converges to the true
solution in the whole physical region of the strip at the rate of only algebraic
order of PML thickness; they further revealed that the PML solution due to the
cylindrical incidence for a flat surface decays exponentially at infinity of a
rectangular strip.
However, it remains unclear how the PML solution radiates at infinity of the
more generally curved strip under consideration. On the other hand, no literally
rigorous theory has been developed to clearly understand why this PML-truncated
strip can further be laterally truncated to a bounded domain by the
aforementioned techniques without introducing artificial ill-posedness; in other
words, the well-posedness of scattering problems in exterior regions of the
truncated domain is unjustified.

To address these questions, we first prove in this paper that under a
Dirichlet-condition PML, the PML solution due to the cylindrical incidence,
i.e., \tcr{Green's function of the strip}, converges to the true solution in the physical
subregion of the strip at an algebraic order of the PML thickness. Next, we
split the strip into three regions: a bounded region containing the defect and
two semi-waveguide regions of \tcr{a }single-directional periodic surface, separated by
two vertical line segments. By use of \tcr{Green's function of the strip}, transparent
boundary conditions can be developed to truncate the unbounded semi-waveguides.
Based on this, we apply the method of variational formulation and Fredholm
alternative to prove the well-posedness of the scattering problem in
either semi-waveguide so as to define a Neumann-to-Dirichlet (NtD) operator on
its associated vertical segment. The two NtD operators serve exactly as lateral
boundary conditions to terminate the strip and to reformulate the unbounded
strip problem as a boundary value problem on the defected region. Due to the
periodicity of the semi-waveguides, \tcr{both NtD operators turn out} to be closely
related to a Neumann-marching operator, which solves a nonlinear \tcr{Riccatti}
equation. It is proved that the Neumann-marching operators are contracting,
indicating that the PML solution decays exponentially fast along \tcr{both lateral
directions} even for the curved strip. The consequences \tcr{culminate in two opposite aspects}. Positively, from a numerical perspective, the \tcr{Riccati} equations can be
efficiently solved by \tcr{an} RDP method so that the strip can
be laterally truncated with ease. Negatively, the
PML solution shall never exponentially converge to the true solution in the
whole physical region of the strip. Nevertheless, as conjectured in
\cite{chamon09}, exponential convergence is optimistically expected to be
realizable in any compact subdomain of the strip.

To validate the above conjecture numerically, we employ a high-accuracy
PML-based boundary integral equation (BIE)
method \cite{luluqia18} to execute the RDP so that the two \tcr{Riccati} equations can
be accurately solved for the two Neumann-marching operators, respectively, and
hence the two NtD operators terminating the strip can be obtained. With the two
NtD operators well-prepared, the boundary value problem in the defected region
can be accurately solved by the PML-based BIE method again. By carrying out several
numerical experiments, we observe that the PML truncation error for wave field
\tcr{over} the defected part of the surface decays exponentially fast as PML absorbing
strength or thickness increases. \tcr{This} indicates that there is a chance that the
PML solution still converges to the true solution exponentially in
any compact subdomain of the strip, the justification of which remains open.

The rest of this paper is organized as follows. In section 2, we introduce the
half-space scattering problem and present some known well-posedness results. In
section 3, we introduce a Dirichlet-condition PML, prove the well-posedness of
the PML-truncated problem and study the prior error estimate of the PML
truncation. In section 4, we study well-posedness of the semi-waveguide
problems. In section 5, we establish lateral boundary conditions, prove the
exponentially decaying property of the PML solution at infinity of the strip,
and develop \tcr{an} RDP technique to get the lateral boundary conditions. In
section 6, we present a PML-based BIE method to numerically solve the
scattering problem. In section 7, numerical experiments are carried out to
demonstrate the performance of the proposed numerical method and to validate the
proposed theory. We draw our conclusion finally in section 8 and propose
some future plans.

\section{Problem formulation} 
\begin{figure}[!ht]
  \centering
  (a)\includegraphics[width=0.341\textwidth]{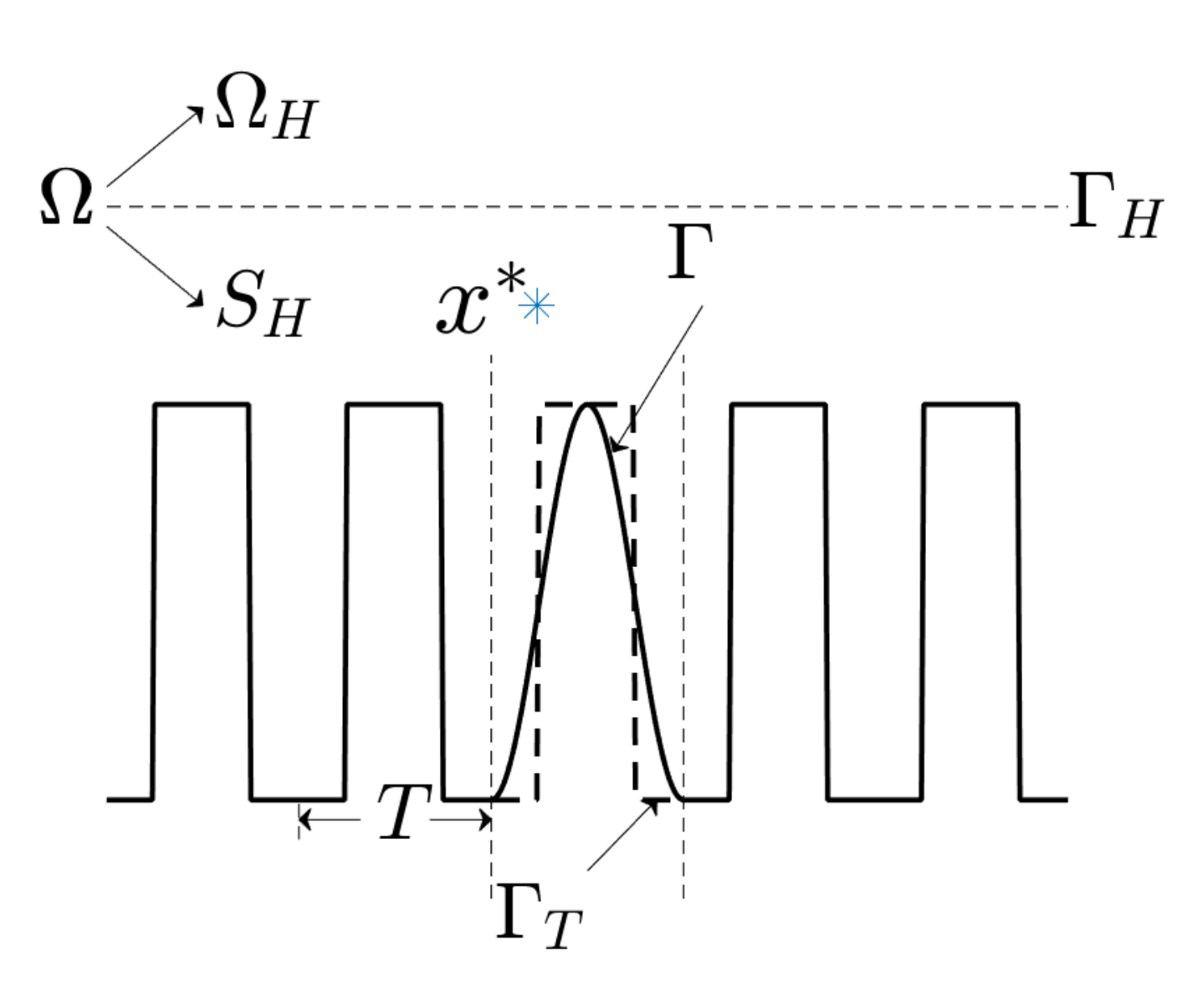}
  (b)\includegraphics[width=0.341\textwidth]{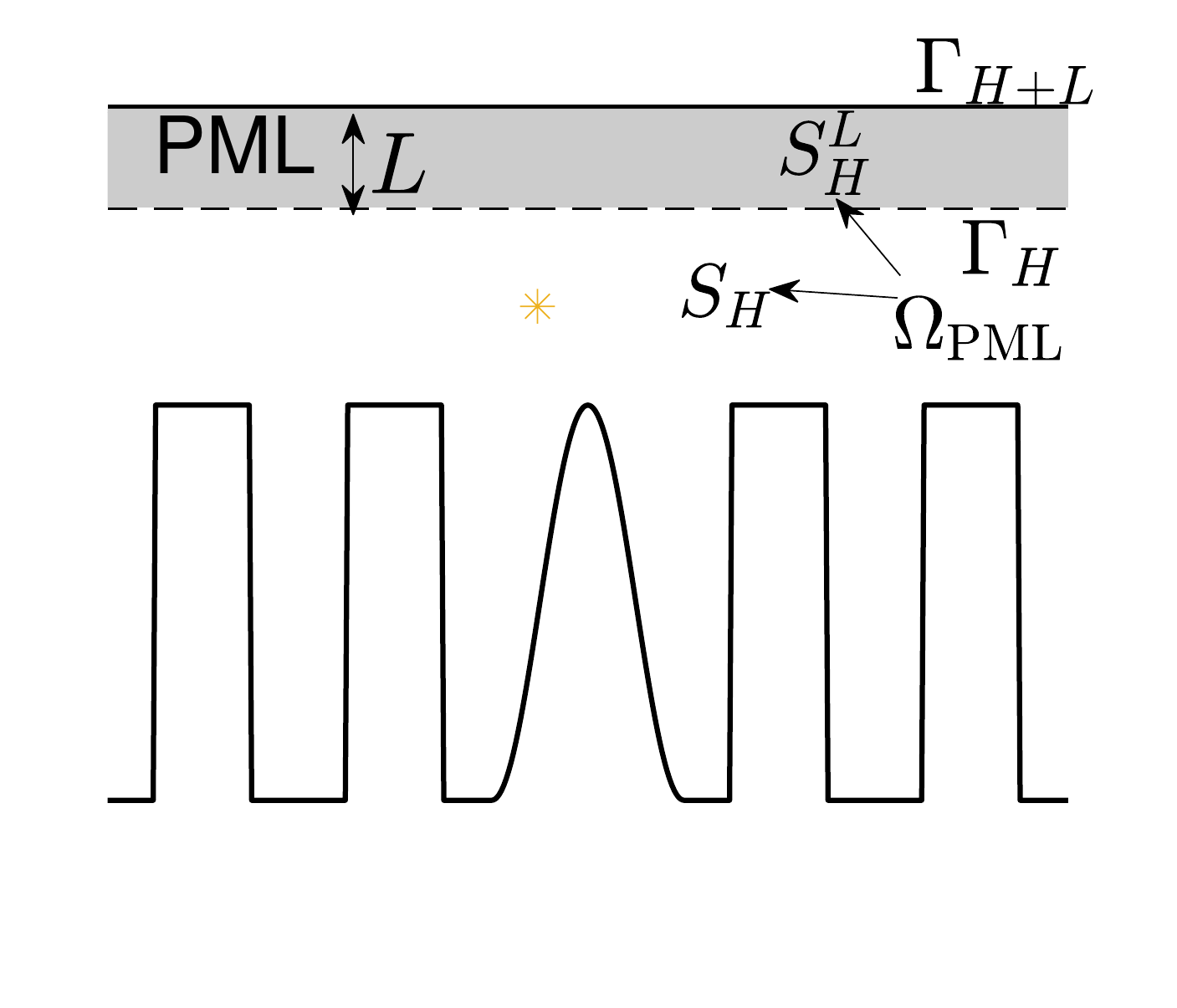}
	\caption{(a) A sketch of the half-space scattering problem. (b) A PML placed
    above $\Gamma$. The scattering surface $\Gamma$ locally perturbs the
    periodic curve $\Gamma_T$ of period $T$. $x^*$ represents the exciting
    source. $\Gamma_H$ is an artificial interface, on which a DtN map is defined
  or a PML is placed.}
  \label{fig:org:problem}
\end{figure}

{\color{red} Let $\Omega\times
  \mathbb{R}\subset\mathbb{R}^3$ be an $x_3$-invariant domain bounded by a
  perfectly-conducting surface $\Gamma\times\mathbb{R}$, where $\Gamma\subset
  \mathbb{R}^2$, bounding domain $\Omega\subset\mathbb{R}^2$, is a local perturbation of a $T$-periodic curve $\Gamma_T\subset
  \mathbb{R}^2$ periodic in $x_1$-direction, as shown in
  Figure~\ref{fig:org:problem}(a).}
\tcr{We denote the Cartesian coordinate system of
  $\mathbb{R}^3$ by $(x_1,x_2,x_3)$ and let $x=(x_1,x_2)\in\Omega$.}
Throughout this paper, we shall assume that $\Gamma$ is Lipschitz and that
$\Omega$ satisfies the following geometrical condition
\begin{equation*}
  {\rm (GC1): }\quad (x_1,x_2)\in \Omega\Rightarrow (x_1,x_2+a)\in \Omega,\quad\forall a\geq 0.
\end{equation*}
For simplicity, suppose $\Gamma$ only perturbs one periodic part of
$\Gamma_T$, say $|x_1|<\frac{T}{2}$.

Let the unbounded domain \tcr{$\Omega\times\mathbb{R}$} be filled by a homogeneous medium of
refractive index $n$. For a time-harmonic transverse-electric (TE) polarized electro-magnetic wave, of
time dependence $e^{-\bi\omega t}$ for the angular frequency $\omega$, the
$x_3$-component of the electric field, denoted by $u^{\rm tot}$, \tcr{is
  $x_3$-invariant and} satisfies the
following two-dimensional (2D) Helmholtz equation
\begin{align}
  \label{eq:gov:utot}
  \Delta u^{\rm tot} + k^2 u^{\rm tot} =& 0,\quad{\rm on}\quad \Omega,\\
  \label{eq:gov:bc:te}
  u^{\rm tot}=& 0,\quad{\rm on}\quad\Gamma,
\end{align}
where $\Delta = \partial_{x_1}^2+\partial_{x_2}^2$ is the 2D Laplacian and
$k=k_0n$ with $k_0=\frac{2\pi}{\lambda}$ denoting the free-space wavenumber for
wavelength $\lambda$.

Let an incident wave $u^{\rm inc}$ be specified in $\Omega$ \tcr{and let
  $x=(x_1,x_2)\in\Omega$}. In this paper, we shall mainly focus on the following
two cases of incidences: (i) a plane wave $u^{\rm inc}(x)=e^{\bi k (\cos\theta
  x_1 -\sin\theta x_2)}$ for the incident angle $\theta\in(0,\pi)$; (ii) a
cylindrical wave $u^{\rm inc}(x;x^*) =
G(x;x^*)=\frac{\bi}{4}H_0^{(1)}(k|x-x^*|)$ \tcr{excited by a source at
$x^*=(x_1^*,x_2^*)\in \Omega$}. In the latter case, equation (\ref{eq:gov:utot})
should be replaced by
\begin{equation}
  \label{eq:utot:ps}
  \Delta u^{\rm tot} + k^2 u^{\rm tot} = -\delta(x-x^*),
\end{equation}
so that $u^{\rm tot}(x;x^*)$ in fact represents the Green function excited by
the source point $x^*$. For simplicity, we assume that $|x_1^*|<T/2$ so
that $x^*$ is right above the perturbed part of $\Gamma$.

Let $u^{\rm sc}=u^{\rm tot}-u^{\rm inc}$ denote the scattered wave. One may enforce the following
UPRC:
\begin{equation}
  \label{eq:uprc}
  u^{\rm sc}(x) = 2\int_{\Gamma_H} \frac{\partial G(x;y)}{\partial y_2} u^{\rm sc}(y) ds(y),
\end{equation}
where $\Gamma_H=\{(x_1,H):x_1\in\mathbb{R}\}$ denotes a straight line strictly
above $\Gamma$ for some $H>0$ \tcr{and $y=(y_1,y_2)$}. According to \cite{chamon05}, the UPRC helps to define a
Dirichlet-to-Neumann map ${\cal T}: H^{1/2}(\Gamma_H)\to H^{-1/2}(\Gamma_H)$ for
the domain $\Omega_H=\{x\in \Omega:x_2>H\}$, such that for any $\phi\in H^{1/2}(\Gamma_H)$, 
\begin{equation}
  \label{eq:def:T}
  {\cal T}\phi = {\cal F}^{-1}M_z\hat{\phi},
\end{equation}
where $\hat{\phi}(H;\xi)=[{\cal F}\phi](H;\xi)$ denotes the following normalized Fourier
transform
\begin{equation}
  \tcr{[{\cal F}\phi](H;\xi)} = \frac{1}{\sqrt{2\pi}}\int_{\mathbb{R}}\phi(x_1,H) e^{-\bi\xi x_1}dx_1,
\end{equation}
and \tcr{the operator $M_z$ in the space of Fourier transforms is the operator of multiplication by} 
\begin{equation}
  \label{eq:def:z}
  z(\xi) = \left\{
    \begin{array}{ll}
      -\bi\sqrt{k^2-\xi^2},&{\rm for}\ |\xi|\leq k,\\
      \sqrt{\xi^2-k^2},&{\rm for}\ |\xi|>k.
    \end{array}
  \right.
\end{equation}
Then, we may enforce
\begin{equation}
  \label{eq:dtn}
  \partial_{\nu} u^{\rm sc} = -{\cal T}u^{\rm sc},\quad{\rm on}\quad \Gamma_H,
\end{equation}
where, unless otherwise indicated, $\nu$ always denotes the outer unit normal
vector on $\Gamma_H$. The UPRC guarantees the well-posedness of our scattering
problem \cite{chamon05}, but \tcr{allows} $u^{\rm sc}$ containing incoming waves,
largely limiting its applications in designing numerical algorithms.
Nevertheless, our recent work \cite{hulurat20} has shown a stronger
Sommerfeld-type condition for the aforementioned two incidences, which still
preserves the well-posedness. Note that \cite{hulurat20} assumes further the
following condition:
\[
\textrm {(GC2): some (and hence any) period of $\Gamma_T$ contains a line segment,}
\]
which guarantees a local behavior of the Green function $u^{\rm tot}(x;y)$ for
any $x, y$ sufficiently close to each line segment. \tcr{Let $S_H=\Omega\cap
  \{x:x_2<H\}$ be the strip between $\Gamma_H$ and $\Gamma$.} The radiation
condition reads as follows:
\begin{itemize}
\item[(i).] For the \tcr{plane-wave} incidence, $u^{\rm og}:=u^{\rm tot}-u^{\rm
    tot}_{\rm ref}$, \tcr{where $u^{\rm tot}_{\rm ref}$ is the reference
    scattered field for the unperturbed \tcr{scattering} curve
    $\Gamma=\Gamma_T$,} satisfies the following half-plane Sommerfeld radiation
  condition (hSRC): \tcr{for some sufficiently large $R>0$ and any $\rho<0$},
  \begin{equation}
    \label{eq:src}
    \tcr{\lim_{r\to\infty}\sup_{\alpha\in[0,\pi]}\sqrt{r}\left|  \partial_{r}u^{\rm og}(x)-\bi ku^{\rm og}(x)\right| = 0,\ \sup_{r\geq R}r^{1/2}|u^{\rm og}(x)|<\infty,\ {\rm and}\ u^{\rm og}\in H_\rho^{1}(S_H^R),}
  \end{equation}
  where $x=(r\cos\alpha,H+r\sin\alpha)$, \tcr{$S_H^R=S_H\cap\{x:|x_1|>R\}$, and
    $H_\rho^1(\cdot)=(1+x_1^2)^{-\rho/2}H^1(\cdot)$ denotes a weighted Sobolev
    space}. We defer the computation of $u^{\rm tot}_{\rm ref}$ to section 6.3.
\item[(ii).] For the cylindrical incidence, the total field $u^{\rm og}:=u^{\rm tot}$ itself satisfies the hSRC
  (\ref{eq:src}) in $\Omega_H$. Thus, the scattered field $u^{\rm sc}$ satisfies (\ref{eq:src}) as
  well since $u^{\rm inc}$ satisfies (\ref{eq:src}).
\end{itemize}
Certainly, $u^{\rm og}$ satisfies the UPRC condition (\ref{eq:uprc}) such that
(\ref{eq:dtn}) holds for $u^{\rm og}$ in place of $u^{\rm sc}$ \cite[Them.
2.9(ii)]{chazha98}. In the following, we shall consider the cylindrical
incidence only and the \tcr{plane-wave incidence} case can be analyzed
similarly.

We recall some important results from \cite{chamon05}. To remove the singularity
of the right-hand side of (\ref{eq:utot:ps}), let
\begin{equation}
  \label{eq:uogr}
  u^{\rm og}_r(x;x^*)=u^{\rm og}(x;x^*)-\chi(x;x^*) u^{\rm inc}(x;x^*),
\end{equation}
\tcr{where} the cut-off function $\chi(x;x^*)=1$ in a neighborhood of $x^*$ and
\tcr{has} a sufficiently small support enclosing $x^*$. Let
$V_H=\{\phi|_{S_H}:\phi\in H_0^1(\Omega)\}$. Then, it is equivalent to seek
$u^{\rm og}_r\in V_H$ that \tcr{satisfies} the following boundary value problem
\begin{align}
  \Delta u^{\rm og}_r + k^2 u^{\rm og}_r =& g,\quad{\rm on}\quad S_H,\\
  \partial_{\nu}u^{\rm og}_r=& -{\cal T}\tcr{u^{\rm og}_r},\quad{\rm on}\quad \Gamma_H, 
\end{align}
where \tcr{$g = -[\Delta\chi]u^{\rm
    inc}-2\sum_{j=1}^2\partial_{x_j}\chi\partial_{x_j}u^{\rm inc}\in L^2(S_H)$}
  such that \tcr{${\rm supp}\ g$ is in the neighborhood of $x^*$ contained in
    $\overline{S_H}$}. An equivalent variational formulation reads as follows:
Find $u_r^{\rm og}\in V_H$, such that for any $\phi\in V_H$,
\begin{equation}
  b(u_r^{\rm og},\phi) = -(g,\phi)_{S_H},
\end{equation}
where the sesqui-linear form $b(\cdot,\cdot): V_H\times V_H\to \mathbb{C}$ is given by
\begin{equation}
  b(\phi,\psi) = \int_{S_H}(\nabla\phi\cdot \nabla\bar{\psi} - k^2\phi\bar{\psi})dx + \int_{\Gamma_H}{\cal T} \phi \bar{\psi} ds.
\end{equation}
It has been shown in \cite{chamon05} that $b$ satisfies the following inf-sup
condition: for all $v\in V_H$,
\begin{equation}
  \label{eq:b:infsup}
  \gamma ||v||_{V_H}\leq \sup_{\phi\in V_H}\frac{|b(v,\phi)|}{||\phi||_{V_H}},
\end{equation}
where $\gamma>0$ depends on $H$, $k$ and $\Omega$. Furthermore, $b$ defines an
invertible operator ${\cal A}: V_H\to V_H^*$ such that $({\cal
  A}\phi,\psi)=b(\phi,\psi)$ and $||{\cal A}^{-1}||\leq \gamma^{-1}$. Thus,
$u_r^{\rm og}=-{\cal A}^{-1}g$ so that $u^{\rm og}=-{\cal A}^{-1}g+\chi u^{\rm inc}$.

The hSRC (\ref{eq:src}) suggests to compute the outgoing wave $u^{\rm og}$
numerically, as the PML technique\cite{ber94,chamon09} could apply now to
truncate the $x_2$-direction. In the following sections, we shall first
introduce the setup of a PML to truncate $x_2$ and then develop an accurate
lateral boundary condition to truncate $x_1$.

\section{PML Truncation}
Mathematically, the PML truncating $x_2$ introduces a complexified coordinate
transformation
\begin{equation}
  \label{eq:x2:pml}
  \tx_2 = x_2 + \bi S\int_0^{x_2}\sigma(t)dt,
\end{equation}
where $\sigma(x_2)=0$ for $x_2\leq H$ and $\sigma(x_2)\geq 0$ for $x_2\geq H$;
\tcr{note that such a tilde notation can also be used to define $\ty_2$ and
  $\tx_2^*$ in the following}. As shown in Figure~\ref{fig:org:problem}(b), the
planar strip $S_H^L=\mathbb{R}\times [H,H+L]$ with nonzero $\sigma$ is called
the PML region so that $L$ represents its thickness. \tcr{In this paper, we
  choose an $m\geq 0$ and},
\begin{equation}
	\label{eq:sigma}
	\sigma(x_2) = \left\{
    \begin{array}{lc}
    \frac{2 f_2^m}{f_1^m+f_2^m}, & x_2\in[H,H+L/2]\\
    2, & \tcr{x_2\geq H+L/2, m\neq 0}\\
    1, & \tcr{x_2\geq H+L/2, m=0}\\
    \end{array}
     \right.
\end{equation}
where we note that $\sigma\equiv 1$ \tcr{if} $m=0$, and
\[
	f_1=\left(\frac{1}{2}-\frac{1}{\tcr{m}}\right)\xi^3+\frac{\xi}{\tcr{m}}+\frac{1}{2},\quad
	f_2 = 1-f_1,\quad \xi = \frac{2 x_2 - (2H+L/2)}{L/2}.
\]
\tcr{Let $\tL \ :=\   \tx_2(H+L) - H=L+\bi \tS L$,
where $\tS = \frac{S}{L}\int_{H}^{H+L}\sigma(t)dt$.} Both the real part and
imaginary part of $\tL$ affect the absorbing strength of the PML
\cite{chewee94}.

Now, let $\tx=(x_1,\tx_2)$. \tcr{For $x^*\in S_H$, according to
(\ref{eq:uprc}), we can define by analytic continuation that
\[
  u^{\rm og}(\tx;x^*) := 2\int_{\Gamma_H} \frac{\partial G(\tx;y)}{\partial y_2} u^{\rm og}(y;x^*) ds(y), 
\]
that satisfies
\[
  \tilde{\Delta} u^{\rm og}(\tx;x^*)+ k^2u^{\rm og}(\tx;x^*) = -\delta(x-x^*),
\]
where $\tilde{\Delta}=\partial_{x_1}^2+\partial_{\tx_2}^2$. } By chain rules, we
see that $\tilde{u}^{\rm og}(x;x^*):=u^{\rm og}(\tilde{x};x^*)$ \tcr{satisfies}
\begin{align}
	\label{eq:helm:pml}
        \nabla\cdot({\bf A} \nabla \tilde{u}^{\rm og}) + k^2\alpha\tilde{u}^{\rm og} =& -\delta(x-x^*), \quad{\rm on\ }\Omega_{\rm PML},\\
  \label{eq:bc1}
        \tilde{u}^{\rm og} =& 0,\quad {\rm on}\quad \Gamma.
\end{align}
where ${\bf A}={\rm diag}\{\alpha(x_2), 1/\alpha(x_2))\}$,
\tcr{$\alpha(x_2)=1+\bi S \sigma(x_2)$}, and the \tcr{PML} region $\Omega_{\rm
  PML}=\Omega\cap\{x:x_2\leq H+L\}$ consists of the physical region $S_H$ and
the PML region $S_H^L$. On the PML boundary $x_2=H+L$, we use the
\tcr{homogeneous} Dirichlet boundary condition
\begin{equation}
  \label{eq:pml:bc}
  \tu^{\rm og} = 0,\quad{\rm on}\quad \Gamma_{H+L}=\{x:x_2=H+L\}.
\end{equation}
\tcr{The authors} in \cite{chamon09} adopted a Neumann condition on the PML
boundary $\Gamma_{H+L}$ and proved the well-posedness of the related PML
truncation problem. Here, we choose the Dirichlet condition (\ref{eq:pml:bc})
since, as we shall see, our numerical results indicate that the Dirichlet-PML
seems more stable than the Neumann-PML. \tcr{Furthermore}, we need \tcr{Green's
  function of the strip} $\tilde{u}^{\rm og}(x;x^*)$ for any $x^*\in\Omega_{\rm
  PML}$ but not limited to $S_H$ to establish lateral boundary conditions. For
completeness, we shall, following the idea of \cite{chamon09}, study the
well-posedness of the problem (\ref{eq:helm:pml}-\ref{eq:pml:bc}) \tcr{for any
  $x^*\in\Omega_{\rm PML}$}.

The fundamental solution of the anisotropic Helmholtz equation
(\ref{eq:helm:pml}) is \cite{luluqia18}
\begin{equation}
	\label{eq:fun:helm2}
	\tilde{G}(x;y) = G(\tx;\ty) =
	\frac{\bi}{4} H_0^{(1)}[k\rho(\tx;\ty)],
\end{equation}
where \tcr{$\ty=(y_1,\ty_2)$}, the complexified distance function
$\rho$ is defined to be 
\begin{equation}
	\label{eq:comp:dist}
	\rho(\tilde{x},\tilde{y}) = [(x_1-y_1)^2 + (\tilde{x}_2-\tilde{y}_2)^2]^{1/2},
\end{equation}
and the half-power operator $z^{1/2}$ is chosen to be the branch of
$\sqrt{z}$ with nonnegative real part for $z\in
\mathbb{C}\backslash(-\infty,0]$ such that ${\rm arg}(z^{1/2})\in[0,\pi)$. The special choice of $\sigma$ in (\ref{eq:sigma}) ensures that
\begin{equation}
 \tG(x;y) = \tG(x;y_{\rm imag}),
\end{equation}
for any $x\in\Gamma_{H+L}$, \tcr{when} $y=(y_1,y_2)$ and $y_{\rm
  imag}=(y_1,2(H+L)-y_2)$, \tcr{the mirror image of $y$ w.r.t line $\Gamma_{H+L}$}, are
sufficiently close to $\Gamma_{H+L}$ so that $\rho(\tx;\ty)= \rho(\tx;\ty_{\rm
  imag})$.

To remove the singularity of the right-hand side of (\ref{eq:helm:pml}), we
introduce
\begin{equation}
  \label{eq:tuogr}
\tu^{\rm og}_r(x;x^*)=\tu^{\rm og}(x;x^*)-\chi(x;x^*) \tu^{\rm inc}(x;x^*),
\end{equation}
with the same cut-off function $\chi$ as in \eqref{eq:uogr}, where $\tu^{\rm inc}(x;x^*)=u^{\rm inc}(\tx;\tx^*)$. Then,
$\tu^{\rm og}_r$ \tcr{satisfies}
\begin{align}
  \label{eq:pm:u1}
 &\nabla\cdot({\bf A} \nabla \tu^{\rm og}_r) + k^2\alpha \tu^{\rm og}_r =
\tg^{\rm inc},\quad{\rm on}\quad \Omega_{\rm PML},\\
  \label{eq:pm:u2}
  &\tu^{\rm og}_r=0,\quad{\rm on}\quad \Gamma,\\
  \label{eq:pm:u3}
  &\tu^{\rm og}_r=0,\quad{\rm on}\quad \Gamma_{H+L},
\end{align}
where $\tg^{\rm inc} = [\nabla\cdot({\bf A} \nabla ) + k^2\alpha]
(1-\chi(x;x^*)) \tu^{\rm inc}(x;x^*)\in L^2(\Omega_{\rm PML})$ with $\tcr{{\rm
    supp}\ \tg^{\rm inc}}\subset\overline{\Omega_{\rm
    PML}}=\overline{S_H}\cup\overline{S_H^L}$. \tcr{Considering that $x^*$ can
  be situated in $S_H^L$, ${\rm supp}\ \tg^{\rm inc}$ may not completely lie in
  the physical domain $S_H$.} To establish a Dirichlet-to-Neumann map on
$\Gamma_H$ like (\ref{eq:dtn}), we need to study the following boundary value
problem in the PML strip $S_H^L$: given $q\in H^{1/2}(\Gamma_H)$, $s\in
H^{1/2}(\Gamma_{H+L})$, and $\tg^{\rm inc}_{\rm PML}=\tg^{\rm inc}|_{S_H^L}\in
L^2(S_H^L)$ with $\tcr{{\rm supp}\ \tg^{\rm inc}_{\rm PML}}\subset\overline{S_H^L}$, find $v\in
H^1(S_H^L)$ such that
\begin{align}
 &\nabla\cdot({\bf A} \nabla v) + k^2\alpha v =
\tg_{\rm PML}^{\rm inc},\quad{\rm on}\quad S_H^L,\\
  &v=q,\quad{\rm on}\quad \Gamma_H,\\
  &v=s,\quad{\rm on}\quad \Gamma_{H+L}.
\end{align}
Let
\[
v_0(x)=v(x)-\int_{S_H^L}[\tG(x;y)-\tG(x;y_{\rm imag})] \tg^{\rm inc}_{\rm
  PML}(y)dy:=v(x)-v_{{\rm PML}}^{\rm inc}(x),
\]
where we recall that $y_{\rm imag}$ is the \tcr{ mirror image of $y$ w.r.t line
  $ \Gamma_{L+H}$}. Thus, $v_0$ \tcr{satisfies}
\begin{align}
 &\nabla\cdot({\bf A} \nabla v_0) + k^2\alpha v_0 = 0,\quad{\rm on}\quad S_H^L,\\
  &v_0=q_n,\quad{\rm on}\quad \Gamma_H,\\
  &v_0=s_n,\quad{\rm on}\quad \Gamma_{H+L},
\end{align}
where $q_n = q-v_{\rm PML}^{\rm inc}|_{\Gamma_H}\in H^{1/2}(\Gamma_H)$ and $s_n =
s-v_{\rm PML}^{\rm inc}|_{\Gamma_{H+L}}\in H^{1/2}(\Gamma_{H+L})$.

Looking for $v_0$ in terms of only complexified plane waves, we get
\begin{equation}
  \hat{v}_0(x_2;\xi) = [{\cal F}v_0](x_2;\xi) = A \exp(z(\xi)(\tx_2-H)) + B\exp (-z(\xi)(\tx_2-H)),
\end{equation}
\tcr{where we recall that $z$ has been defined in (\ref{eq:def:z}),} 
\begin{align*}
  A(\xi) = \frac{\hat{s}_n(\xi)-\exp(-z\tcr{(\xi)}\tL)\hat{q}_n(\xi)}{\exp(z\tcr{(\xi)}\tL)-\exp(-z\tcr{(\xi)}\tL)},\quad{\rm and}\quad
  B(\xi) = \frac{-\hat{s}_n(\xi)+\exp(z\tcr{(\xi)}\tL)\hat{q}_n(\xi)}{\exp(z\tcr{(\xi)}\tL)-\exp(-z\tcr{(\xi)}\tL)},
\end{align*}
$\hat{s}_n(\xi)=[{\cal F}s_n]\tcr{(H+L;\xi)}$ and $\hat{q}_n(\xi)=[{\cal F}q_n]\tcr{(H;\xi)}$. Here, to make $A$ and $B$
well-defined, we could let $\xi$ travel through a Sommerfeld integral path
$-\infty+0\bi\to 0\to\infty-0\bi$ instead of $\mathbb{R}$ \cite{lu21} such that
$z\neq 0$. 
Consequently, 
\begin{equation}
  \left.-\frac{\partial \hat{v}_0}{\partial x_2}\right|_{x_2=H} = z\frac{-2}{\exp(z\tL)-\exp(-z\tL)}\hat{s}_n + z\frac{\exp(z\tL)+\exp(-z\tL)}{\exp(z\tL)-\exp(-z\tL)}\hat{q}_n.
\end{equation}
Now define two bounded operators ${\cal T}_p: H^{1/2}(\Gamma_H)\to H^{-1/2}(\Gamma_H)$ by 
\[
  \tcr{{\cal F}[{\cal T}_pq_n](H;\xi)} = z\frac{\exp(z\tL)+\exp(-z\tL)}{\exp(z\tL)-\exp(-z\tL)}\hat{q}_n,
\]
and ${\cal N}_p: H^{1/2}(\Gamma_{H+L})\to H^{-1/2}(\Gamma_H)$ by
\[
  \tcr{{\cal F}[{\cal N}_ps_n](H+L;\xi)}= z\frac{-2}{\exp(z\tL)-\exp(-z\tL)}\hat{s}_n;
\]
note that the above definitions allow $\xi\in\mathbb{R}$ now, since limits can
be considered when $z=0$. Returning back to the PML-truncated problem
(\ref{eq:pm:u1}-\ref{eq:pm:u3}), we reformulate it as an equivalent boundary
value problem on the physical region $S_H$: Find $\tu^{\rm og}_r\in V_H$ that \tcr{satisfies}
\begin{align}
 \nabla\cdot({\bf A} \nabla \tu^{\rm og}_r) + k^2\alpha \tu^{\rm og}_r =&
\tcr{\tg^{\rm inc}|_{S_H}},\quad{\rm on}\quad S_H,\\
  \label{eq:dtn:pml}
  \partial_{\nu} \tu^{\rm og}_r =& -{\cal T}_{p} \tu^{\rm og}_r|_{\Gamma_H} +f_p,\quad\tcr{{\rm on}\quad \Gamma_H},
\end{align}
where
\begin{equation}
  f_p = {\cal N}_p (v_{\rm PML}^{\rm inc}|_{\Gamma_{H+L}}) +{\cal T}_p (v_{\rm PML}^{\rm inc}|_{\Gamma_H}) + \partial_{\nu} v_{\rm PML}^{\rm inc}|_{\Gamma_H}\in H^{-1/2}(\Gamma_H).
\end{equation}
The associated variational formulation reads as follows: Find $\tu^{\rm og}_r\in V_H$, such
that for any $\psi\in V_H$,
\begin{equation}
  \label{eq:vf:pml}
  b_p(\tu^{\rm og}_r,\psi) = -\int_{S_H}\tg^{\rm inc}|_{S_H}\bar{\psi}dx+\int_{\Gamma_H}f_p\bar{\psi}ds,
\end{equation}
where the sesquilinear form $b_p(\cdot,\cdot): V_H\times V_H\to \mathbb{C}$ is
given by
\begin{equation}
  \label{eq:def:bp}
  b_p(\phi,\psi) = \int_{S_H}(\nabla\phi\cdot\nabla\bar{\psi}-k^2\phi\bar{\psi})dx + \int_{\Gamma_H}\bar{\psi}{\cal T}_p\phi ds.
\end{equation}
\tcr{As in \cite{chamon09}, we define the following $k$-dependent norm
\[
  ||\phi||_{H^{s}(\mathbb{R})}^2 = \int_{\mathbb{R}}(k^2+\xi^2)^{s}|[{\cal F}\phi](\xi)|^2d\xi
\]
for $H^s(\mathbb{R})$}. Then, the following lemma characterizes a rough \tcr{difference} of ${\cal T}_p$ away from ${\cal T}$.
\begin{mylemma}
    \label{lem:difTTp}
  We have for any \tcr{$k\tS L>0$},
  \begin{equation}
    ||{\cal T}-{\cal T}_p||\leq \frac{1}{\tcr{k\tS L}}.
  \end{equation}
  \begin{proof}
    By a simple analysis, it can be seen that 
    \begin{align*}
      ||{\cal T}-{\cal T}_p||=&\sup_{\xi\in\mathbb{R}}\frac{|z(\xi)|}{\sqrt{k^2+\xi^2}}|1-\coth(z(\xi)\tL)|\\
      =&\sup_{\xi\in\mathbb{R}}\frac{2|z(\xi)\exp(-2z\tcr{(\xi)}\tL)|}{\sqrt{k^2+\xi^2}|1-\exp(-2z\tcr{(\xi)}\tL)|} =\max\{S_1,S_2\},
    \end{align*}
    where we recall that $\tL = L + \bi \tS L$,
    \begin{align*}
      S_1 =& \sup_{0\leq t\leq 1}\frac{2t\exp(-2tk\tcr{\tS}L)}{\sqrt{2-t^2}\sqrt{1+\exp(-4tkL)-2\cos(2tk\tS L)\exp(-2tkL)}}\\
      =&\sup_{0\leq t\leq 1}\frac{2t\exp(-2tk\tcr{\tS}L)}{\sqrt{2-t^2}\sqrt{(1-\exp(-2tkL))^2+4\exp(-2tkL)\sin^2(tk\tS L)}}
    \end{align*}
    and
    \begin{align*}
      S_2=&\sup_{t\geq 1}\frac{2t\exp(-2tk\tcr{\tS}L)}{\sqrt{2+t^2}\sqrt{(1-\exp(-2tkL))^2+4\exp(-2tkL)\sin^2(tk\tS L)}}.
    \end{align*}
    Clearly, $S_2\leq 2\exp(-2kL)$. Since for $t\geq 0$,
    \[
      \tcr{f(t)=\frac{t\exp(-2tk\tS L)}{1-\exp(-2tk\tS L)}}
    \]
    is nonincreasing, it is easy to see that \tcr{$S_1\leq
      2f(0)=\frac{1}{k\tS L}$}.
  \end{proof}
\end{mylemma}

We do not intend to study the relation of $||{\cal T}_p-{\cal T}||$ and the
other parameter $\tS$, as was done in \cite{chamon09} to optimize the
performance of the PML, since the estimate in Lemma~\ref{lem:difTTp} is
enough. Clearly, the sesquilinear form $b_p$ in
(\ref{eq:def:bp}) defines a bounded linear functional ${\cal A}_p: V_H\to V_H^*$
such that: for any $\phi\in V_H$,
\[
  (({\cal A}-{\cal A}_p)\phi,\psi) = b(\phi,\psi)-b_p(\phi,\psi)
  =\int_{\Gamma_H}\bar{\psi}({\cal T}-{\cal T}_p)\phi ds.
\]
Analogous to \cite[Sec. 3]{chamon09}, we see immediately that
\[
  ||{\cal A}-{\cal A}_p||\leq 2||{\cal T}-{\cal T}_p||\leq \frac{2}{k\tS L}.
\]
Consequently, ${\cal A}_p$ has a bounded inverse provided that $\tS L$ is
sufficiently large as ${\cal A}$ is invertible. Since the right-hand side of
(\ref{eq:vf:pml}) defines a bounded functional in $V_H^*$, we in fact have
justified the following well-posedness result.
\begin{mytheorem}
  \label{thm:wp:pml}
  Provided that $\tS L$ is sufficiently large, the PML-truncated problem
  (\ref{eq:helm:pml}), (\ref{eq:bc1}) and (\ref{eq:pml:bc}) admits a unique
  solution $\tu^{\rm og}(x;x^*)=\tu^{\rm og}_r(x;x^*)+\chi(x;x^*)\tu^{\rm
    inc}(x;x^*)$ with $\tu^{\rm og}_r\in H_0^1(\Omega_{\rm PML})=\{\phi\in
  H^1(\Omega_{\rm PML}):\phi|_{\Gamma\cup \Gamma_{H+L}}=0\}$ for any $x^*\in
  \Omega_{\rm PML}$ such that $||\tu_r^{\rm og}(\cdot;x^*)||_{H^1(\Omega_{\rm
      PML})}\leq C||\tg^{\rm inc}||_{L^2(\Omega_{\rm PML})}$.
\end{mytheorem}
\begin{myremark}
  \label{rem:wp}
  The well-posedness in Theorem~\ref{thm:wp:pml} holds in general for any
  Lipschitz curve satisfying (GC1).
\end{myremark}
Since for any $\phi\in V_H$,
  \begin{align*}
    b_p(\phi,\psi) =b(\phi,\psi)-\int_{\Gamma_H}\bar{\psi}({\cal T}-{\cal
      T}_p)\phi ds, 
  \end{align*}
the inf-sup condition (\ref{eq:b:infsup}) of $b$ implies the inf-sup condition
of $b_p$: for any $\phi\in V_H$,
\begin{align}
  \label{eq:bf:infsup}
  \sup_{\psi\in V_H}\frac{|b_p(\phi,\psi)|}{||\psi||_{V_H}}\geq \sup_{\psi\in V_H}\frac{|b(\phi,\psi)|}{||\psi||_{V_H}} - \frac{2}{k\tS L}||\phi||_{V_H}\geq (\gamma-\frac{2}{k\tS L})||\phi||_{V_H},
\end{align}
provided $\tS L$ is sufficiently large. As a consequence of (\ref{eq:bf:infsup}), we
immediately obtain the prior error estimate for the PML truncation \tcr{if} $x^*\in
S_H$.
\begin{mycorollary}
Provided that $\tS L$ is sufficiently large, 
\begin{equation}
  \label{eq:err:est}
  ||u^{\rm og}(\cdot;x^*)-\tu^{\rm og}(\cdot;x^*)||_{V_H}\leq \frac{{2}}{\gamma k\tS L-2}||u^{\rm og}_r\tcr{(\cdot;x^*)}||_{V_H}.
\end{equation}
whenever $x^*\in S_H$.
\begin{proof}
  Since for $x^*\in S_H$, $(u^{\rm og}-\tu^{\rm og})|_{S_H}=(u^{\rm
    og}_r-\tu^{\rm og}_r)|_{S_H}\in V_H$, we have for any $\phi\in V_H$,
  \begin{align*}
    b_p(u^{\rm og}_r-\tu^{\rm og}_r,\phi)=-\int_{\Gamma_H}\bar{\phi}({\cal T}-{\cal T}_p)u^{\rm og}_rds,
  \end{align*}
  so that by the inf-sup condition (\ref{eq:bf:infsup}),
  \begin{align*}
||u^{\rm og}_r-\tu^{\rm og}_r||_{V_H}\leq& (\gamma-\frac{{2}}{k\tS L})^{-1}\sup_{\phi\in
                     V_H}\frac{|b_p(u^{\rm og}_r-\tu^{\rm og}_r,\phi)|}{||\phi||_{V_H}} =(\gamma-\frac{{2}}{k\tS L})^{-1}\sup_{\phi\in V_H}\frac{|\int_{\Gamma_H}\bar{\phi}({\cal T}-{\cal T}_p)u^{\rm og}_rds|}{||\phi||_{V_H}}\\
    \leq& \frac{{2}}{(\gamma-\frac{{2}}{k\tS L})k \tS L}||u^{\rm og}_r||_{V_H}.
  \end{align*}
\end{proof}
\end{mycorollary}

\section{Semi-waveguide problems}
\begin{figure}[!ht]
  \centering
  (a)\includegraphics[width=0.341\textwidth]{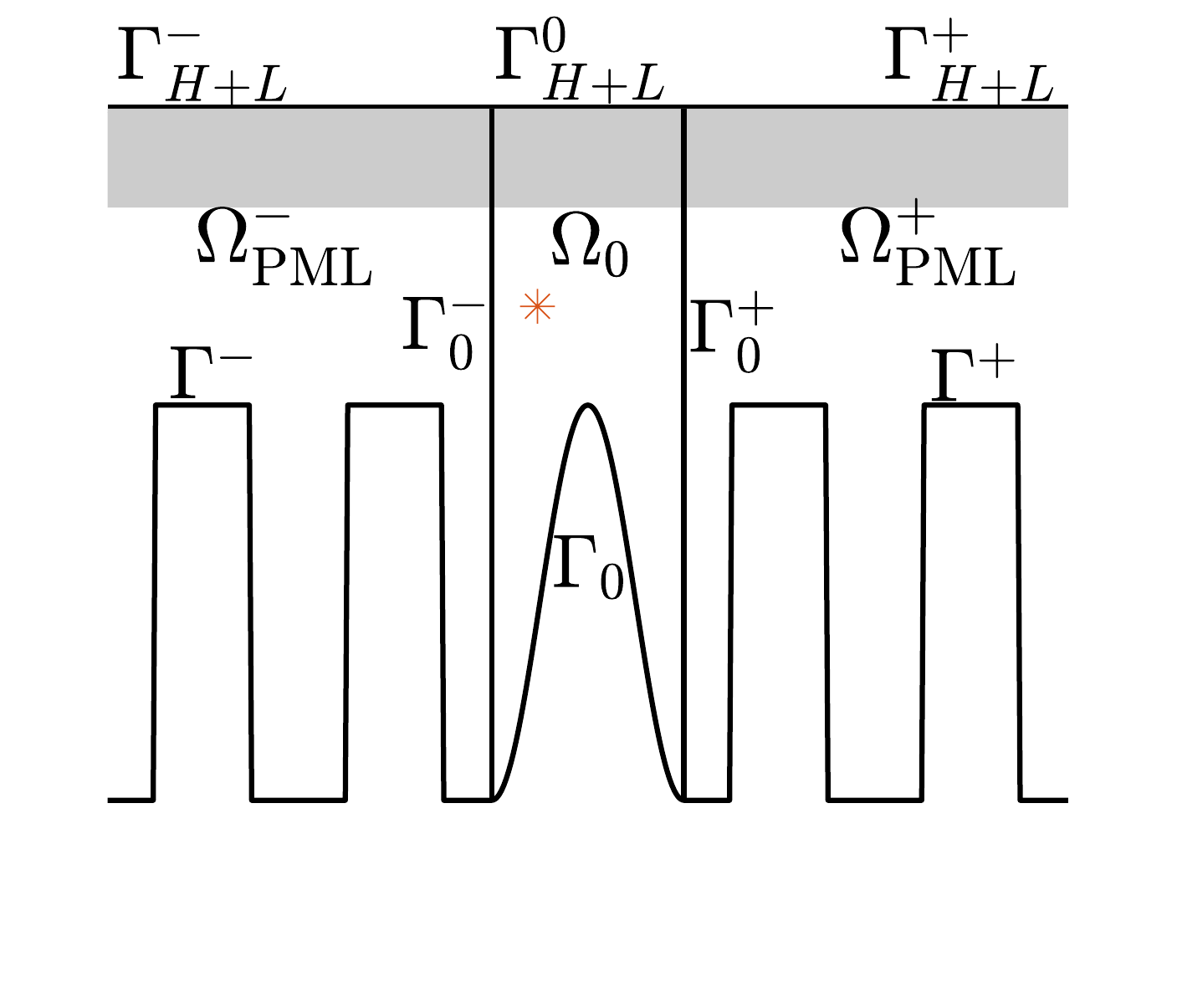}
  (b)\includegraphics[width=0.341\textwidth]{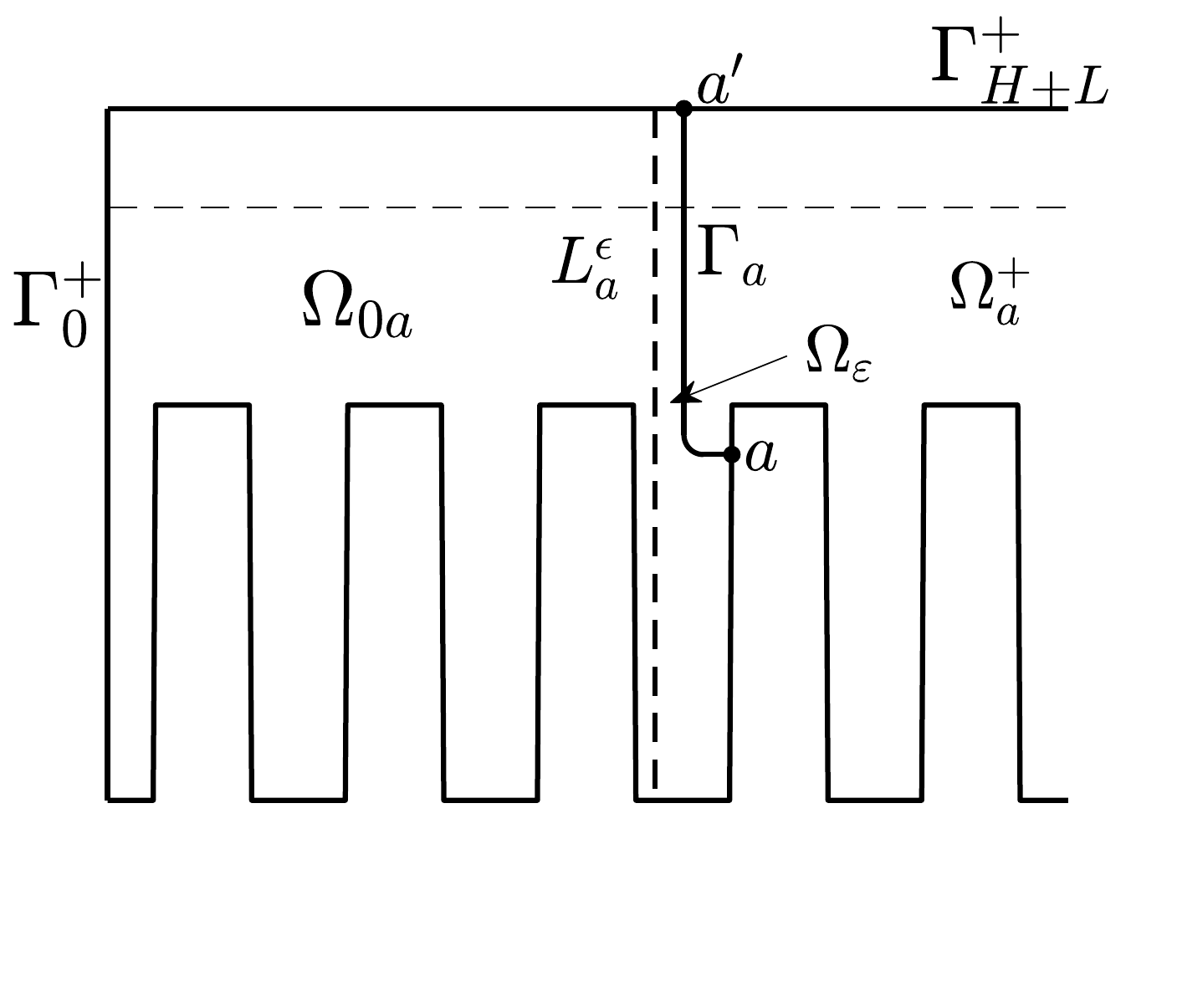}
	\caption{(a) Region $\Omega_{\rm PML}$ is divided into three regions
    $\Omega_{\rm PML}^-$, $\Omega_0$ and $\Omega_{\rm PML}^+$; semi-waveguide
    problems $(P^\pm)$ are defined in $\Omega^\pm$ bounded by
    $\Gamma_{H+L}^\pm$, $\Gamma_0^\pm$ and $\Gamma^\pm$. (b) Domain $\Omega_{\rm
      PML}^+$ is further truncated onto $\Omega_{0a}$ by a smooth curve
    $\Gamma_a$ \tcr{intersecting $\Gamma^+$ and $\Gamma_{H+L}^+$
      perpendicularly at $a$ and $a'$},
    respectively. $\Omega_a^+=\Omega_{\rm
      PML}^+\backslash\overline{\Omega_{0a}}$ and the \tcr{auxiliary line
    $L_a^\epsilon$} is chosen such that the domain
    $\Omega_\epsilon$ is sufficiently narrow. }
  \label{fig:semi:problem}
\end{figure}

Unlike the exponential convergence results in \cite{chewu03,zhowu18},
(\ref{eq:err:est}) indicates only a poor convergence of the PML method over
$S_H$. We however believe that exponential convergence can be realized in a
compact subset of $S_H$, which is indeed true \tcr{if} $\Gamma$ is flat
\cite{chamon09}. This leads to an essential question after the vertical PML
truncation: how to accurately truncate $\Omega_{\rm PML}$ in the lateral
$x_1$-direction? To address this question, as inspired by \cite{jollifli06} and
as illustrated in Figure~\ref{fig:semi:problem} (a), we shall consider the
following two semi-waveguide problems:
\begin{equation*}
  (P^{\pm}):\quad \left\{
    \begin{array}{l}
        \nabla\cdot({\bf A} \nabla \tilde{u}) + k^2\alpha\tilde{u} = 0, \quad{\rm on}\quad{\Omega}_{\rm PML}^{\pm}:=\Omega_{\rm PML}\cap\left\{x: \pm x_1>\frac{T}{2} \right\},\\
  \tilde{u} = 0,\quad {\rm on}\quad \Gamma^{\pm}:=\Gamma\cap\left\{x: \pm x_1>\frac{T}{2} \right\},\\
  \tu=0,\quad{\rm on}\quad  \Gamma^{\pm}_{L+H}:=\Gamma_{L+H}\cap\left\{x: \pm x_1>\frac{T}{2} \right\},\\
  \partial_{\nu_c}\tu=g^{\pm},\quad{\rm on}\quad \Gamma_0^{\pm}:=\Omega_{\rm PML}\cap \left\{x:x_1=\pm \frac{T}{2}\right\},\\
\end{array}
\right.
\end{equation*}
for given Neumann data $g^{\pm}\in H^{-1/2}(\Gamma_0^{\pm})$, where $\nu_c={\bf
  A}\nu$ denotes the co-normal vector with $\nu$ pointing towards $\Omega_{\rm PML}^\pm$, $\tu$ denotes a generic field, and we note that $\Gamma^{\pm}\subset
\Gamma_T$ does not contain the defected part $\Gamma_0$.
In this section, we shall study the well-posedness of the semi-waveguide
problems ($P^{\pm}$).

By Theorem~\ref{thm:wp:pml}, \tcr{the following uniqueness result is easy to obtain}.
\begin{mylemma}
  \label{lem:uniq:swp}
  Provided that $\tS L$ is sufficiently large, problem~($P^{\pm}$) has at most one
  solution in $H^1(\Omega_{\rm PML}^{\pm})$.
  \begin{proof}
    Suppose $\tu\in H^1(\Omega_{\rm PML}^{+})$ \tcr{satisfies} ($P^+$) with $g^+=0$.
    Let
    \begin{align*}
      \Omega_{\rm PML}^e =& \{x\in\mathbb{R}^2:(x_1,x_2)\in \Omega_{\rm PML}^{+}\ {\rm or}\ (T-x_1,x_2)\in \Omega_{\rm PML}^{+}\}\cup \Gamma_0^+,\\
      \Gamma^e =& \{x\in\mathbb{R}^2:(x_1,x_2)\in \Gamma^{+}\ {\rm or}\ (T-x_1,x_2)\in \Gamma^{+}\ {\rm or}\ (T/2,x_2)\in\Gamma\}.
    \end{align*}
    Then, 
    \[
      \tu^e(x_1,x_2) = \left\{
        \begin{array}{ll}
        \tu(x_1,x_2),& x_1\geq T/2,\\
        \tu(T-x_1,x_2),& x_1<T/2,\\
        \end{array}
      \right.
    \]
    in $H^1(\Omega_{\rm PML}^e)$ \tcr{satisfies} problem~(\ref{eq:pm:u1}-\ref{eq:pm:u3}) 
    with $\tg$, $\Omega_{\rm PML}$ and $\Gamma$ replaced by $0$, $\Omega_{\rm
      PML}^e$ and $\Gamma^e$, respectively. Theorem~\ref{thm:wp:pml} and
    Remark~\ref{rem:wp} 
    imply that $\tu^e= 0$ on $\Omega_{\rm PML}^e$ so that
    $\tu=\tu^e|_{\Omega_{\rm PML}^+}=0$. The uniqueness of problem ($P^-$) can
    be established similarly.
  \end{proof}
\end{mylemma}
We are ready to study the well-posedness of problem ($P^{\pm}$) by the Fredholm
alternative. Without loss of generality, we shall study ($P^+$) only. To make
use of Fredholm theory, we need first to truncate $\Omega_{\rm PML}^+$ by an
exact transparent boundary condition. Under condition (GC2), there exists a line
segment $L_a\subset \Gamma_T\cap\Gamma$ with the midpoint $a=(a_1,a_2)\in L_a$
for $a_1>T/2$. For \tcr{a small fixed constant $\epsilon>0$}, we can find a
vertical line segment $L_a^\epsilon$ and a simple and smooth curve
$\Gamma_{a}\subset\Omega_{\rm PML}$ connecting $L_a$ and $\Gamma_{H+L}^+$ such
that the distance of $L_a^\epsilon$ and $\Gamma_a$ is $\epsilon$ and that
$\Gamma_a$ \tcr{intersecting $L_a$ and $\Gamma_{H+L}^+$ perpendicularly at $a$
  and $a'$}, respectively, as shown in Figure~\ref{fig:semi:problem} (b). Let
$\Omega_\epsilon$ be the domain bounded by $\Gamma_a$, $L_a^\epsilon$, $L_a$ and
$\Gamma_{H+L}^+$ and $\Omega_a^+$ be the unbounded domain bounded by $\Gamma_a$,
$\Gamma_{H+L}^+$ and $\Gamma$. For sufficiently small $\epsilon$, the above
choice of $L_a^\epsilon$ and $\Gamma_a$ guarantees that $k>0$ is not an
eigenvalue of
\begin{align}
  \label{eq:pro:eps}
  -\nabla\cdot({\bf A}\nabla \tu) =& k^2\alpha \tu,\quad{\rm on}\quad\Omega_\epsilon,\\
  \label{eq:bc:eps}
  \tu =& 0,\quad{\rm on}\quad\partial\Omega_\epsilon.
\end{align}
Now for the unbounded domain $\Omega_\epsilon^+=\Omega_\epsilon\cup\Gamma_a\cup
\Omega_a^+$, by a \tcr{symmetrical reflection w.r.t the line containing
$L_a^\epsilon$}, the partial boundary $\partial\Omega_\epsilon^+\cap \Gamma$ can
be extended to a Lipschitz boundary, denoted by $\Gamma_\epsilon$, satisfying
(GC1). Then, Theorem~\ref{thm:wp:pml}, with $\Gamma_\epsilon$ in place of
$\Gamma$, can help to construct the Dirichlet Green function of
$\Omega_\epsilon^+$ by $\tG_{\rm D}(x;y)=\tu^{\rm og}(x;y)-\tu^{\rm og}(x;y_{\rm
  imag}^\epsilon)$ satisfying \tcr{$\tG_{\rm D}(\cdot;y)|_{\partial\Omega^+_\epsilon}=0$}, where
$y^\epsilon_{\rm imag}$ is \tcr{the mirror image of the source point $y$ w.r.t line
$L_a^\epsilon$.} Choosing $\Gamma_{a}$ in such a special way\tcr{, the following local regularity property of $\tG_{\rm D}(x;y)$ can be ensured.} 
\begin{myprop}
  \label{prop:lbg}
  Under \tcr{the }geometrical conditions (GC1) and (GC2), for sufficiently large values of $L$ and $m$ in
  (\ref{eq:sigma}), $\tG_{\rm D}(x;y)$ admits the following decomposition
  \begin{equation}
    \label{eq:dec:pmlG}
    \tG_{\rm D}(x;y) = \tG(x;y) - \tG(x;y^l_{\rm imag}) + R_l(x,y),y\in \overline{\Omega_l},
  \end{equation}
  such that $R_l(x,y)$ is a sufficiently smooth function of $x$ and $y$ for
  $(x,y)\in\overline{\Omega_{l}\cup \Omega_{c}}\times\overline{\Omega_l}$, where \tcr{$\tilde{G}$ is defined by (\ref{eq:fun:helm2}), }
  $\Omega_l$ is a sufficiently small neighborhood of point $l$ in
  $\Omega_{\epsilon}^+$ and $y^l_{\rm imag}$ is \tcr{the mirror image of $l$ w.r.t
  line $L_l$} for
  $l=a,a'$, and $\Omega_c$ can be any bounded subset of $\Omega_{\epsilon}^+$.
  \begin{proof}
    We consider $y$ close to point $a'$ only. Define
    \[
      \tcr{\tu_{a'}(x;y) := \tu^{\rm og}(x;y) - \chi_{a'}(x)\left[  \tG(x;y) - \tG(x;y^{a'}_{\rm imag})\right]},
    \]
    where the cut-off function $\chi_{a'}=1$ in a neighborhood of \tcr{$a'$} and
    has a small support that is independent of $y$. Then, it can be
    seen that $\tu_{a'}$ \tcr{satisfies} (\ref{eq:pm:u1}-\ref{eq:pm:u3}) with
    \tcr{$\tg^{\rm inc}$} replaced
    by 
    \[
      [\nabla\cdot({\bf A} \nabla ) + k^2\alpha](1- \chi_{a'}(x))\left[
        \tG(x;y) - \tG(x;y^{a'}_{\rm imag})\right]\in\tcr{ C^{m-1}_{\rm
        comp}}(\Omega_{\rm PML}\times \overline{\Omega_l}),
    \]
    where \tcr{$C^{m-1}_{\rm comp}$ consists of $m-1$ times
      differentiable functions with compact supports}, and we note that $m$
    defined in (\ref{eq:sigma}) determines the smoothness of $\sigma$. By
    arguing the same way as in \cite[Lem 2.4]{hulurat20} and by choosing $m$
    sufficiently large, $R_{a'}(x,y)=\tu_{a'}(x;y)-\tu^{\rm og}(x;y_{\rm
      imag}^\epsilon)$ becomes a sufficiently smooth function for
    $(x,y)\in\overline{\Omega_{l}\cup\Omega_c}\times\overline{\Omega_l}$.
  \end{proof}
\end{myprop}
On $\Gamma_{a}$, we now define the following two integral operators:
\begin{align}
  \label{eq:S}
  [{\cal S}_{a} \phi](x) &= 2\int_{\Gamma_{a}} \tG_{\rm D}(x;y) \phi(y) ds(y),\\
  [{\cal K}_{a} \phi](x) &= \tcr{2\int_{\Gamma_{a}} \partial_{\nu_c(y)}\tG_{\rm D}(x;y) \phi(y) ds(y)},
\end{align}
Proposition~\ref{prop:lbg} reveals that classic mapping properties hold for the
above two integral operators on the open arc $\Gamma_{a}$.
\begin{mylemma}
  \label{lem:ext:SK} 
  We can uniquely extend the operator $\cS_{a}$ as a bounded operator from
  $H^{-1/2}(\Gamma_{a})$ to $\wtd{H^{1/2}}(\Gamma_{a})$, the operator
  $\cK_{a}$ as a compact (and certainly bounded) operator from
  $\wtd{H^{1/2}}(\Gamma_{a})$ to $\wtd{H^{1/2}}(\Gamma_{a})$.
  Moreover, we
  have the decomposition $\cS_{a} = \cS_{p,a} + {\cal L}_{p,a}$ such that
  $\cS_{p,a}: H^{-1/2}(\Gamma_{a})\to \wtd{H^{1/2}}(\Gamma_{a})$ is positive and
  bounded below, i.e., for some constant $c>0$,
  \[
    {\rm Re}\left( \int_{\Gamma_a}{\cal S}_{p,a}\phi\bar{\phi} ds
    \right)\geq c ||\phi||^2_{H^{-1/2}(\Gamma_{a})},
  \]
  for any $\phi\in H^{-1/2}(\Gamma_{a})$, and ${\cal
    L}_{p,a}:H^{-1/2}(\Gamma_{a})\to \wtd{H^{1/2}}(\Gamma_{a})$ is compact.
  \begin{proof}
    By Proposition~\ref{prop:lbg}, the proof follows from similar arguments as
    in \cite[Sec. 2.3]{hulurat20} \tcr{but relies on Fredholm of the
      single-layer potential and compactness of the double-layer potential of
      kernels relating to $\tG$, the fundamental solution of the strongly
      elliptic Helmholtz equation (\ref{eq:helm:pml}), as has been studied in
      \cite[Thm. 7.6]{mcl00} and \cite{lassom01}.} We omit the details.
  \end{proof}
\end{mylemma}
\tcr{Analogous to \cite[Lem. 5.1]{lassom01}, one gets the following the Green's
representation}
\begin{equation}
  \label{eq:rep}
{\tu(x) = \int_{\Gamma_a}\left[\partial_{\nu_c(y)}\tG_{\rm D}(x;y) \tu(y)-\tG_{\rm D}(x;y)\partial_{\nu_c(y)}\tu(y)  \right]ds(y).}
\end{equation}
\tcr{By the jump relations \cite[Thm. 5.1]{lassom01},} letting $x$ approach $\Gamma_a$, we get the following transparent
boundary condition (TBC) 
\begin{equation}
  \label{eq:tbc}
  {\tu-{\cal K}_a\tu=-{\cal S}_a\partial_{\nu_c}\tu,\quad {\rm on}\quad\Gamma_a.}
\end{equation}
As indicated in Figure~\ref{fig:semi:problem} (b), let $\Omega_{0a}$ be the domain bounded by
$\Gamma_0^+$, $\Gamma_{a}$, $\Gamma_{H+L}^+$ and $\Gamma_+$,
\[
H_D^1(\Omega_{0a})=\{v|_{H^1(\Omega_{0a})}: v\in H^1(\Omega_{\rm PML}^+),
v|_{\Gamma_+}=0, v|_{\Gamma_{H+L}}=0\},
\]
and $V_{a}=H_D^1(\Omega_{0a})\times H^{-1/2}(\Gamma_{a})$ be equipped with the
natural cross-product norm. ($P^+$) can be equivalently formulated as the
following boundary value problem: find $(\tu,\phi)\in V_{a}$ solving
\begin{align}
  &\nabla\cdot({\bf A} \nabla \tilde{u}) + k^2\alpha\tilde{u} = 0, \quad{\rm on}\quad{\Omega}_{0a},\\
  &\partial_{\nu_c}\tu|_{\Gamma_0^+}=g^{+},\quad{\rm on}\quad\Gamma_0^+,\\
  &\partial_{\nu_c}\tu|_{\Gamma_a}=\phi,\quad{\rm on}\quad\Gamma_a,\\
  &\tu-{\cal K}_a\tu=-{\cal S}_a\phi,\quad {\rm on}\quad\Gamma_a.
\end{align}
An equivalent variational formulation reads: find $(\tu,\phi)\in V_a$ such that
\begin{equation}
  \label{eq:vp:bps}
  b_{ps}((\tu,\phi),(v,\psi)) = \int_{\Gamma_0^+}g^+\bar{v}ds,
\end{equation}
for all $(v,\psi)\in V_a$, where the sesquilinear form $b_{ps}(\cdot,\cdot):
V_a\times V_a\to\mathbb{C}$ is given by
\tcr{
\begin{align*}
  b_{ps}((\tu,\phi),(v,\psi)) =& \int_{\Omega_{0a}}\left[({\bf A}\nabla \tu)^{T}\overline{\nabla v} - k^2\alpha \tu\bar{v}\right]dx - \int_{\Gamma_a}\left[\phi\bar{v} -\left(\tu-{\cal K}_a\tu + {\cal S}_a\phi\right)\bar{\psi}\right] ds.
\end{align*}
}
We are now ready to establish the well-posedness of problems (${P}^{+}$).
\begin{mytheorem}
  \label{thm:sw:wp}
  Under \tcr{the }geometrical conditions (GC1) and (GC2), provided that $L$ is sufficiently large, the
  semi-waveguide problem ($P^{\pm}$) has a unique solution $\tu\in
  H^{1}(\Omega_{\rm PML}^{\pm})$ such that $||\tu||_{H^1(\Omega_{\rm PML}^\pm)}\leq C||g^\pm||_{H^{-1/2}(\Gamma_0^\pm)}$ for any $g^\pm\in
  H^{-1/2}(\Gamma_0^{\pm})$, respectively, where $C$ is independent of $g^\pm$.
  \begin{proof}
    We study ($P^+$) only. For the variational problem (\ref{eq:vp:bps}), we can
    decompose $b_{ps}=b_1+b_2$ where
    \begin{align*}
      b_1((\tu,\phi),(v,\psi)) =& \int_{\Omega_{0a}}\left[({\bf A}\nabla \tu)^{T}\overline{\nabla v}- k^2\alpha\tu\bar{v}\right]dx -\int_{\Gamma_a}\left[\phi\bar{v} -\tu\bar{\psi} -{\cal S}_{p,a}\phi\bar{\psi}  \right] ds,\\ 
      b_2((\tu,\phi),(v,\psi)) =&  \int_{\Gamma_a}\left[{\cal L}_{p,a} \phi -\cK_a\phi\right] \bar{\psi}   ds.
    \end{align*}
    According to Lemma~\ref{lem:ext:SK}, $b_1$ is coercive on $V$ as
    \begin{align*}
      {\rm Re}(b_1((\tu,\phi),(\tu,\phi)))=& \int_{\Omega_{0a}}\left[|\tu_{x_1}|^2+ (1+\sigma^2(x_2))^{-1}|\tu_{x_2}|^2 - k^2|\tu|^2   \right]dx+{\rm Re}\left( \int_{\Gamma_a}{\cal S}_{p,a}\phi\bar{\phi} ds \right)\\
      \geq& c||\tu||_{H^1(\Omega_{0a})}^2 - C||\tcr{\tilde{u}}||_{L^2(\Omega_{0a})}^2 + c||\phi||_{H^{-1/2}(\Gamma_a)}^2,
    \end{align*}
    and the bounded linear operator associated with $b_2$ is compact.
    Consequently, $b_{ps}$ is Fredholm of index zero \cite[Thm. 2.34]{mcl00}. 

    Now, we prove $\tu=0$ and $\phi=0$ \tcr{when }
     $g^+=0$. By (\ref{eq:rep}), we can
    directly extend $\tu$ to $\Omega_{\epsilon}^+$, denoted by $\tu^{\rm ext}$.
    Then, the TBC (\ref{eq:tbc}) implies $\gamma^+\tu^{\rm
      ext}|_{\Omega_a^+}=\tu|_{\Gamma_a}$ so that by the \tcr{jump} relations,
    \tcr{$\gamma^-\tu^{\rm ext}|_{\Omega_\epsilon} = 0$ where $\gamma^{+}$
    ($\gamma^-$) }defines the trace operator of $\tu^{\rm ext}$ onto $\Gamma_a$
    from $\Omega_a^+$ ($\Omega_\epsilon$). Thus, $\tu^-=\tu^{\rm
      ext}|_{\Omega_\epsilon}$ \tcr{satisfies} (\ref{eq:pro:eps}) and (\ref{eq:bc:eps}).
    But the special choice of $\epsilon$ and $\Omega_\epsilon$ has ensured that
    $\tu^-\equiv 0$ on $\Omega_\epsilon$ so that \tcr{the trace of }$\partial_{\nu_c}\tu^{\rm ext}$
    taken from $\Omega_\epsilon$ is $0$.
    The \tcr{jump} conditions then imply that \tcr{the trace of } $\partial_{\nu_c}\tu^{\rm ext}$ taken
    from $\Omega_a^+$ is $\phi$. Consequently,
    \begin{align*}
      w(x) = \left\{
      \begin{array}{ll}
        \tu(x),& x\in\Omega_{0a},\\
        \tu^{\rm ext}(x),& x\in\Omega_{a}^+,
      \end{array}
      \right.
      \end{align*}
      belongs to $H^1(\Omega_{\rm PML}^+)$ and \tcr{satisfies} ($P^+$) with $g^+=0$. But
      Lemma~\ref{lem:uniq:swp} already justifies that $w$ must be $0$ on
      $\Omega_{\rm PML}^+$, which indicates that $\tu=0$ and $\phi=0$. The proof
      then follows from the fact that the right-hand side of (\ref{eq:vp:bps})
      defines a bounded anti-linear functional in $V_a^*$.
  \end{proof}
\end{mytheorem}
\begin{myremark}
  Like Therorem~\ref{thm:wp:pml},
  Theorem~\ref{thm:sw:wp} also holds for any Lipschitz curves
  $\Gamma^\pm$, \tcr{which are not
  necessarily periodic}, satisfying \tcr{the }geometrical conditions (GC1) and (GC2).
\end{myremark}
\section{Lateral boundary conditions}
According to \tcr{Theorem}~\ref{thm:wp:pml}, $\partial_{\nu_c}\tu^{\rm
  og}(\cdot;x^*)|_{\Gamma_0}^{\pm}\in H^{-1/2}(\Gamma_0^{\pm})$ for any $x^*\in
S_H$ with $|x_1^*|<T/2$. Thus, $\tu=\tu^{\rm og}(\cdot;x^*)|_{\Omega_{\rm
    PML}^{\pm}}$ \tcr{satisfies} ($P^{\pm}$) with $g^\pm=\partial_{\nu_c}\tu^{\rm
  og}(\cdot;x^*)|_{\Gamma_0}^{\pm}$ in the distributional sense, respectively.
Theorem~\ref{thm:sw:wp} then implies that we can define two vertical
Neumann-to-Dirichlet (vNtD) operators ${\cal N}^{\pm}:
H^{-1/2}(\Gamma_0^{\pm})\to \widetilde{H^{1/2}}(\Gamma_0^{\pm})$ satisfying
$\tu^{\rm og}|_{\Gamma_0^{\pm}}={\cal N}^{\pm}\partial_{\nu_c}\tu^{\rm
  og}|_{\Gamma_0^{\pm}}$. Such transparent boundary conditions can serve as
exact {\it lateral boundary conditions} to terminate the $x_1$-variable for the
PML-truncated problem (\ref{eq:helm:pml}) and (\ref{eq:bc1}). Consequently, the
original unbounded problem (\ref{eq:gov:utot}) and (\ref{eq:gov:bc:te}) equipped
with the hSRC condition (\ref{eq:src}) can be truncated onto the perturbed cell
$\Omega_{0}:=\Omega_{\rm PML}\cap\left\{x: |x_1|<\frac{T}{2} \right\}$ and be
reformulated as the following boundary value problem:
\begin{equation*}
  ({\rm BVP1}): \left\{
    \begin{array}{ll}
      \nabla\cdot({\bf A} \nabla \tilde{u}^{\rm og}) + k^2\alpha\tilde{u}^{\rm og} = -\delta(x-x^*), &{\rm on}\ \Omega_{0},\\
      \tilde{u}^{\rm og} = 0,&{\rm on}\ \Gamma_{0}=\Gamma\cap\{x:|x_1|<T/2\} ,\\
      \tu^{\rm og}=0,&{\rm on}\ \Gamma^{0}_{H+L}=\Gamma_{H+L}\cap\{x:|x_1|<T/2\},\\
      \tu^{\rm og}={\cal N}^{\pm}\partial_{\nu_c}\tu^{\rm og},&{\rm on}\ \Gamma_0^{\pm}.\\
    \end{array}
       \right.
\end{equation*}
Theorems~\ref{thm:wp:pml} and \ref{thm:sw:wp} directly imply that
\tcr{(BVP1) admits the following unique solution
  \[
    \tu^{\rm og}(\cdot;x^*)=\tu^{\rm og}_r(\cdot;x^*)|_{\Omega_0}+\chi(\cdot;x^*)|_{\Omega_0}\tu^{\rm inc}(x;x^*)|_{\Omega_0},
\]
with $\tu^{\rm og}_r$ defined in Theorem~\ref{thm:wp:pml}. }Nevertheless, it is
challenging to get ${\cal N}^\pm$ by directly solving the unbounded problem
($P^{\pm}$) in practice. To overcome this difficulty, in this section, we shall
define two closely related Neumann-marching operators, derive the governing
\tcr{Riccati} equations, and design an efficient RDP to accurately approximate
${\cal N}^\pm$.

\subsection{Neumann-marching operators ${\cal R}_p^{\pm}$}
\begin{figure}[!ht]
  \centering
(a)\includegraphics[width=0.35\textwidth]{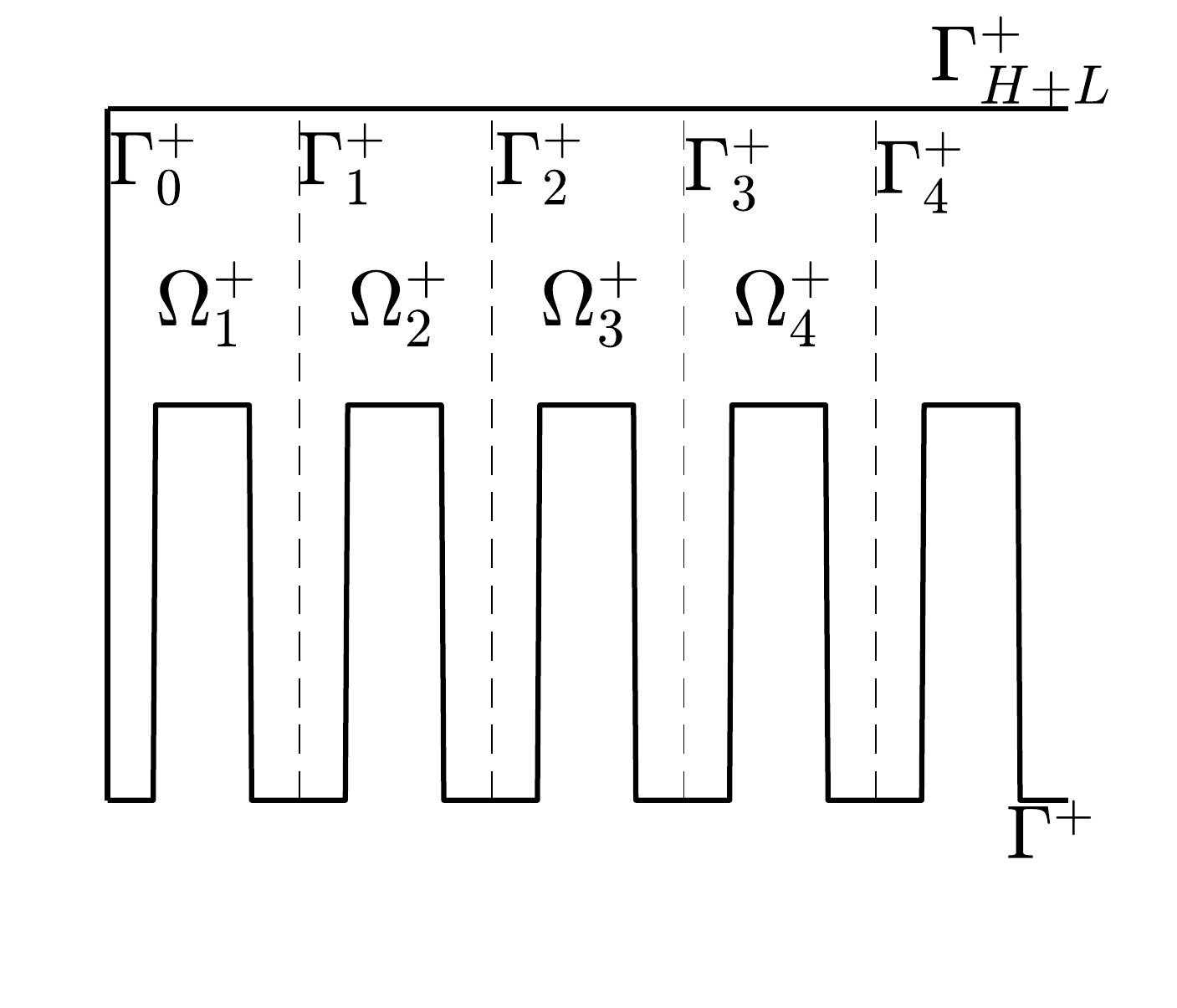}
(b)\includegraphics[width=0.35\textwidth]{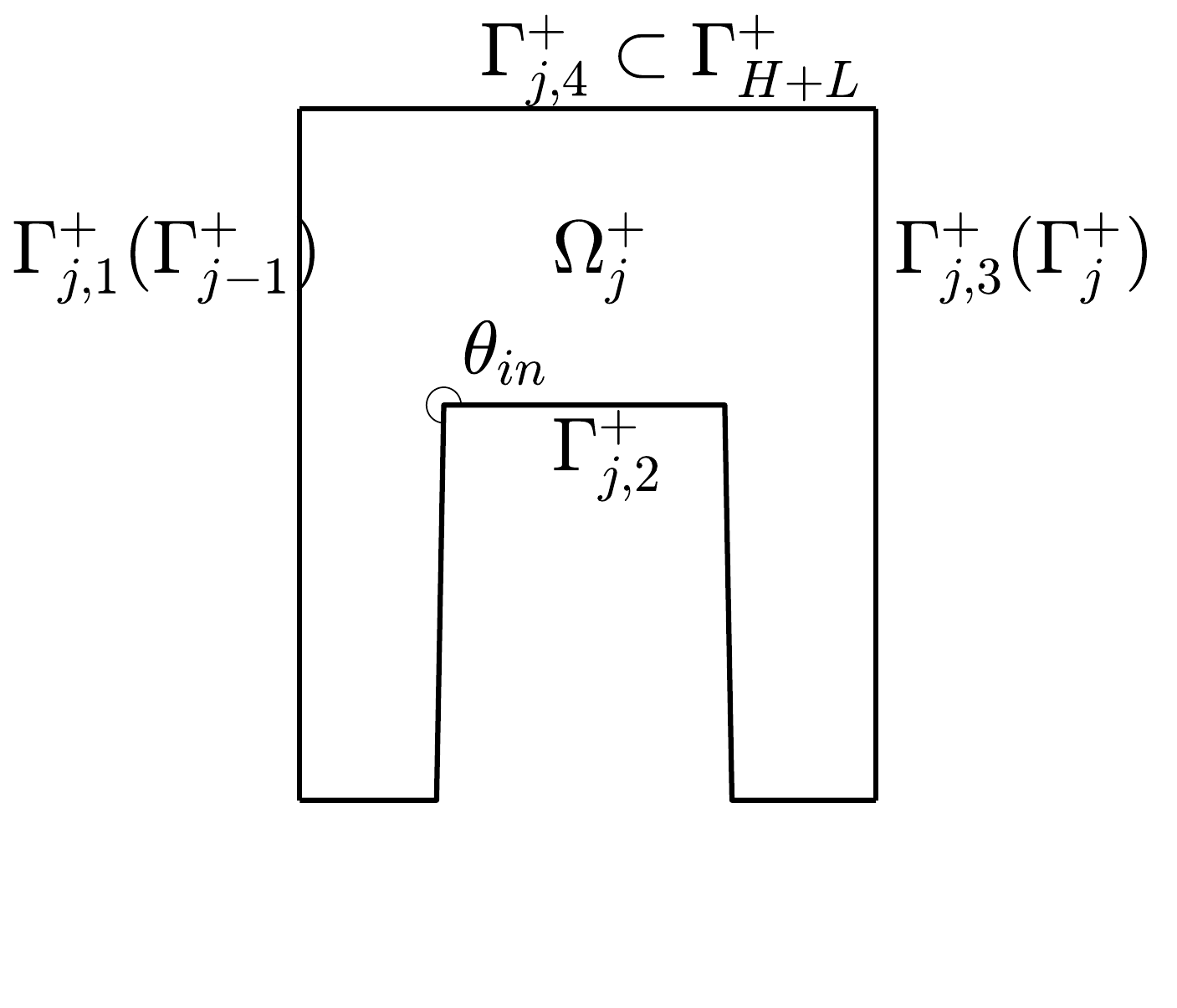}
\caption{(a) The semi-waveguide region $\Omega_{\rm PML}^+$ is divided into
  $\{\Omega_{j}^+\}_{j=1}^{\infty}$ of the same shape. The operator ${\cal
    R}_p^+$ can then march Neumann data through the vertical line segments
  $\{\Gamma_j^+\}_{j=0}^{\infty}$. (b) The boundary of $\Omega_j$ consists of
  four parts: $\Gamma_{j,1}^+$ (left), $\Gamma_{j,2}^+$ (bottom), $\Gamma_{j,3}^+$
  (right), $\Gamma_{j,4}^+$ (top). Here, $\theta^{\rm in}$ indicates the interior
  angle at a corner, as will be used in (\ref{eq:k01:exa}).}
  \label{fig:lateral}
\end{figure}

Now, let
\begin{align*}
\Gamma_j^\pm=&\{(x_1\pm jT,x_2):, x=(x_1,x_2)\in
\Gamma_0^\pm\},\\
\Omega_{{\rm PML},j}^\pm=&\{x\in\Omega_{\rm PML}^\pm:\pm x_1>T/2+(j-1)T\},\\
\tcr{\Omega_{j}^\pm}=&\Omega_{{\rm PML},j}^\pm\backslash\overline{\Omega^{\pm}_{{\rm PML},j+1}}
\end{align*}
for $j\in\mathbb{N}^*$, as illustrated in Figure~\ref{fig:lateral}(a) for \tcr{the}
notations of superscript $+$.

As insipired by \cite{jollifli06}, the well-posedness of
($P^{\pm}$) well defines two bounded Neumann-marching operators ${\cal
  R}_p^{\pm}:H^{-1/2}(\Gamma_0^{\pm})\to H^{-1/2}(\Gamma_1^{\pm})$ such that
$\partial_{\nu_c^\pm}\tu^{\rm og}|_{\Gamma_1^{\pm}}={\cal
  R}_p^{\pm}\partial_{\nu_c^\pm} \tu^{\rm og}|_{\Gamma_0^{\pm}}$, where
$\nu_c^\pm={\bf A}\nu^\pm$ with $\nu^\pm=(\pm 1,0)^{T}$. We have the following
properties of ${\cal R}_p^{\pm}$, analogous to \cite[Thm. 3.1]{jollifli06}.
\begin{myprop}
  Under the conditions that (GC2) holds and $kL$ is sufficiently large, we can
  choose $\Gamma_0^{\pm}$ intersecting $\Gamma$ at a smooth point
  such that ${\cal R}_p^{\pm}$ are compact operators and 
  \begin{equation}
    \label{eq:R+-}
    \partial_{\nu_c^\pm}\tu^{\rm og}|_{\Gamma_{j+1}^{\pm}}={\cal R}_p^{\pm}\partial_{\nu_c^\pm} \tu^{\rm og}|_{\Gamma_j^{\pm}},
  \end{equation}
  holds for any $j\geq 0$. Furthermore, 
  \begin{equation}
    \label{eq:radius}
    \rho({\cal R}_p^\pm)<1,
  \end{equation}
  where $\rho$ denotes the spectral radius.
  \begin{proof}
    We study only the property of ${\cal R}_p^+$. The choice of
    $\Gamma_{0}^{\tcr{+}}$ and the interior regularity theory of elliptic operators
    directly imply the compactness of ${\cal R}_p^+$.

    It is clear that (\ref{eq:R+-}) holds for $j=0$. We need only justify the
    case $j=1$ as all others can be done by induction. Consider the
    semi-waveguide problem ($P^+$) with $g^+=-\partial_{\nu_c^+} \tu^{\rm
      og}|_{\Gamma_1^{+}}$, where the negative sign appears since $\nu_c^+=-\nu_c$. Theorem~\ref{thm:sw:wp} implies that $\tu_n^{\rm
      og}(x)=\tu^{\rm og}(x_1+T,x_2)$ for $x\in \Omega_{\rm PML}^+$ is the unique
    solution. Then $\partial_{\nu_c^+}\tu^{\rm og}_n|_{\Gamma_1^+}={\cal R}_p^+
    \partial_{\nu_c^+}\tu^{\rm og}_n|_{\Gamma_0^+}$, which reads exactly
    $\partial_{\nu_c^+}\tu^{\rm og}|_{\Gamma_{2}^{+}}={\cal R}_p^{+}\partial_{\nu_c^+}
    \tu^{\rm og}|_{\Gamma_1^{+}}$.

    Now we prove (\ref{eq:radius}) by contradiction. Suppose otherwise there
    exists $0\neq g\in H^{-1/2}(\Gamma^+_0)$ such that ${\cal R}_p^+g=\lambda_0g$ with
    $|\lambda_0|\geq 1$. Suppose $\tu$ \tcr{satisfies} ($P^+$) with $g^+=g$ on
    $\Gamma_0^+$. Then, for any \tcr{$v\in H^1(\Omega^+_{\rm PML})$
    so that $v(\cdot-jT,\cdot)\in H^1(\Omega^+_{{\rm PML},j+1})$} for any $j\geq 0$, we have
    by Green's identity that,
    \begin{align*}
      |\lambda_0|^j\left|\tcr{\int_{\Gamma_0}g\bar{v}ds}\right|=& \left|\tcr{\int_{\Gamma_j}({\cal R}_p^+)^jg \overline{v(\cdot-jT,\cdot)}ds}\right| \\
      =& \left|  \int_{\Omega_{{\rm PML},j+1}}\left[ ({\bf A}\nabla \tu^{\rm og})^{T} \overline{\nabla v(\cdot-jT,\cdot)}- k^2\alpha \tu^{\rm og}\overline{v(\cdot-jT,\cdot))} \right]dx\right|\\
      \leq & C ||\tu^{\rm og}||_{H^1(\Omega^+_{{\rm PML}, j+1})}||v||_{H^1(\Omega^+_{\rm PML})}\to 0,\quad j\to\infty,
    \end{align*}
    which is impossible.
  \end{proof}
\end{myprop}
By the following identity \cite{kato95},
\[
  \rho({\cal R}_p^\pm) = \lim_{j\to\infty} ||({\cal R}_p^\pm)^j||^{1/j},
\]
it can be seen that there exists a sufficiently large integer $N_0>0$ such that
$({\cal R}_p^\pm)^{N_0}$ is contracting, i.e.,
\begin{equation}
  \label{eq:cont:Rp}
||({\cal R}_p^\pm)^{N_0}||<1.
\end{equation}
Let $\Omega_j^{\pm,N_0}$ be the interior of $N_0$ consecutive cells
$\cup_{j'=1}^{N_0}\overline{\Omega^{\pm}_{(j-1)N_0+j'}}$. As a corollary, the
above results indicate that $\tu^{\rm og}$ decays exponentially at infinity of
the strip.
\begin{mycorollary}
  \label{cor:exp:conv}
  Under the conditions that (GC2) holds and $kL$ is sufficiently large, 
  \begin{equation}
    \label{eq:exp:conv}
    ||\tu^{\rm og}(\cdot;x^*)||_{H^1(\Omega_{j}^{\pm,N_0})}\leq C||({\cal
      R}_p^{\pm})^{N_0}||^{j-1}||\tg^{\rm inc}||_{L^2(\Omega_{\rm PML})},
  \end{equation}
  where we recall that $\tg^{\rm inc} = [\nabla\cdot({\bf A} \nabla ) +
  k^2\alpha] (1-\chi(x;x^*)) \tu^{\rm inc}(x;x^*)$, and $C$ is independent of
  $j\geq 0$. In other words, the PML truncated solution $\tu^{\rm og}(x;x^*)$
  decays exponentially fast to $0$ in the strip as $|x_1|\to\infty$ for any
  $x^*\in\Omega_{\rm PML}$.
\end{mycorollary}
\begin{myremark}
  Authors in \cite{chamon09} have revealed a similar result as
  (\ref{eq:exp:conv}) for $\Gamma$ being a flat surface. The above corollary
  indicates that such an exponentially decaying property for the PML truncated
  solution holds even for locally defected periodic curves. As a consequence,
  this reveals that the PML truncation cannot realize an exponential convergence
  to the true solution for numerical solutions at regions sufficiently away from
  the source or local defects since the true solution is expected to decay only
  of an algebraic rate at infinity: \cite{chamon05} has indicated that $u^{\rm
    og}$ behaves as ${\cal O}(x_1^{-3/2})$ as $x_1\to\infty$.
\end{myremark}

Though Corollary~\ref{cor:exp:conv} provides hopeless results, we point out
that (\ref{eq:exp:conv}) holds for $L$ being fixed but $j\to\infty$. If, on the
contrary, $j$ is fixed but $L\to\infty$, we believe exponential convergence can
still be achieved. In doing so, we need a more \tcr{effective} description \tcr{of} the
Neumann-marching operators ${\cal R}^{\pm}_p$, as was done in \cite{jollifli06}.
Take ${\cal R}^+_p$ as an example. As shown in Figure~\ref{fig:lateral}(b), recall that
$\Omega_j^+$ denotes the $j$-th unit cell on the right of $\Gamma_0^+$, which is
unperturbed for $j\geq 1$, and to simplify the presentation, we further denote
the four boundaries of $\Omega_j^+$ by
\begin{align*}
  \Gamma_{j,1}=\Gamma_{j-1}^+,\quad \Gamma_{j,3} =\Gamma_{j}^+,\quad
  \Gamma_{j,2} = \overline{\Omega_j^+}\cap\Gamma,\quad \Gamma_{j,4}=\overline{\Omega_j^+}\cap\Gamma_{H+L}^+.
\end{align*}
Consider the following boundary value problem for a generic field $\tu$:
\begin{equation*}
  {\noindent\rm (BVP2)}:\quad \left\{
    \begin{array}{ll}
      \nabla\cdot({\bf A} \nabla \tu) + k^2\alpha \tu =0,&{\rm on}\ \Omega_j^+,\\
      \tu = 0,&{\rm on}\ \Gamma_{j,2}\cup \Gamma_{j,4},\\
      \partial_{\nu_c}\tu = g_{i},&{\rm on}\ \Gamma_{i}^+,i=j-1,j,
    \end{array}
\right.
\end{equation*}
for $g_{i}\in H^{-1/2}(\Gamma_{i}^+), i=j-1,j$. We have the following well-posedness theorem.
\begin{mytheorem}
  \label{thm:wp:bvpntd}
  \tcr{Provided that $kT/\pi\notin {\cal E}:=\{i'/2^{j'}|j'\in
    {\mathbb N}, i'\in{\mathbb N}^*\}$},
 and $L$ is sufficiently large, (BVP2) is well-posed. \tcr{The well-posedness even holds} with $\Omega_j^+$ replaced by the interior
  domain of $2^l$ consecutive cells, say $\cup_{j=1}^{2^l}\overline{\Omega_j^+}$,
  for any number $l\geq 0$.
  \begin{proof}
    It is clear that only uniqueness is needed \cite[Thm. 4.10]{mcl00}. Suppose
    $j=1$ and $g_i=0, i=0,1$. Then, by first \tcr{an even extension over
      $\Gamma_0^+$} and then a $2T$-periodic extension, we get a $2T$-periodic
    solution $\tu^e$ (corresponding to a normal incidence) in a strip bounded in
    the $x_2$-direction by a $2T$-periodic grating surface, \tcr{possibly
      different from $\Gamma$}, and $\Gamma_{H+L}$. However, according to the
    well-posedness theory \cite[Cor. 5.2]{chaels10} for the half-space
    scattering by the grating, the PML convergence theory in \cite[Thm.
    2.4]{chewu03} can be readily adapted here to show that $\tu^e\equiv 0$,
    considering that $kT/\pi\notin{\cal E}$ has excluded horizontally
    propagating Bloch modes.
  \end{proof}
\end{mytheorem}
\begin{myremark}
  We note that the condition $kT/\pi\notin{\cal E}$ is not necessary \tcr{for
    the well-posedness of (BVP2)}. Even \tcr{if} $kT/\pi\in{\cal E}$, one may impose zero
  Neumann condition on $\Gamma_{4,j}$ to guarantee the uniqueness of \tcr{the modified} (BVP2)
  \cite{luluson19,zhowu18}.
\end{myremark}
By Theorem~\ref{thm:wp:bvpntd}, we can define a bounded Neumann-to-Dirichlet operator ${\cal
  N}^{(0)}: H^{-1/2}(\Gamma_{j-1}^+)\times H^{-1/2}(\Gamma_{j}^+)\to
\widetilde{H^{1/2}}(\Gamma_{j-1}^+)\times \widetilde{H^{1/2}}(\Gamma_{j}^+)$
such that
\begin{align}
  \label{eq:ntd:j}
  \left[
  \begin{array}{ll}
    \tu|_{\Gamma_{j-1}^+}\\
    \tu|_{\Gamma_{j}^+}\\
  \end{array}
  \right]={\cal N}^{(0)}\left[
  \begin{array}{ll}
    \partial_{\nu_c^-}\tu|_{\Gamma_{j-1}^+}\\
    \partial_{\nu_c^+}\tu|_{\Gamma_{j}^+}\\
  \end{array}
  \right],
\end{align}
for all $j\geq 1$. Due to the invariant shape of $\Omega_j^+$ with respect to
$j$, ${\cal N}^{(0)}$ is in fact independent of $j$. Suppose $j=1$. Then, by \tcr{the}
linearity principle, ${\cal N}^{(0)}$ can be rewritten in the following matrix
form
\[
{\cal N}^{(0)} = \left[
  \begin{array}{ll}
    {\cal N}^{(0)}_{00} & {\cal N}^{(0)}_{01}\\
    {\cal N}^{(0)}_{10} & {\cal N}^{(0)}_{11}\\
  \end{array}
  \right],
\]
where the bounded map ${\cal N}^{(0)}_{i'j'}: H^{-1/2}(\Gamma_{j'}^+)\to
\widetilde{H^{1/2}}(\Gamma_{i'})$ maps $\partial_{\nu_c}\tu|_{\Gamma_{j'}^+}=g_{j'}$ to
$\tu|_{\Gamma_{i'}^+}$ \tcr{if} $g_{1-j'}=0$ for $i',j'=0,1$. 

Due to the shape invariance of $\Gamma_{j}^{+}$, we shall identify
$H^{-1/2}(\Gamma_{j}^+)$ for all $j\geq 0$ as the same space
$H^{-1/2}(\Gamma_{0}^+)$, \tcr{and, similarly, the $\widetilde{H^{1/2}}(\Gamma_{j}^+)$ shall} all be
identified as the dual space of $H^{-1/2}(\Gamma_{0}^+)$. 

Returning back to the semi-waveguide problems ($P^\pm$), we have, by the
definition of ${\cal R}_p^+$ and (\ref{eq:ntd:j}) for $j=1$ and $2$, that
\begin{align}
  {\cal N}^{(0)}_{10}\partial_{\nu_c^-}\tu^{\rm og}|_{\Gamma_{0}^+} -{\cal N}^{(0)}_{11}{\cal R}_p^+\partial_{\nu_c^-}\tu^{\rm og}|_{\Gamma_{0}^+}= \tu^{\rm og}|_{\Gamma_{1}^+} = {\cal N}^{(0)}_{00}{\cal R}_p^+\partial_{\nu_c^-}\tu^{\rm og}|_{\Gamma_{0}^+} -{\cal N}^{(0)}_{01}({\cal R}_p^+)^2\partial_{\nu_c^-}\tu^{\rm og}|_{\Gamma_{0}^+}.
\end{align}
\tcr{Here and in the following}, the product of two operators should be regarded as
their composition. Thus,
\[
\left[  {\cal N}^{(0)}_{10}+ {\cal N}^{(0)}_{11}{\cal R}_p^++{\cal N}^{(0)}_{00}{\cal R}_p^++{\cal
    N}^{(0)}_{01}({\cal R}_p^+)^2\right]\partial_{\nu_c}\tu^{\rm og}|_{\Gamma_{0}^+}=0,
\]
for any $\partial_{\nu_c}\tu^{\rm og}|_{\Gamma_{0}^+}\in H^{-1/2}(\Gamma_0^+)$,
so that we end up with the following \tcr{Riccati} equation for ${\cal R}_p^+$:
\begin{equation}
  \label{eq:ric:Rp+}
{\cal N}^{(0)}_{10}+ [{\cal N}^{(0)}_{11}+{\cal N}^{(0)}_{00}]{\cal R}_p^++{\cal N}^{(0)}_{01}({\cal R}_p^+)^2 = 0.
\end{equation}
One similarly obtains the governing equation for ${\cal R}_p^-$:
\begin{equation}
  \label{eq:ric:Rp-}
{\cal N}^{(0)}_{01}+ [{\cal N}^{(0)}_{11}+{\cal N}^{(0)}_{00}]{\cal R}_p^-+{\cal N}^{(0)}_{10}({\cal R}_p^-)^2 = 0.
\end{equation}
Analogous to \cite{jollifli06}, the previous results in fact indicate that the two
\tcr{Riccati} equations (\ref{eq:ric:Rp+}) and (\ref{eq:ric:Rp-}) must be uniquely
solvable under the condition that $\rho(R_p^\pm)<1$. 
The vNtD operators ${\cal N}^{\pm}$ mapping $\partial_{\nu_c}\tu^{\rm og}|_{\Gamma_0^{\pm}}$ to $\tu^{\rm og}|_{\Gamma_0^{\pm}}$ are respectively
given by
\begin{align}
  {\cal N}^{+} =& {\cal N}^{(0)}_{00} - {\cal N}^{(0)}_{01}{\cal R}_p^+,\\
  {\cal N}^{-} =& {\cal N}^{(0)}_{11} -{\cal N}^{(0)}_{10}{\cal R}_p^-.
\end{align}
However, due to the nonlinearity of the \tcr{Riccati} equations (\ref{eq:ric:Rp+}) and
(\ref{eq:ric:Rp-}), it is not that easy to get ${\cal N}^\pm$ in practice \cite{jollifli06}. To
tackle this difficulty, we shall develop \tcr{an} RDP to effectively approximate ${\cal R}_p^\pm$.

\subsection{Recursive doubling procedure}

Take ${\cal R}_p^+$ as an example. We first study the NtD
operator
\begin{equation}
{\cal N}^{(l)}=\left[
  \begin{array}{ll}
    {\cal N}^{(l)}_{00} & {\cal N}^{(l)}_{01}\\
    {\cal N}^{(l)}_{10} & {\cal N}^{(l)}_{11}\\
  \end{array}
    \right]
\end{equation}
on the boundary of $\cup_{j=1}^{2^l}\overline{\Omega_j^+}$ for $l\geq 1$,
\tcr{where ${\cal N}^{(l)}_{i'j'}$ is bounded from $H^{-1/2}(\Gamma_{0}^+)$ to
  $\widetilde{H^{1/2}}(\Gamma_{0})$ for $i',j'=0,1$}. \tcr{If} $l=1$, we need to
compute ${\cal N}^{(1)}$ on the boundary of
$\overline{\Omega_1^+\cup\Omega_2^+}$. Using (\ref{eq:ntd:j}) for $j=1$ and $2$
and eliminating $\tu^{\rm og}$ and $\partial_{\nu_c} \tu^{\rm og}$ by the
continuity condition on $\Gamma_1^+$, one gets
\begin{equation}
  \label{eq:ndr02to1}
     ({\cal N}^{(l-1)}_{00}+{\cal N}^{(l-1)}_{11} )\partial_{\nu_c^+}\tu^{\rm og}|_{\Gamma_{1}^+}= -{\cal N}^{(l-1)}_{10}\partial_{\nu_c^-}\tu^{\rm og}|_{\Gamma_{0}^+}+{\cal N}^{(l-1)}_{01}\partial_{\nu_c^+}\tu^{\rm og}|_{\Gamma_{2}^+}.
\end{equation}
By Theorem~\ref{thm:wp:bvpntd} , the well-posedness of \tcr{the modified }(BVP2) for $l=1$,
indicates that there exist two bounded operators ${\cal A}_{l-1},{\cal
  B}_{l-1}:H^{-1/2}(\Gamma_0^+)\to H^{-1/2}(\Gamma_0^+)$ such that
\[
  \partial_{\nu_c^+}\tu^{\rm og}|_{\Gamma_{1}^+}= -{\cal
    A}_{l-1}\partial_{\nu_c^-}\tu^{\rm og}|_{\Gamma_{0}^+}+{\cal
    B}_{l-1}\partial_{\nu_c^+}\tu^{\rm og}|_{\Gamma_{2}^+}.
\]
Equation (\ref{eq:ndr02to1}) implies that
\[
{\cal A}_{l-1} = ({\cal N}^{(l-1)}_{00}+{\cal N}^{(l-1)}_{11} )^{-1}{\cal N}_{10}^{(l-1)},\quad{\cal B}_{l-1} = ({\cal N}^{(l-1)}_{00}+{\cal N}^{(l-1)}_{11} )^{-1}{\cal N}_{01}^{(l-1)},
\]
\tcr{where $({\cal N}^{(l-1)}_{00}+{\cal N}^{(l-1)}_{11} )^{-1}$ is a
  generalized inverse from $\widetilde{H^{1/2}}(\Gamma_{0})$ to
  $H^{-1/2}(\Gamma_{0}^+)$.} Thus, one obtains
\begin{align}
  \label{eq:rdp:1}
  {\cal N}_{00}^{(l)}=&{\cal N}^{(l-1)}_{00}-{\cal N}^{(l-1)}_{01}{\cal A}_{l-1},\quad 
  {\cal N}_{01}^{(l)}={\cal N}^{(l-1)}_{01}{\cal B}_{l-1},\\
  \label{eq:rdp:4}
  {\cal N}_{10}^{(l)}=&{\cal N}^{(l-1)}_{10}{\cal A}_{l-1},\quad
  {\cal N}_{11}^{(l)}={\cal N}^{(l-1)}_{11}-{\cal N}^{(l-1)}_{10}{\cal B}_{l-1}.
\end{align}
Equations (\ref{eq:rdp:1}-\ref{eq:rdp:4}) can be recursively applied
to get ${\cal N}^{(l)}$ for all $l\geq 1$,
and the number of consecutive
cells $\{\Omega_j\}$ doubles after each iteration, which form the origin of the term
``{\it recursive doubling procedure}'' (RPD) in the literature \cite{yualu07, ehrhanzhe09}. In the
following, we shall see that RDP provides a simple approach for solving
(\ref{eq:ric:Rp+}) and (\ref{eq:ric:Rp-}). 

Now, analogous to (\ref{eq:ric:Rp+}), we obtain from ${\cal N}^{(l)}$ and
(\ref{eq:R+-}) the following equations
\begin{align}
  \label{eq:ric:Rpl+}
&{\cal N}^{(l)}_{10}+ [{\cal N}^{(l)}_{11}+{\cal N}^{(l)}_{00}]({\cal R}_p^+)^{2^l}+{\cal N}_{01}^{(l)}({\cal R}_p^+)^{2^{(l+1)}} = 0,\\
  \label{eq:vNtD:l}
&{\cal N}^{+} = {\cal N}^{(l)}_{00} - {\cal N}^{(l)}_{01}({\cal R}_p^+)^{2^l}.
\end{align}
Since $||({\cal R}_p^+)^{N_0}||<1$, the third term in (\ref{eq:ric:Rpl+}) is expected to
be exponentially small for $l\gg \log_2{N_0}$, so that we approximate
\begin{align}
&({\cal R}_p^+)^{2^l} \approx -[{\cal N}^{(l)}_{11}+{\cal N}^{(l)}_{00}]^{-1} {\cal N}^{(l)}_{10},\\
  \label{eq:N+:iter}
&{\cal N}^{+} \approx {\cal N}^{(l)}_{00} +{\cal N}^{(l)}_{01}[{\cal N}^{(l)}_{11}+{\cal N}^{(l)}_{00}]^{-1} {\cal N}^{(l)}_{10},
  \end{align}
and we get ${\cal R}_p^+$ iteratively from
\begin{align}
  \label{eq:Rp+:iter}
({\cal R}_p^+)^{2^{j}}=-[{\cal N}^{(j)}_{11}+{\cal N}^{(j)}_{00}]^{-1}\left[  {\cal N}^{(j)}_{10}-{\cal N}_{01}^{(j)}({\cal R}_p^+)^{2^{j+1}}\right],\quad j=l-1,\cdots,0.
\end{align}
One similarly obtains ${\cal N}^-$ and ${\cal R}_p^-$ from
\begin{align}
&({\cal R}_p^-)^{2^l} \approx -[{\cal N}^{(l)}_{11}+{\cal N}^{(l)}_{00}]^{-1} {\cal N}^{(l)}_{01},\\
  \label{eq:N-:iter}
&{\cal N}^{-} \approx{\cal N}^{(l)}_{11} +{\cal N}^{(l)}_{01}[{\cal N}^{(l)}_{11}+{\cal N}^{(l)}_{00}]^{-1} {\cal N}^{(l)}_{10},\\
  \label{eq:Rp-:iter}
&({\cal R}_p^-)^{2^{j}}=-[{\cal N}^{(j)}_{11}+{\cal N}^{(j)}_{00}]^{-1}\left[  {\cal N}^{(j)}_{01}-{\cal N}_{10}^{(j)}({\cal R}_p^-)^{2^{j+1}}\right],\quad j=l-1,\cdots,0.
\end{align}

From the above, it can be seen that the essential step to approximate ${\cal
  N}^\pm$ is to get the NtD operator ${\cal N}^{(0)}$ on the boundary of any
unperturbed unit cell $\Omega_j^\pm$ for $j\in\mathbb{Z}^+$. As no information
of the field $\tu^{\rm og}$ in $\Omega_j^\pm$ is required, it is clear that
\tcr{the} BIE method is an optimal choice, as it treats only the boundary of
$\Omega_j^\pm$. Since PML is involved in domain $\Omega_j^\pm$, the
high-accuracy PML-based BIE method developed in our previous work
\cite{luluqia18} straightforwardly provides an accurate approximation of ${\cal
  N}^{(0)}$, so as to effectively drive RDP to get ${\cal N}^\pm$. \tcr{We shall
  present the details in the next section.}

\section{The PML-based BIE method}
In this section, we shall first review the PML-based BIE method in
\cite{luluqia18} to approximate the NtD operator on the boundary of any unit
cell, perturbed or not, by an NtD matrix. Then, we shall use these NtD matrices
to approximate the two vNtD operators ${\cal N}^{\pm}$ on $\Gamma_0^{\pm}$ and
to solve (BVP1) finally. From now on, we shall assume that the scattering
surface $\Gamma$ is piecewise smooth and satisfies (GC1) only. Though the
previous well-posedness theory relies on (GC2), our numerical solver does not
rely on such an assumption, and we believe (GC2) can be weakened to at least
accept piecewise smooth curves, which we shall investigate in a future work.
\subsection{Approximating ${\cal N}^{\pm}$}
Without loss of generality, consider (BVP2) in an unperturbed cell, say
$\Omega_1^+$, and we need to approximate ${\cal N}^{(0)}$ first. According to
\cite{luluqia18}, for any $\tu$ satisfying
\begin{equation}
  \label{eq:helm:hom}
  \nabla\cdot({\bf A}\nabla \tu) + k^2\alpha \tu = 0,
\end{equation}
on $\Omega_1^+$, we have the following Green's representation theorem
\begin{align}
	\label{eq:grerep}
	\tilde{u}(x)= \int_{\partial \Omega_1^+}\{\tilde{G}(x,y)\partial_{\nu_c}\tilde{u}(y)-\partial_{\nu_c}\tilde{G}(x,y)\tilde{u}(y)\}ds(y),
\end{align} 
for all $x\in \Omega_1^+$; we recall that $\nu$ denotes the outer unit
normal vector on $\partial\Omega_1^+$. Moreover, as $x$ approaches $\partial \Omega_1^+=\cup_{j=1}^{4}\overline{\Gamma_{j,1}}$,
the usual \tcr{jump} conditions imply \tcr{\cite{luluqia18}}
\begin{align}
  \label{eq:ntd:bie}
	{\cal K}[\tu](x) -{\cal K}_0 [1](x)\tu(x)  = {\cal S} 
\partial_{{\nu}_c} [\tilde{u}] (x),
\end{align}
where we have defined the following integral operators
\begin{align}
  \label{eq:def:S}
	{\cal S}[\phi](x) &= 2\int_{\partial \Omega_1^+} \tilde{G}(x,y) \phi(y) ds(y),\\
  \label{eq:def:K}
	{\cal K}[\phi](x) &= 2{\rm p.v.}\int_{\partial \Omega_1^+} \partial_{{\nu_c}} \tilde{G}(x,y) \phi(y) ds(y),\\
  \label{eq:def:K0}
	{\cal K}_0[\phi](x) &= 2{\rm p.v.}\int_{\partial \Omega_1^+} \partial_{{\nu_c}} \tilde{G}_0(x,y) \phi(y) ds(y),
\end{align}
where p.v. indicates the Cauchy principle value, and
\begin{equation}
	\label{eq:green:complap}
	\tilde{G}_0(x,y) = -\frac{1}{2\pi}\log\rho(\tilde{x},\tilde{y}),
\end{equation}
is the fundamental solution of the complexified Laplace equation
\begin{equation}
	\label{eq:cmp:lap}
	\nabla\cdot({\bf A} \nabla \tilde{u}_0(x)) = 0.
\end{equation}
Note that theoretically,
\begin{equation}
  \label{eq:k01:exa}
  {\cal K}_0 [1](x) = - \frac{\theta^{\rm in}(x)}{\pi}.
\end{equation}
where $\theta^{\rm in}(x)$ is defined as the interior angle at $x$, as indicated
in Figure~\ref{fig:lateral}(b). However, numerically evaluating ${\cal K}_0[1]$ near corners
is more advantageous as has been illustrated in the literature
\cite{colkre13,lulu14}. Thus, $\tu = ({\cal K}-{\cal K}_0 [1])^{-1}{\cal
  S}\partial_{\nu_c}\tu$ on $\partial \Omega_1^+$. Consequently, the NtD
operator ${\cal N}_{u}$ for any unperturbed domain can be defined as
\[
{\cal N}_u = ({\cal K}-{\cal K}_0 [1])^{-1}{\cal S}.
\]

To approximate ${\cal N}_u$, we need \tcr{to} discretize the three integral
operators on the right-hand side. Suppose now the piecewise smooth curve $\partial \Omega_1^+$ is
parameterized by $x(s)=\{(x_1(s), x_2(s))|0\leq s\leq \tcr{L_1}\}$, where $s$ is the
arclength parameter. Since corners may exist, $\tilde{u}(x(s))$ can have corner
singularities in its derivatives at corners. To smoothen $\tilde{u}$, we
introduce a \tcr{grading} function $s=w(t), 0\leq t\leq 1$. For a smooth segment
of $\partial \Omega_1^+$ corresponding to $s\in[s^0, s^1]$ and $t\in[t^0, t^1]$
such that $s^i=w(t^i)$ for $i=0,1$, where $s^0$ and $s^1$ correspond to two
corners, we take \cite[Eq. (3.104)]{colkre13}
\begin{equation}
	\label{eq:gfun}
	s=w(t) = \frac{s^0w_1^p + s^1w_2^p}{w_1^p+w_2^p},\quad t\in[t^0,t^1],
\end{equation}
where the positive integer $p$ ensures \tcr{that} the derivatives of $w(t)$ vanish at the corners up to order $p$,
\[
	w_1=\left(\frac{1}{2}-\frac{1}{p}\right)\xi^3+\frac{\xi}{p}+\frac{1}{2},\quad
	w_2 = 1-w_1,\quad \xi = \frac{2 t - (t^0+t^1)}{t^1-t^0}.
\]

To simplify \tcr{notation}, we shall use $x(t)$ to denote $x(w(t))$, and $x'(t)$
to denote $\frac{dx}{ds}(w(t))w'(t)$ in the following. Assume that $t\in[0,1]$
is uniformly sampled by an even number, denoted by $N$, of grid points
$\{t_j=jh\}_{j=1}^{N}$ with grid size $h=1/N$, and that the grid points contain
\tcr{all the} corner points. Thus, ${\cal S}[\partial_{\nu_c}\tu]$ at point
$x=x(t_j)$ can be parameterized by
\begin{align}
	\label{eq:intS}
	{\cal S}[\partial_{{\nu}_c}\tilde{u}] (x(t_j)) &= \int_0^1
	S(t_j,t) \phi^{\rm s}(t) dt,
\end{align}
where \tcr{$S(t_j,t)=\frac{\bi}{2} H_0^{(1)}(k\rho(x(t_j),x(t)))$, and the scaled co-normal vector $\phi^{\rm s}(t)=
\partial_{\nu_c}\tilde{u}(x(t))|x'(t)|$, smoother than
$\partial_{\nu_c}\tilde{u}(x(t))$, is introduced to regularize the approximation
of ${\cal N}_u$.}

Considering the logarithmic singularity of $S(t_j,t)$ at $t=t_j$, we can
discretize the integral in (\ref{eq:intS}) by Alpert's 6th-order hybrid Gauss-trapezoidal
quadrature rule \cite{alp99} and then by trigonometric interpolation to get
\begin{equation}
	{\cal S}[\partial_{\nu_c}\tilde{u}^s]\left[
		\begin{array}{c}
			x(t_1)\\
			\vdots\\
			x(t_N)
		\end{array}
	\right] \approx {\bf S} \left[
		\begin{array}{c}
			\phi^{\rm s}(t_1)\\
			\vdots\\
			\phi^{\rm s}(t_N)
		\end{array}
	\right],
\end{equation}
where the $N\times N$ matrix ${\bf S}$ approximates ${\cal S}$. One similarly
approximates ${\cal K}[\tu](x(t_j))$ and ${\cal K}_0[1](x(t_j))$ for
$j=1,\cdots, N$, so that we obtain, on the boundary of $\partial\Omega_1^+$,
\begin{equation}
\label{eq:ntd:cell}
  \left[
    \begin{array}{c}
      {\bm u}_{1,1}\\
      {\bm u}_{1,2}\\
      {\bm u}_{1,3}\\
      {\bm u}_{1,4}
    \end{array}
  \right]
={\bm N}_u\left[
    \begin{array}{c}
      {\bm \phi}_{1,1}^{\rm s}\\
      {\bm \phi}_{1,2}^{\rm s}\\
      {\bm \phi}_{1,3}^{\rm s}\\
      {\bm \phi}_{1,4}^{\rm s}
    \end{array}
  \right],
\end{equation}
where ${\bm u}_{1,j'}$ and ${\bm \phi}^s_{1,j'}$ represent $N_{j'}\times 1$ column
vectors of $\tu$ and $\phi^{\rm s}$ at the $N_{j'}$ grid points of $\Gamma_{1,j'}$,
respectively for $j'=1,2,3,4$; note that $N=\sum_{j'=1}^4N_{j'}$ and the grid points
on $\Gamma_{1,3}$ are obtained by horizontally translating the grid points on
$\Gamma_{1,1}$ so that $N_1=N_3$. Clearly, the $N\times N$ matrix ${\bm N}_u$
approximates the scaled NtD operator ${\cal N}^{\rm s}_u$ related to ${\cal
  N}_u$ by ${\cal N}_u\partial_{\nu_c} \tu = {\cal N}_u^{\rm s}\phi^{\rm s}$.
Now, by $\tu|_{\Gamma_{1,2}\cup\Gamma_{1,4}}=0$, \tcr{we eliminate vectors ${\bm
  u}_{1,2}$, ${\bm u}_{1,4}$, ${\bm \phi}_{1,2}$ and ${\bm \phi}_{1,4}$ in
(\ref{eq:ntd:cell})} so that we obtain two $2N_1\times 2N_1$ matrices ${\bm
  N}^{(0)}$ and ${\bm T}$ that satisfy
\begin{align}
  \label{eq:aprNtD}
  \left[
  \begin{array}{ll}
    {\bm u}_{1,1}\\
    {\bm u}_{1,3}\\
  \end{array}
  \right]={\bm N}^{(0)}\left[
  \begin{array}{ll}
    {\bm \phi}^{\rm s}_{1,1}\\
    {\bm \phi}^{\rm s}_{1,3}\\
  \end{array}
  \right],\quad \left[
  \begin{array}{ll}
    {\bm \phi}^{\rm s}_{1,2}\\
    {\bm \phi}^{\rm s}_{1,4}\\
  \end{array}
  \right]={\bm T}\left[
  \begin{array}{ll}
    {\bm \phi}^{\rm s}_{1,1}\\
    {\bm \phi}^{\rm s}_{1,3}\\
  \end{array}
  \right],
\end{align}
where we denote
\[
{\bm N}^{(0)} = \left[
  \begin{array}{ll}
    {\bm N}^{(0)}_{00} & {\bm N}^{(0)}_{01}\\
    {\bm N}^{(0)}_{10} & {\bm N}^{(0)}_{11}\\
  \end{array}
  \right],
\]
with ${\bm N}^{(0)}_{ij}\in {\mathbb C}^{N_1\times N_1}$; \tcr{the above
  elimination is stable due to the well-posedness of (BVP2)
  in Theorem~\ref{thm:wp:bvpntd}}. Note that, different
from \cite{luluqia18}, we no longer simultaneously assume $\tu=\phi^{\rm s}=0$
on ${\Gamma_{1,2}\cup\Gamma_{1,4}}$, which could cause pronounced error in
numerical results. Now compare \eqref{eq:ntd:j} and (\ref{eq:aprNtD}). Like
\tcr{${\bm N}_u$}, ${\bm N}^{(0)}$ approximates the scaled NtD operator ${\cal
  N}^{(0),\rm s}$ on $\Gamma_{1}^+\cup\Gamma_{3}^+$ related to ${\cal N}^{(0)}$
by ${\cal N}^{(0)}\partial_{\nu_c}\tu = {\cal N}^{(0),\rm s}\phi^s$.

Consequently, the previously developed RDP can be
easily adapted here in terms of notationally replacing ${\cal N}$ by ${\bm N}$
for the equations (\ref{eq:rdp:1}-\ref{eq:Rp-:iter}), so that we get two $N_1\times N_1$ matrices
${\bm R}_p^{+}$ and ${\bm N}^{+}$ approximating \tcr{the} (scaled) Neumann-marching
operator ${\cal R}_p^+$ and \tcr{the} (scaled) vNtD operator ${\cal N}^+$ such that ${\bm
  \phi}^{\rm s}_{1,3}=-{\bm R}_p^+{\bm \phi}^{\rm s}_{1,1}$ and ${\bm
  u}_{1,1}={\bm N}^+{\bm \phi}^{\rm s}_{1,1}$. One similarly obtains two $N_1\times N_1$
matrices ${\bm R}_p^-$ and ${\bm N}^-$ approximating ${\cal R}_p^-$ and ${\cal
  N}^-$, respectively.

\subsection{Solving (BVP1)}

We are now ready to use the PML-based BIE method to solve the main problem
(BVP1). For $x^*\in \Omega_0$, to eliminate \tcr{the} $\delta$ function, we consider
$\tu^{\rm sc}(x;x^*) = \tu^{\rm og}(x;x^*)-\tu^{\rm inc}(x;x^*)$, satisfying
(\ref{eq:helm:hom}). For simplicity, we denote \tcr{(cf.
  Fig.\ref{fig:semi:problem} (a))}
\begin{align*}
  \Gamma_{0,1}=\Gamma_0^-,\quad \Gamma_{0,2}=\Gamma_0,\quad \Gamma_{0,3}=\Gamma_0^+,\quad{\rm and}\quad \Gamma_{0,4}=\Gamma_0^{H+L}.
\end{align*}
Then, analogous to (\ref{eq:ntd:cell}), on the four boundaries $\Gamma_{0,j}, j=1,2,3,4$, we apply the
PML-based BIE method in the previous section to approximate the NtD operator for \tcr{$\tu^{\rm sc}$}
and \tcr{$\partial_{\nu_c}\tu^{\rm sc}$} on the boundary of \tcr{the} perturbed cell $\Omega_0$ by a
matrix ${\bm N}_p$,
\begin{align}
  \label{eq:Np}
\left[
    \begin{array}{c}
      {\bm u}^{\rm sc}_{0,1}\\
      {\bm u}^{\rm sc}_{0,2}\\
      {\bm u}^{\rm sc}_{0,3}\\
      {\bm u}^{\rm sc}_{0,4}
    \end{array}
  \right]
={\bm N}_p\left[
    \begin{array}{c}
      {\bm \phi}_{0,1}^{\rm sc,s}\\
      {\bm \phi}_{0,2}^{\rm sc,s}\\
      {\bm \phi}_{0,3}^{\rm sc,s}\\
      {\bm \phi}_{0,4}^{\rm sc,s}
    \end{array}
  \right],
\end{align}
 where ${\bm u}^{\rm sc}_{0,j}$ and ${\bm \phi}^{\rm sc,s}_{0,j}$ represent column vectors
 of $\tu^{\rm sc}$
and $\partial_{\nu_c}\tu^{\rm sc}|x'|$ at the grid points of $\Gamma_{0,j}$,
respectively, for $j=1,2,3,4$.
Rewriting the above in terms of $\tu^{\rm og}$ and $\partial_{\nu_c} \tu^{\rm og}$,
we get
\begin{align}
  \label{eq:linsys:1}
\left[
    \begin{array}{c}
      {\bm u}^{\rm og}_{0,1}\\
      {\bm u}^{\rm og}_{0,2}\\
      {\bm u}^{\rm og}_{0,3}\\
      {\bm u}^{\rm og}_{0,4}
    \end{array}
  \right]
={\bm N}_p\left[
    \begin{array}{c}
      {\bm \phi}_{0,1}^{\rm og,s}\\
      {\bm \phi}_{0,2}^{\rm og,s}\\
      {\bm \phi}_{0,3}^{\rm og,s}\\
      {\bm \phi}_{0,4}^{\rm og,s}
    \end{array}
  \right] + \left[
    \begin{array}{c}
      {\bm u}^{\rm inc}_{0,1}\\
      {\bm u}^{\rm inc}_{0,2}\\
      {\bm u}^{\rm inc}_{0,3}\\
      {\bm u}^{\rm inc}_{0,4}
    \end{array}
  \right]-{\bm N}_p\left[
    \begin{array}{c}
      {\bm \phi}_{0,1}^{\rm inc,s}\\
      {\bm \phi}_{0,2}^{\rm inc,s}\\
      {\bm \phi}_{0,3}^{\rm inc,s}\\
      {\bm \phi}_{0,4}^{\rm inc,s}
    \end{array}
  \right],
  \end{align}
  where ${\bm u}^{\rm inc}_{0,j}$ and ${\bm \phi}^{\rm inc,s}_{0,j}$ represent
  column vectors of $\tu^{\rm inc}(x;x^*)$ and $\partial_{\nu_c}\tu^{\rm
    inc}(x;x^*)|x'|$ at the grid points of $\Gamma_{0,j}$, respectively, etc..
  The boundary conditions in (BVP1) imply that
\begin{align}
  \label{eq:linsys:2}
{\bm u}^{\rm og}_{0,2}=&0,\quad {\bm u}^{\rm og}_{0,4}=0,\\
  \label{eq:linsys:3}
  {\bm u}^{\rm og}_{0,1}=&{\bm N}^{-}{\bm \phi}_{0,1}^{\rm og,s},\quad {\bm u}^{\rm og}_{0,3}={\bm N}^{+}{\bm \phi}_{0,3}^{\rm og,s}.
\end{align}
Solving the linear system (\ref{eq:linsys:1}-\ref{eq:linsys:3}), we get
$\tu^{\rm og}(x;x^*)$ and $\partial_{\nu_c}\tu^{\rm og}(x;x^*)$ on all grid
points of $\partial\Omega_0$.

Now we discuss how to evaluate $\tu^{\rm og}(x;x^*)$ in the physical domain $S_H$. We
distinguish two cases:
\begin{itemize}
\item[1.] $x\in\Omega_0$. Since on the grid points of $\partial\Omega_0$,
  $\tu^{\rm sc}$
  and $\partial_{\nu_c}\tu^{\rm sc}|x'|$ are available, we use Green's representation formula
  (\ref{eq:grerep}) with $\partial\Omega_1^+$ replaced by $\partial\Omega_0$ to compute
  $\tu^{\rm sc}(x;x^*)$ in $\Omega_0$ so that $\tu^{\rm og}(x;x^*)$ becomes available in $\Omega_0$.
\item[2.] $x\in\Omega_j^\pm$. Consider $\Omega_1^+$ first. Suppose ${\bm u}^{\rm
    og}_{1,j'}$ and ${\bm \phi}^{\rm og, s}_{1,j'}$ represent column vectors of
  $\tu^{\rm og}$ and $\partial_{\nu_c}\tu^{\rm og}|x'|$ at the grid points of
  $\Gamma_{1,j'}$, for $1\leq j'\leq 4$. By the continuity of
  $\partial_{\nu_c}\tu^{\rm og}$ on $\Gamma_{1,1}=\Gamma_{0,3}=\Gamma_{0}^+$,
  ${\bm \phi}^{\rm og,s}_{1,1}=-{\bm \phi}^{\rm og,s}_{0,3}$. Since ${\bm
    \phi}^{\rm og,s}_{1,3}=-{\bm R}_p^+{\bm \phi}^{\rm og,s}_{1,1}$, we get
  ${\bm u}^{\rm og}_{1,j'}$ for $j'=1,3$ by (\ref{eq:aprNtD}), and ${\bm
    \phi}^{\rm og,s}_{1,j'}$ for $j'=2,4$. Given that ${\bm u}^{\rm
    og}_{1,2}={\bm u}^{\rm og}_{1,4}=0$, $\tu^{\rm og}(x;x^*)$ and
  $\partial_{\nu_c}\tu^{\rm og}|x'|$ on $\partial\Omega_1^+$ become available,
  so that the Green's representation formula~(\ref{eq:grerep}) applies to get
  $\tu^{\rm og}(x;x^*)$ in $\Omega_1^+$. Repeating the same procedure, one
  obtains $\tu^{\rm og}(x;x^*)$ in $\Omega_j^+$ for $j\geq 2$. The case for
  $x\in\Omega_j^-$ can be handled similarly.
\end{itemize}
Consequently, $u^{\rm tot}(x;x^*)\approx \tu^{\rm og}(x;x^*)$ becomes available
for $x\in S_H\subset \overline{\Omega_0}\cup\left[  \cup_{j=1}^{\infty}\overline{\Omega_{j,+}\cup\Omega_{j,-}}\right].$
\subsection{Computing $u^{\rm tot}$ for \tcr{plane-wave} incidence}
To \tcr{close} this section, we briefly discuss how to compute $u^{\rm tot}$ for
a plane incident wave $u^{\rm inc}=e^{\bi k (\cos\theta x_1 - \sin\theta x_2)}$
for $\theta\in(0,\pi)$. \tcr{First, we consider the non-perturbed case $\Gamma =
  \Gamma_T$ so that $u^{\rm tot}$ becomes the reference solution $u^{\rm
    tot}_{\rm ref}$}. It is clear that $u_{\rm ref}^{\rm sc}=u^{\rm tot}_{\rm
  ref}-u^{\rm inc}$ satisfies the following quasi-periodic condition
\begin{align}
  \label{eq:qpc:1}
  u_{\rm ref}^{\rm sc}(-T/2,x_2) =& \gamma u_{\rm ref}^{\rm sc}(T/2,x_2),\\
  \label{eq:qpc:2}
  \partial_{x_1} u_{\rm ref}^{\rm sc}(-T/2,x_2) =& \gamma \partial_{x_1} u_{\rm ref}^{\rm sc}(T/2,x_2),
\end{align}
where $\gamma = e^{\bi k\cos\theta T}$. On $\Gamma$, we have from
(\ref{eq:gov:bc:te}) that
\begin{equation}
  \label{eq:surf:cond}
  u_{\rm ref}^{\rm sc}=-u^{\rm inc}.
\end{equation}

Due to the quasi-periodicity, above $x_2=H$, we could express $u_{\rm ref}^{\rm
  sc}$ in terms of \tcr{a} Fourier series, i.e.,
\begin{equation}
  u_{\rm ref}^{\rm sc}(x_1,x_2) = \sum_{j=-\infty}^{\infty}R_je^{\bi \alpha_j x_1 + \bi \beta_j x_2},\quad x_2\geq H,
\end{equation}
where $\alpha_j =k\cos\theta + \frac{2\pi j}{T}$ and \tcr{$\beta_j
=\sqrt{k^2-\alpha_j^2}$ if $|\alpha_j| \le k$, otherwise $\beta_j
={\bi}\sqrt{\alpha_j^2-k^2}$}, and $R_j$ denotes the $j$-th reflective coefficient.
Thus, the complexified field $\tilde{u}_{\rm ref}^{\rm sc}(x_1,x_2)=u_{\rm ref}^{\rm sc}(x_1,\tilde{x}_2)$ satisfies on the PML boundary
$x_2=L+H$
\begin{equation}
  \label{eq:sc:pmlbc}
  \tilde{u}_{\rm ref}^{\rm sc}(x_1,L+H)=\sum_{j=-\infty}^{\infty}R_je^{\bi \alpha_j x_1 + \bi \beta_j (H+L) -\beta_j \tS L}. 
\end{equation}
For simplicity, we assume that all $\beta_j$ are sufficiently away from $0$, so
that provided that $L$ and $\tS$ are sufficiently large, we can directly 
impose the following Dirichlet boundary condition
\begin{equation}
  \label{eq:sc:bc1}
  \tilde{u}_{\rm ref}^{\rm sc}(x_1,H+L)=0.
\end{equation}
If $\beta_j$ is quite close to 0, accurate boundary conditions can be developed;
we refer readers to \cite{lulu12josaa,luluson19,zhowu18} for details. Besides, $\tu^{\rm sc}_{\rm
  ref}$ satisfies the quasi-periodic conditions (\ref{eq:qpc:1}) and
(\ref{eq:qpc:2}) and the surface condition (\ref{eq:surf:cond}), but with $u$
replaced by $\tu$.

On the boundary $\partial \Omega_0$, the PML-BIE method gives, analogous to
(\ref{eq:Np}), 
\begin{equation}
\label{eq:scgov:1}
  \left[
    \begin{array}{c}
      {\bm u}^{\rm sc}_1\\
      {\bm u}^{\rm sc}_2\\
      {\bm u}^{\rm sc}_3\\
      {\bm u}^{\rm sc}_4
    \end{array}
  \right]
={\bm N}_p\left[
    \begin{array}{c}
      {\bm \phi}^{\rm sc}_1\\
      {\bm \phi}^{\rm sc}_2\\
      {\bm \phi}^{\rm sc}_3\\
      {\bm \phi}^{\rm sc}_4
    \end{array}
  \right],
\end{equation}
where ${\bm u}^{\rm sc}_{j'}$ and ${\bm \phi}^{\rm sc}_{j'}$ represent vectors of
$\tu_{\rm ref}^{\rm sc}$ and
$\partial_{\nu_c}\tu_{\rm ref}^{\rm sc}|w'|$ at the grid points of $\Gamma_{0,j'}$,
respectively, for $1\leq j'\leq 4$; note that ${\bm N}_p$ is the same as ${\bm N}_u$ in (\ref{eq:ntd:cell}) since $\Gamma=\Gamma_T$.
Equation (\ref{eq:sc:bc1}) directly implies that 
\begin{equation}
\label{eq:scgov:2}
{\bm u}_4^{\rm sc} = 0.
\end{equation}

The quasi-periodic conditions (\ref{eq:qpc:1}) and (\ref{eq:qpc:2}) imply
\begin{equation}
\label{eq:scgov:3}
  {\bm u}^{\rm sc}_3 = \gamma {\bm u}^{\rm sc}_1,\quad {\bm \phi}^{\rm sc}_3 = -\gamma {\bm \phi}^{\rm sc}_1.
\end{equation}
The interface condition (\ref{eq:surf:cond}) indicates
\begin{equation}
\label{eq:scgov:4}
  {\bm u}^{\rm sc}_2 = -{\bm u}^{\rm inc}_2,
\end{equation}
where ${\bm u}_2^{\rm inc}$ represents the vector of $u^{\rm inc}$ at grid
points of $\Gamma_{0,2}$. Solving the linear system
(\ref{eq:scgov:1}-\ref{eq:scgov:4}) gives rise to values of $\tilde{u}_{\rm
  ref}^{\rm sc}$ and $\partial_{\nu_c}\tu_{\rm ref}^{\rm sc}|w'|$ on $\partial
\Omega_0$. The Green's representation formula (\ref{eq:grerep}) can help to
compute $\tilde{u}_{\rm ref}^{\rm sc}$ in $\Omega_0$. The
quasi-periodicity helps to construct $\tilde{u}_{\rm ref}^{\rm sc}$ in any other
cells $\Omega_j^{\pm}$ for $j\in\mathbb{N}^*$. \tcr{Consequently, $u_{\rm ref}^{\rm
  tot}$ becomes available in the physical domain $S_H$.}

\tcr{Now, if $\Gamma$ is} a local perturbation of $\Gamma_T$, as $\tu_{\rm ref}^{\rm
  sc}$ is available now, one follows the same approach developed in section 6.2
to get $\tu^{\rm og}=\tu^{\rm sc}-\tu_{\rm ref}^{\rm sc}$ in any unperturbed
cell and thus $u^{\rm tot}$ in the physical region $S_{H}$. We omit the details here.
\section{Numerical examples}
In this section, we will carry out four numerical experiments to validate the
performance of the PML-based BIE method and also the proposed theory. In all
examples, we set the free-space wavelength $\lambda=1$ so that $k_0=2\pi$, and
the period $T=1$. We consider two types of \tcr{incidence}: (1) a cylindrical
incidence excited by source point $x^*=(0,1.5)$; (2) a \tcr{plane-wave} incidence of angle
$\theta$ to be specified. We suppose that only one unit cell of the background
periodic structure is perturbed. To setup the PML, we let $m=0$ in
(\ref{eq:sigma}) to define $\sigma$ for simplicity. In the RDP iterations
(\ref{eq:N+:iter}), (\ref{eq:Rp+:iter}), (\ref{eq:N-:iter}) and
(\ref{eq:Rp-:iter}), we take $l\!=\!20$. Furthermore, we choose $H=3$ and set
the computational domain to be $[-5.5,5.5]\!\times\![-2, 3]$, which contains 11
cells. To validate the accuracy of our method, we compute the relative error
\[
  E_{\rm rel}:=\frac{||({\bm \phi}^{\rm sc,s}_{2,0})^{\rm num} - ({\bm
      \phi}_{2,0}^{\rm sc,s})^{\rm exa}||_{\infty}}{||({\bm \phi}^{\rm sc,s}_{2,0})^{\rm exa}||_\infty},
\]
for ${\bm \phi}_{2,0}^{\rm sc,s}$ representing the scaled normal derivative
$|w'|\partial_{\nu} u^{\rm sc}$ on $\Gamma_{2,0}$, the perturbed part of
$\Gamma$, and for different values of $S$
and $L$ in the setup of the PML, where superscript ``num'' indicates numerical
solution, superscript ``exa'' indicates a sufficiently accurate numerical
solution or the exact solution if available.

{\bf Example 1: a flat curve.} In the first example, we assume that $\Gamma$ is
the straight line $\{x\!:x_2\!=\!0\}$. {Certainly}, we can regard such a simple
structure as a periodic structure {with period equal to {one} wavelength}. We
regard the line segment between $x_1\!=\!-0.5$ and $x_1\!=\!0.5$ on $\Gamma$ as
segment $\Gamma_{0,2}$, i.e., {as} the ``perturbed'' part. For the cylindrical
incidence, the total {wave field} $u^{\rm tot}$ is given by
\[
   u^{\rm tot}(x;x^*) = \frac{\bi}{4}\left[  H_0^{(1)}({k|x-x^*|}) -
     H_0^{(1)}({k|x-x_{\rm imag}^{*}|})\right],
\]
where the image source point  ${x_{\rm imag}^{*}}\!=\!(0,-1.5)$. Using this to compute the
scaled co-normal derivative on segment $\Gamma_{0,2}$, we get the reference solution and
can check the accuracy of our method. We discretize each smooth segment of the perturbed/unperturbed unit
cell by 600 grid points. To check how the wavenumber condition in
Theorem~\ref{thm:wp:bvpntd} affect the accuracy of our numerical solver, we
consider two values of the refractive index $n$ in $\Omega$: (1) $n=1.03$ so that
\tcr{$kT/\pi=2.06\notin\mathcal{E}$}; (2) $n=1$ so that
\tcr{$kT/\pi=2\in\mathcal{E}$}. For both cases, we compare results of Dirichlet
and Neumann boundary conditions on $\Gamma_{H+L}$.

\begin{figure}[!ht]
  \centering
  (a)\includegraphics[width=0.23\textwidth]{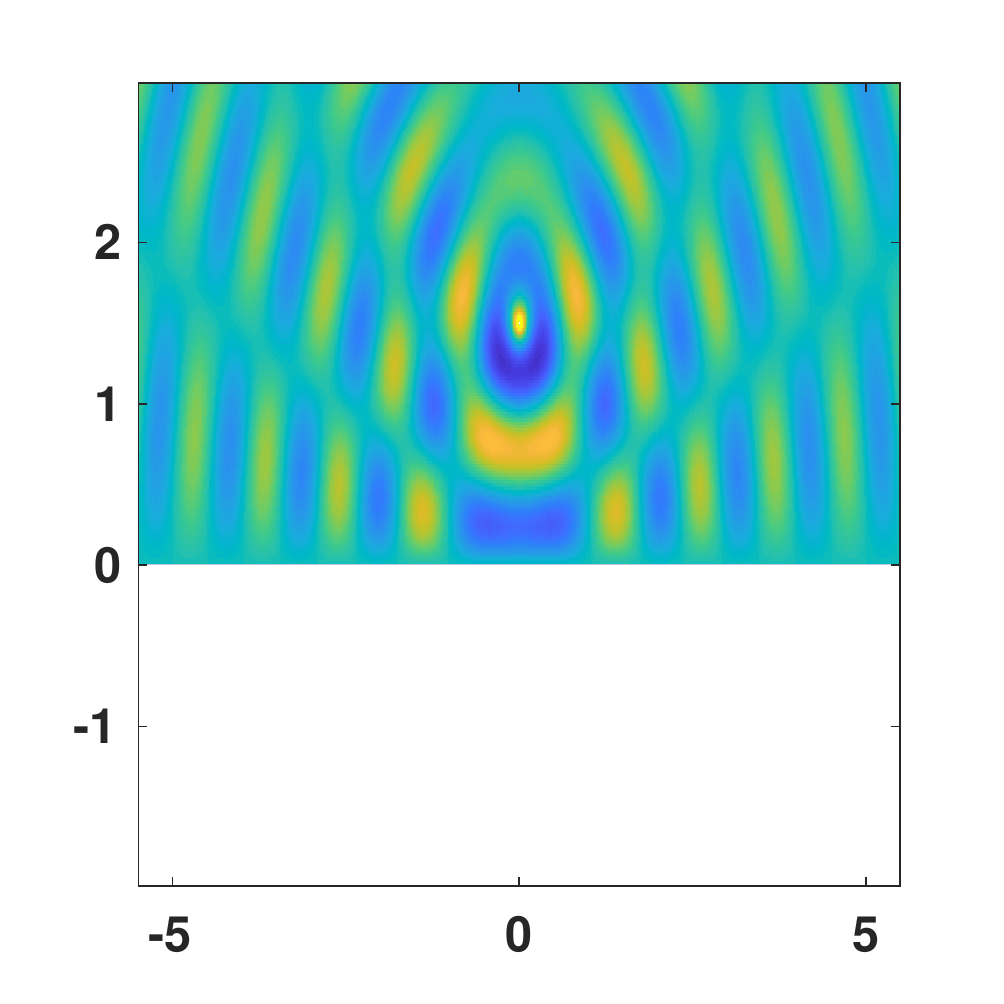}
  (b)\includegraphics[width=0.23\textwidth]{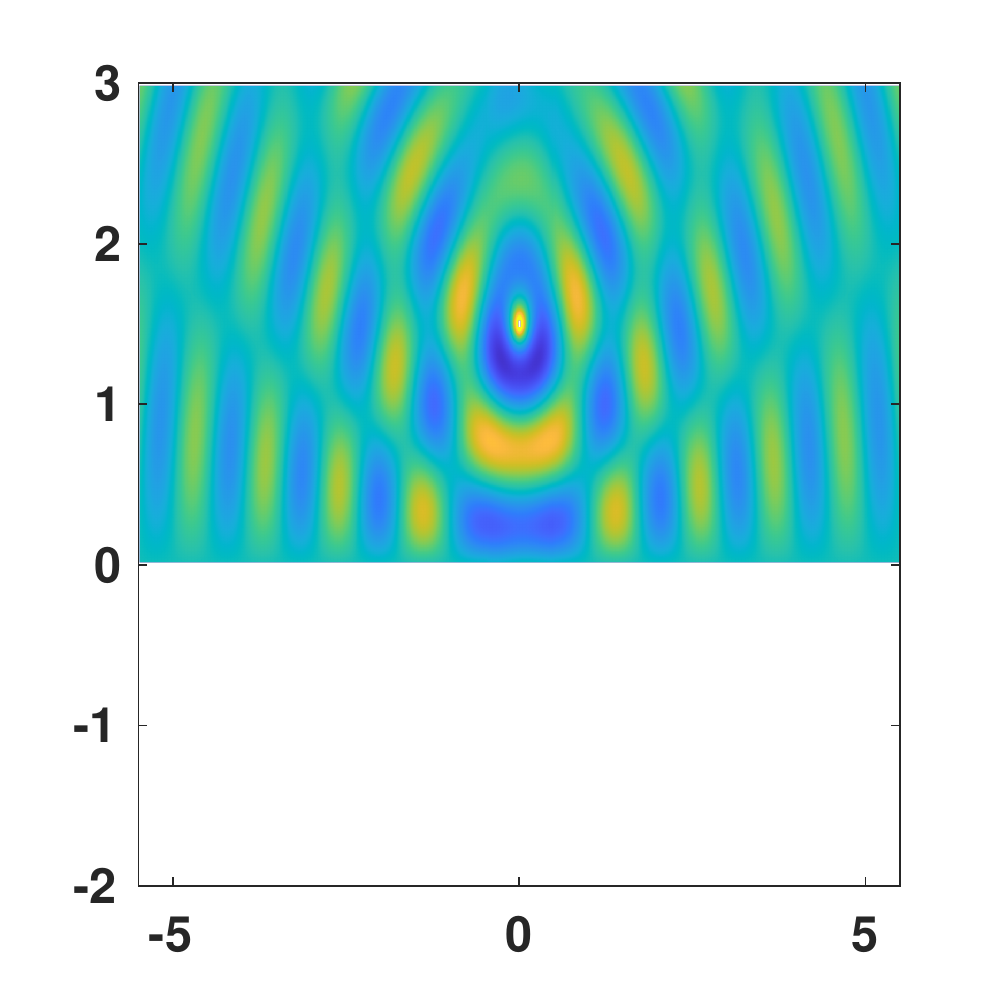}
  (c)\includegraphics[width=0.19\textwidth]{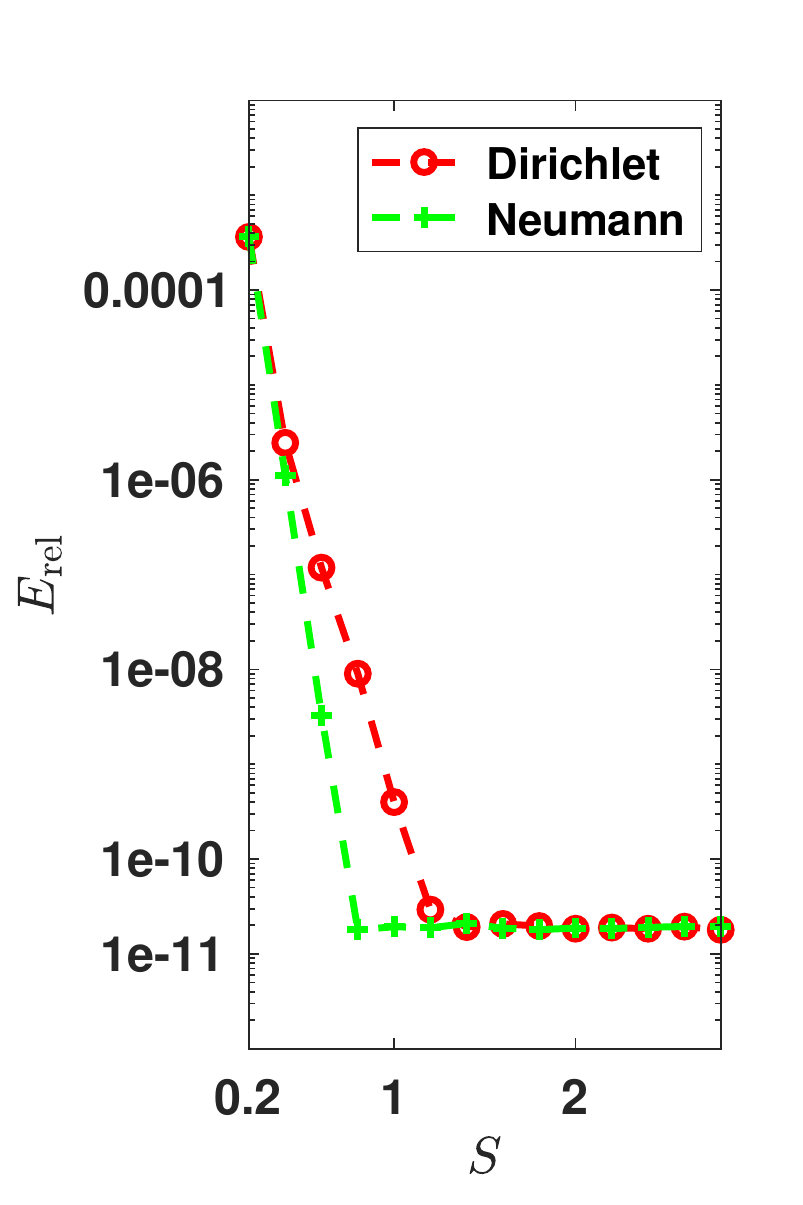}
  (d)\includegraphics[width=0.19\textwidth]{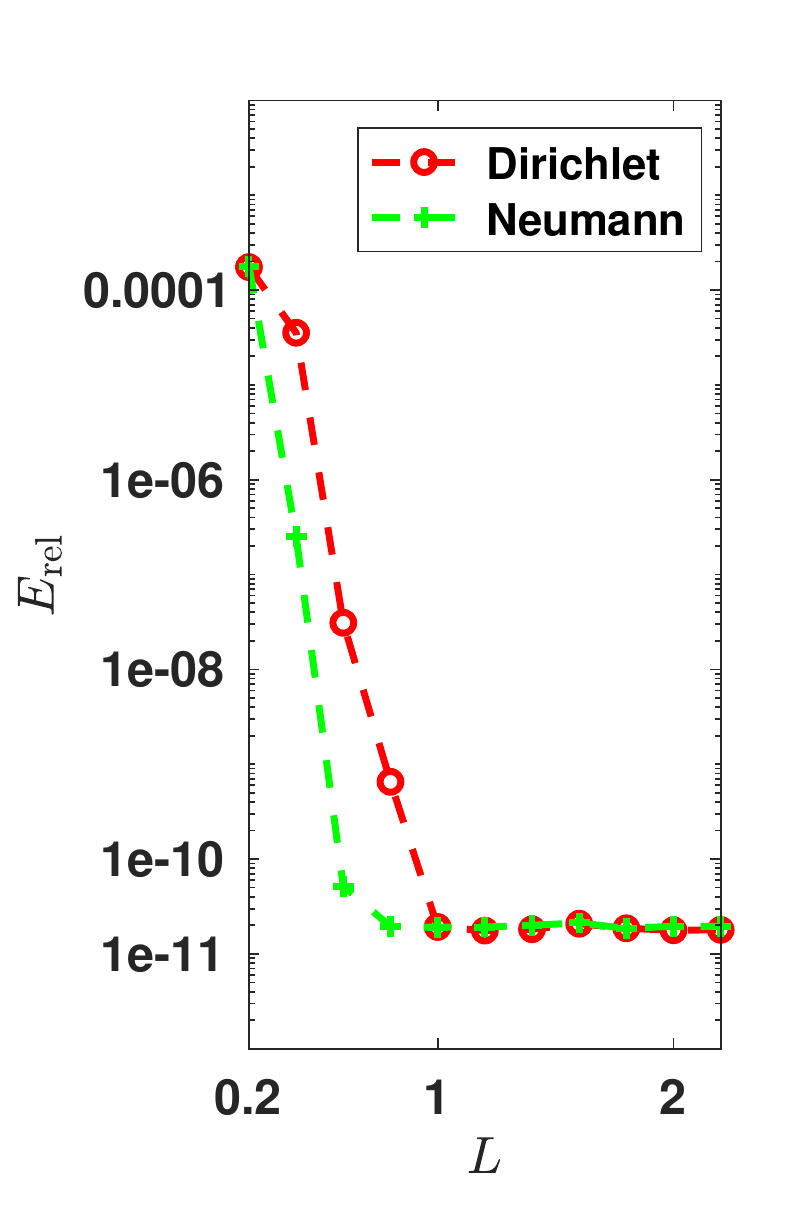}
	\caption{Example 1: Real part of $u^{\rm tot}$ in $[-5.5,5.5]\!\times\![-2.0,3.0]$
    excited by the point source ${y}\!=\!(0,1.5)$: (a) exact solution; (b) numerical
    solution. Convergence history of relative error {$E_{\rm rel}$} versus: (c) PML
    absorbing constant $S$; (d) Thickness of the PML $L$, for both Dirichlet and
    Neumann conditions on $\Gamma_{H+L}$.}
  \label{fig:ex1:1}
\end{figure}

\begin{figure}[!ht]
  \centering
  (a)\includegraphics[width=0.29\textwidth]{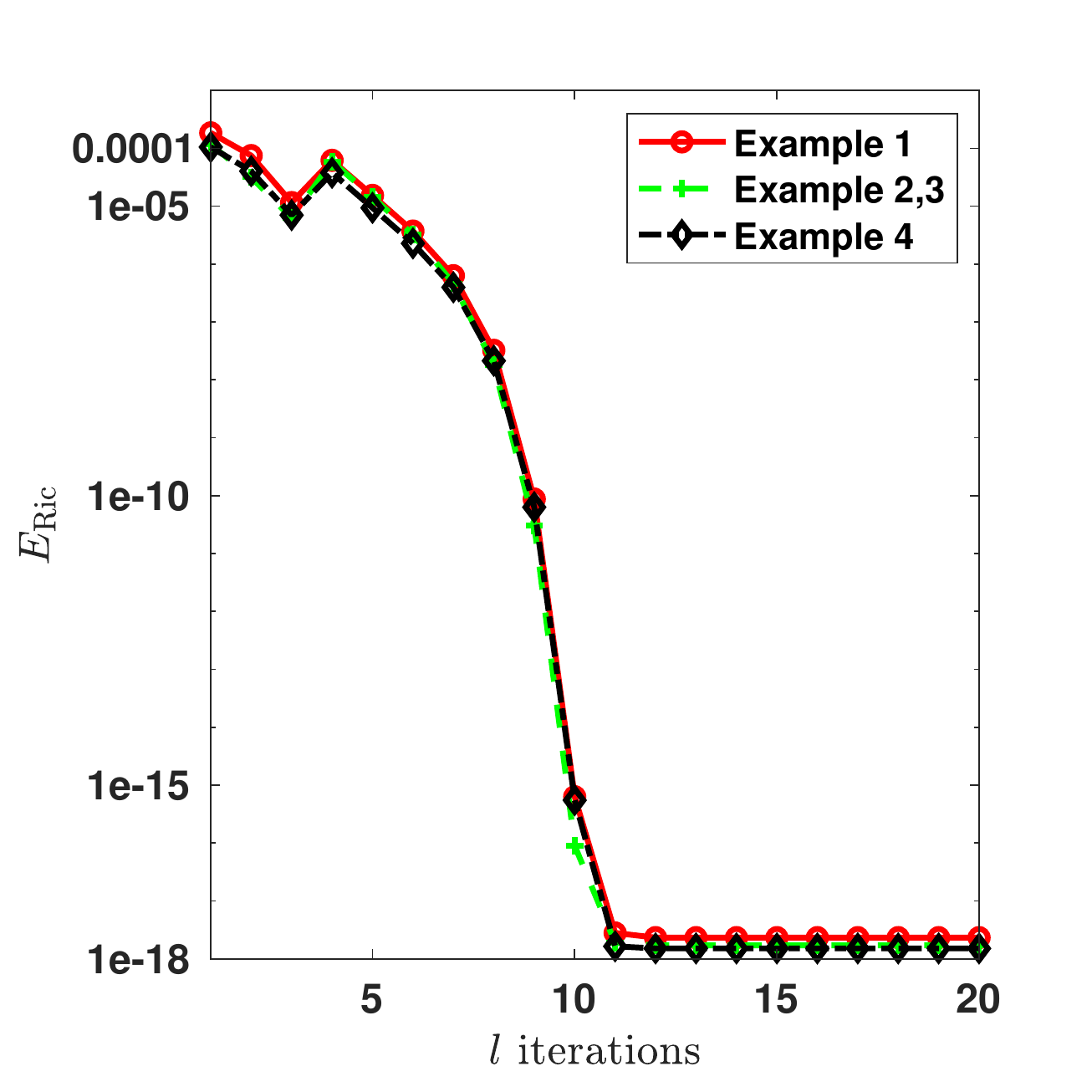}
  (b)\includegraphics[width=0.29\textwidth]{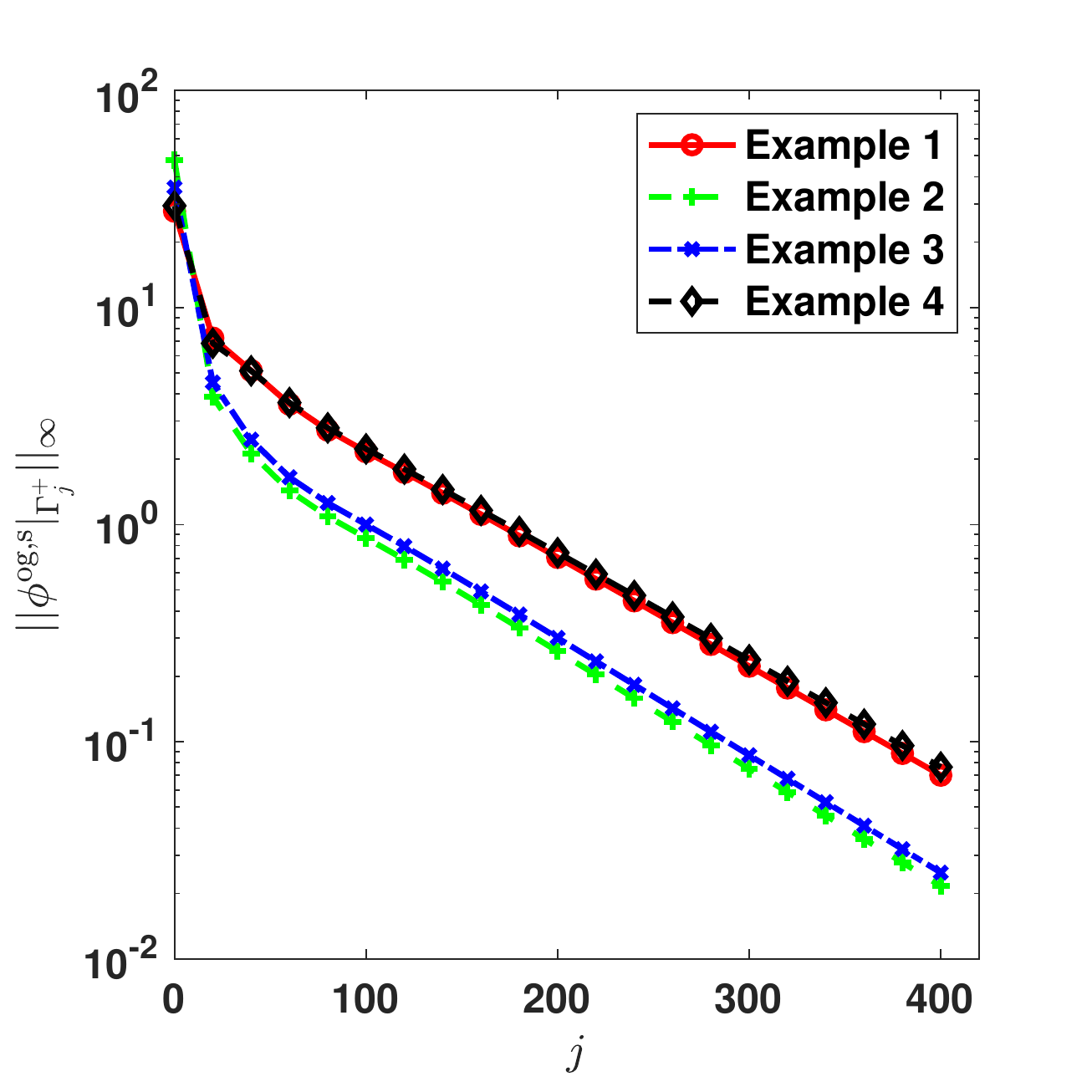}
  (c)\includegraphics[width=0.29\textwidth]{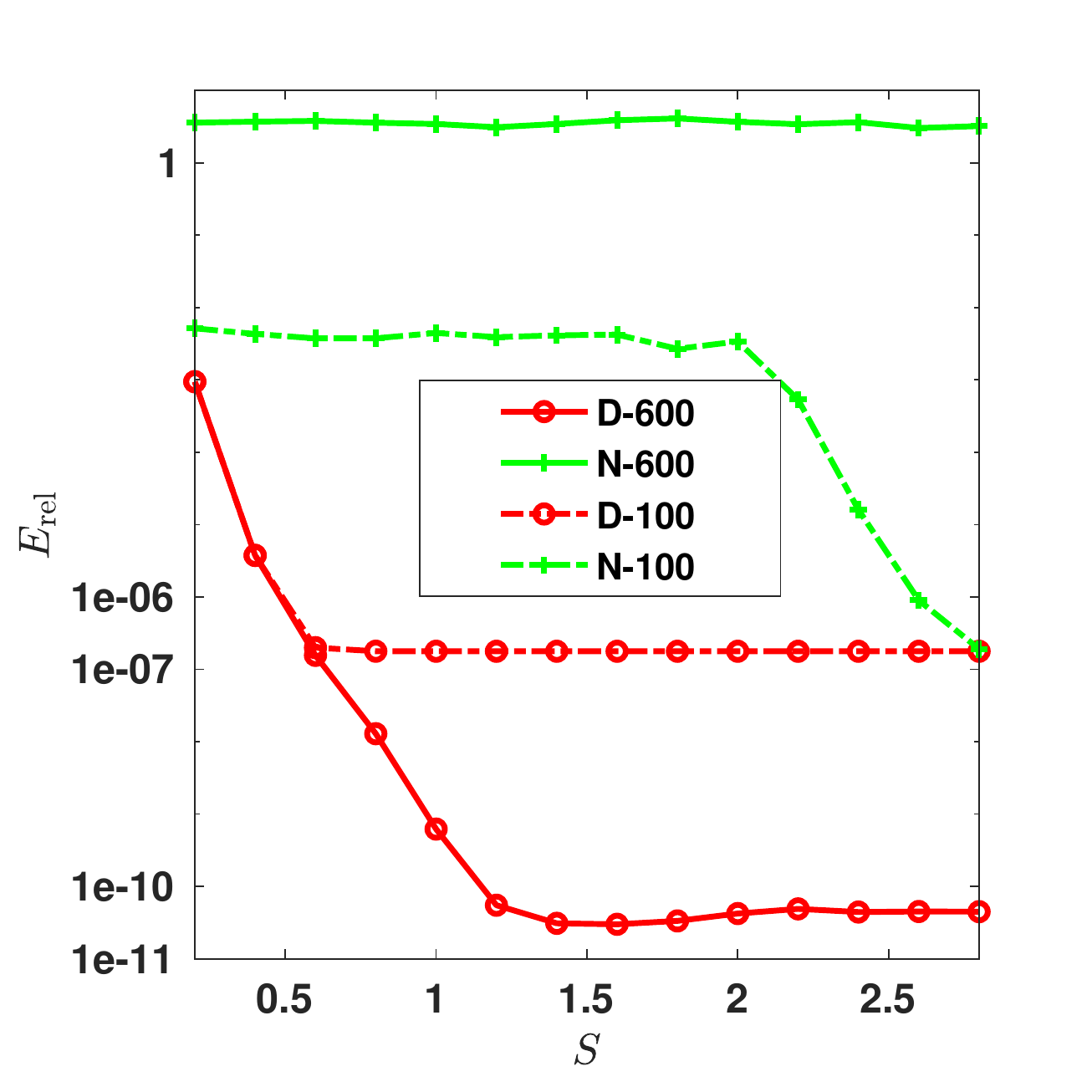}
	\caption{All four examples: (a) Convergence history of $E_{\rm Ric}$ in
    (\ref{eq:eric}) against the number of iterations $l$; (b) Radiation behavior
    of $\phi^{\rm og,s}|_{\Gamma_j^+}$ as $j\to\infty$. (c): Performance of
    Dirichlet and Neumann conditions in Example 1 for $n=1$, at which
    $kT/\pi\in{\mathbb Q}^+$; here 'D' stands for Dirichlet and 'N' for Neumann,
    and $100$ indicates $100$ grid points are used to discretize each smooth
    segment of the unit cells, etc.. }
  \label{fig:ex1:2}
\end{figure}

For $n=1.03$, {Figure}~\ref{fig:ex1:1} (a) and (b) compare the exact
solution and our numerical solution for $L\!=\!2.2$ and $S\!=\!2.8$. The two solutions are indistinguishable.
To give a detailed comparison, {Figure}~\ref{fig:ex1:1} (c) and (d) show how
{the relative error $E_{\rm rel}$} decays as one of the two PML parameters, the
absorbing constant $S$ and {the} thickness $L$, increases for either zero
Dirichlet or zero Neumann condition on $\Gamma_{H+L}$. In
{Figure}~\ref{fig:ex1:1}(c), we take $L\!=\!2.2$ and let $S$ vary between $0.2$
and $2.8$, while in {Figure}~\ref{fig:ex1:1}(d), we take $S\!=\!2.8$ and let
{the} PML thickness $L$ vary between $0.2$ and $2.2$. In both figures, the
vertical axis is logarithmically scaled so that {the vertical dashed lines}
indicate that the relative error {$E_{\rm rel}$} decays exponentially as {$L$}
or $S$ increases for both conditions. On the other hand, Neumann condition gives
faster convergence rate than Dirichlet condition. The convergence curves
indicate nearly $11$ significant digits are revealed by the proposed PML-based
BIE method. The 'o' lines in Figure~\ref{fig:ex1:2}(a) show the convergence
curve of
\begin{equation}
  \label{eq:eric}
  E_{\rm Ric} = ||{\bm N}^{(0)}_{10}+ [{\bm N}^{(0)}_{11}+{\bm
    N}^{(0)}_{00}]{\bm R}_p^++{\bm N}^{(0)}_{01}({\bm R}_p^+)^2||_{\infty}
\end{equation}
against the number of iterations $l$. It can be seen that after only $11$
iterations, ${\bm R}_p^+$ \tcr{satisfies} its governing \tcr{Riccati} equation
(\ref{eq:ric:Rp+}) up to round-off errors. The 'o' lines in
Figure~\ref{fig:ex1:2}(b) show the curve of $||\phi^{\rm og,s}|_{\Gamma_{j}^+}||_{\infty}$ against $j$. It can be seen that $\phi^{\rm
  og,s}$ and hence $\partial_{\nu_c^+}u^{\rm og}$ indeed decay exponentially as
$j$ or $x_1$ increases, as has been illustrated in Corollary~\ref{cor:exp:conv}.

In Figure~\ref{fig:ex1:2}(c), we compare Dirichlet and Neumann conditions for
$n=1$. We take $L=2.2$ and let $S$ vary from $0.2$ to $2.8$.
Among the four convergence curves, solid lines indicate $600$ grid points 
chosen on each smooth segment of each unit cell, while dashed lines indicate
$100$ grid points; '+' indicates Neumann condition on $\Gamma_{H+L}$ while 'o'
indicates Dirichlet condition. \tcr{If} $100$ grid points are used, $E_{\rm rel}$
for Neumann condition starts decreasing after $S\geq 2$ whereas $E_{\rm rel}$ for
Dirichlet condition has already reached its minimum error; \tcr{if} $600$ grid
points are used, Neumann condition \tcr{does not} make $E_{\rm rel}$ converge at all
for $S\in[0.2,2.8]$, but Dirichlet condition still possesses the same
convergence rate and accuracy as in case $n=1.03$. Consequently, Dirichlet
condition outperforms Neumann condition for $n=1$.

\begin{figure}[!ht]
  \centering
  (a)\includegraphics[width=0.28\textwidth]{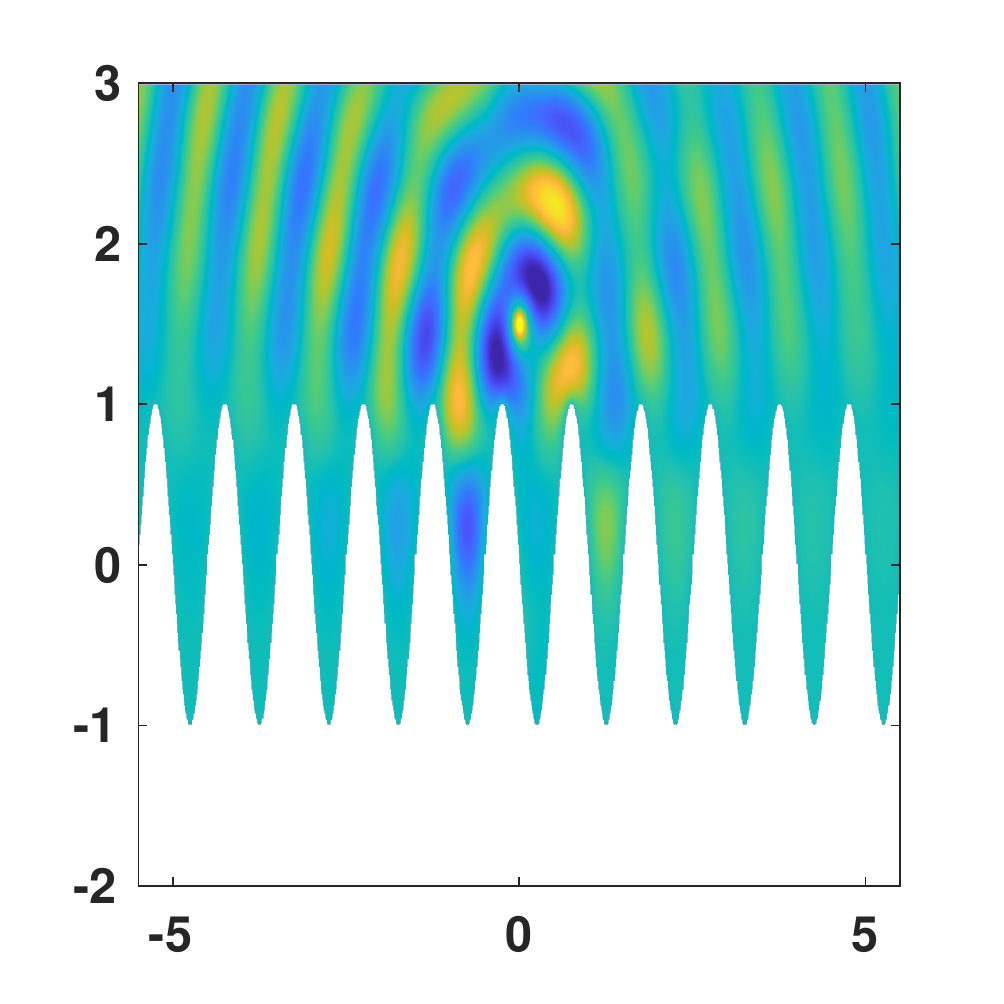}
  (b)\includegraphics[width=0.19\textwidth]{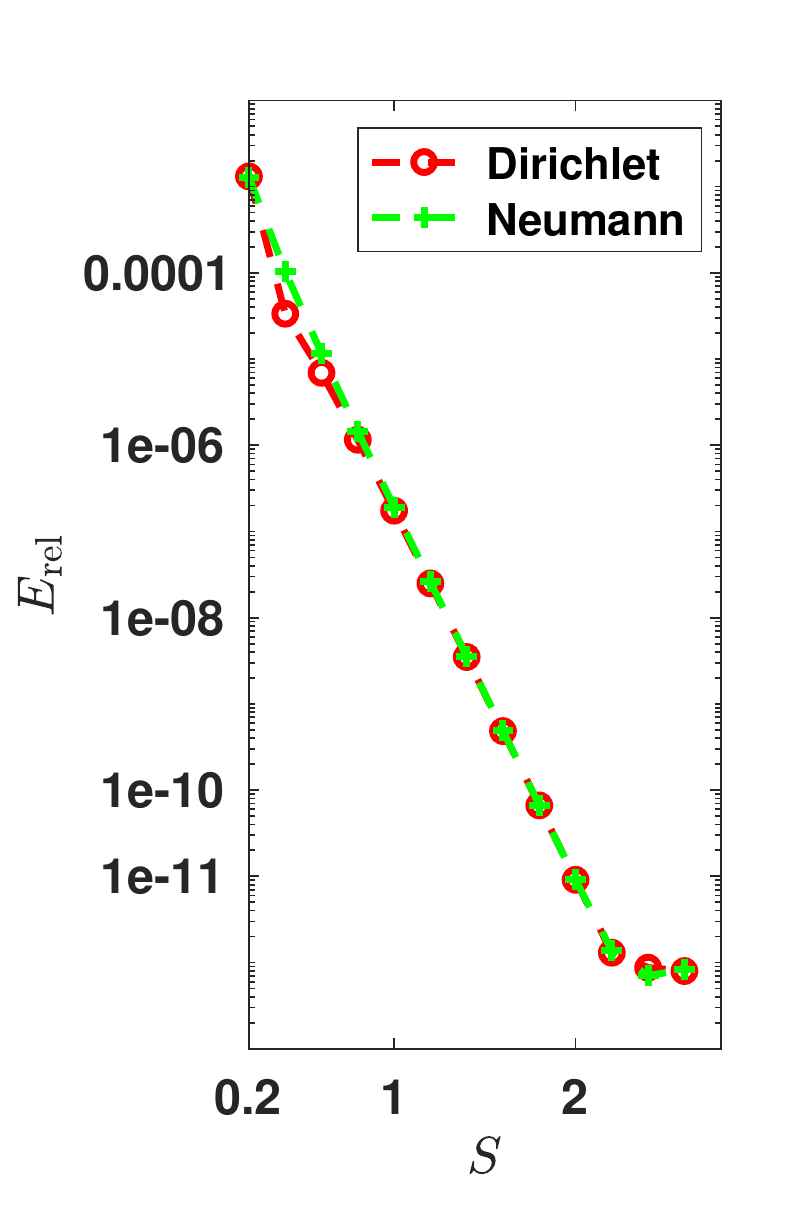}
  (c)\includegraphics[width=0.19\textwidth]{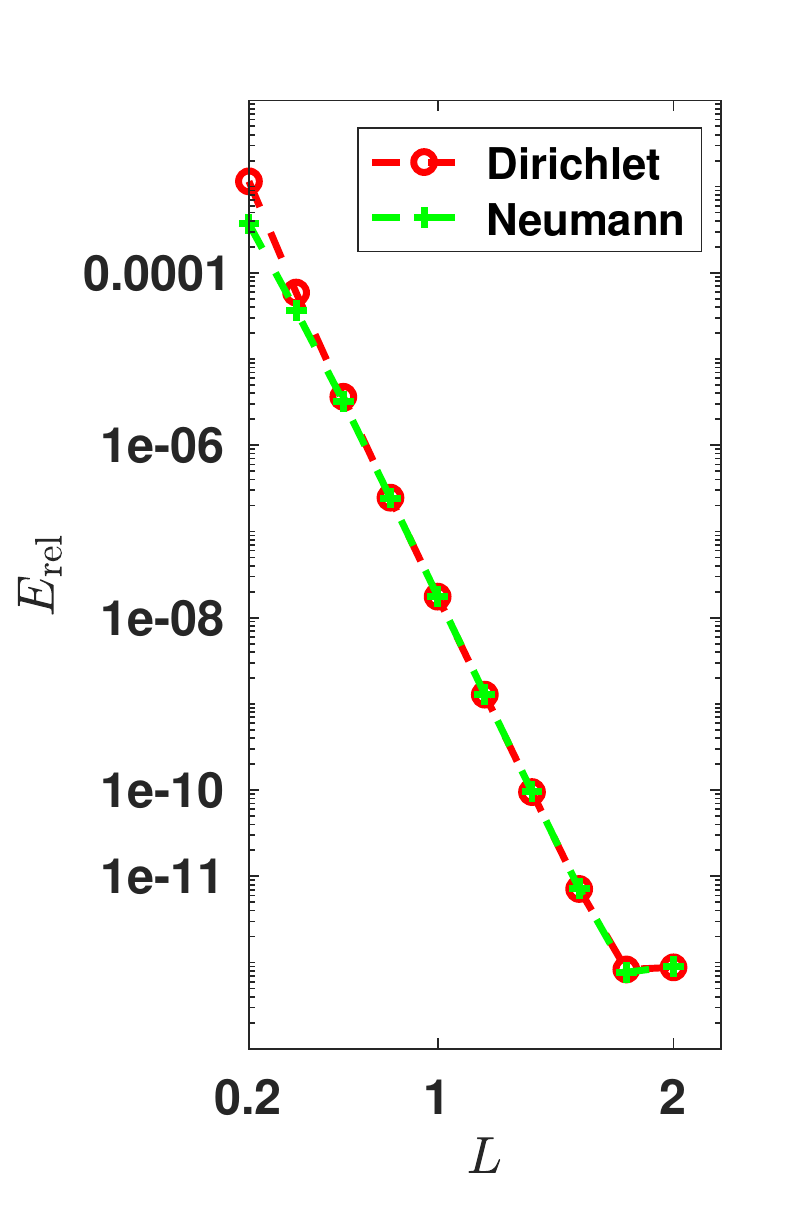}
	\caption{Example 2: (a) Numerical solution of {real} part of the total wave
    field $u$ in $[-5.5,5.5]\!\times \![-2.0,3.0]$ excited by the point source
    ${y}\!=\!(0,1.5)$. Convergence history of relative error {$E_{\rm rel}$}
    versus: (b) PML absorbing constant $S$ for fixed PML Thickness
    $\tcr{L}\!=\!2$, (c) PML Thickness $\tcr{L}$ for fixed PML absorbing constant
    $S\!=\!2.8$; vertical axes are logarithmically scaled.}
  \label{fig:ex2:ps}
\end{figure}

{\bf Example 2: a sine curve.} In the second example, we assume that $\Gamma$ is the sine curve,
$x_2\!=\!\sin(2\pi x_1+\pi)$, as shown in {Figure}~\ref{fig:ex2:ps}(a) and that
$n=1.03$ to make $kT/\pi\notin{\mathbb Q}^+$. 
For the cylindrical incidence, we discretize each smooth segment of any unit
cell by 600 grid points, and compare results of Dirichlet and Neumann boundary
conditions on $\Gamma_{H+L}$. Taking $S\!=\!2.8$ and $L\!=\!2.2$, we evaluate
the wave field in $[-5.5,5.5]\!\times\! [-2.0,3.0]$ and use this as the
reference solution {since the} exact solution is no longer {available.} In
{Figure}~\ref{fig:ex2:ps}, (a) shows the field pattern of the reference
solution, and (b) and (c) show the convergence history of relative error
{$E_{\rm rel}$} versus one of the two PML parameters $S$ and $L$, respectively.
Again, we observe that {$E_{\rm rel}$} decays exponentially {as} $S$ or $L$
increases. Unlike the flat surface in Example 1, we no longer observe a faster
convergence rate of Neumann condition, but find that both conditions share the
same convergence rate and accuracy. Considering its worse result for
$kT/\pi\in{\mathbb Q}$ and unimpressive improvement for $kT/\pi\notin{\mathbb
  Q}$, we conclude that Neumann condition is less superior than Dirichlet
condition, and thus shall only use the latter one in the rest experiments. With
Dirichlet condition, the '+' lines in Figure~\ref{fig:ex1:2} (a) show the
convergence curve of $E_{\rm Ric}$ in (\ref{eq:eric}) against the number of
iterations $l$. The '+' lines in Figure~\ref{fig:ex1:2} (b) show the curve of
$||\phi^{\rm og,s}|_{\Gamma_{j}^+}||_{\infty}$ against $j$.

\begin{figure}[!ht]
  \centering
  (a)\includegraphics[width=0.28\textwidth]{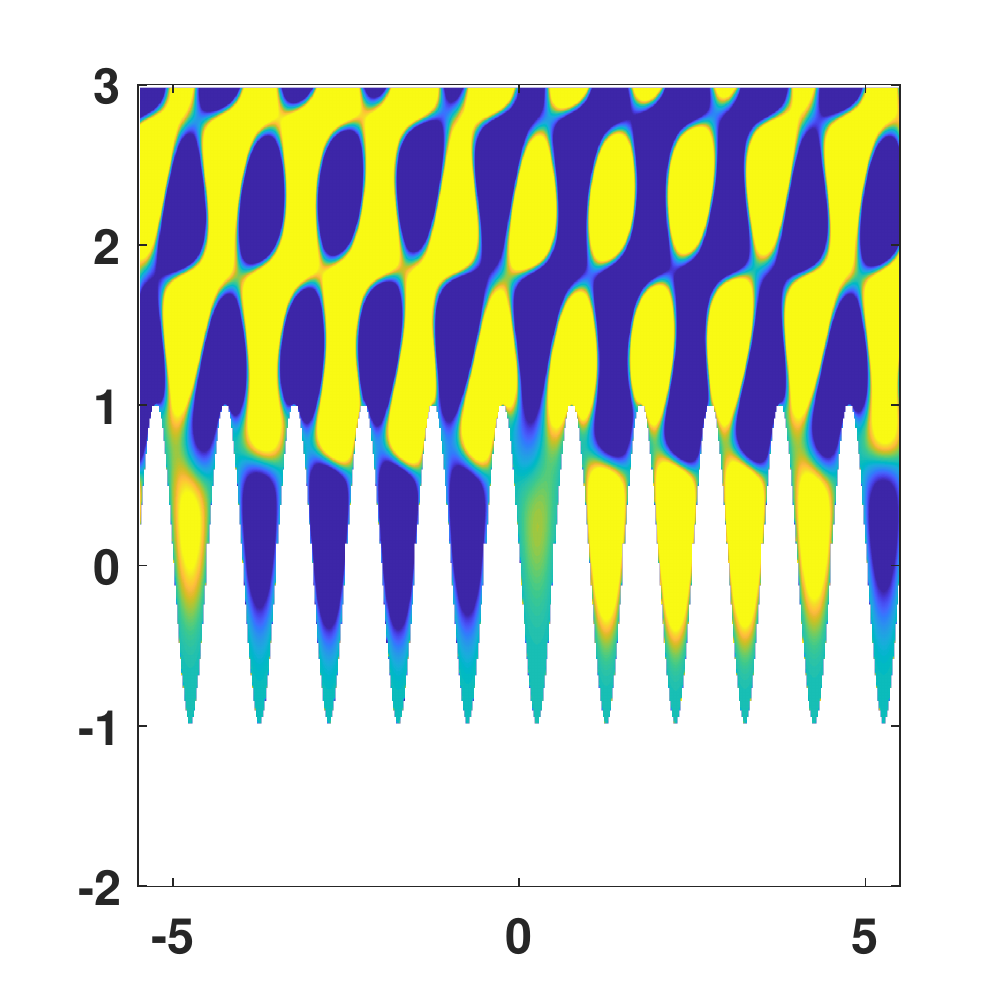}
  (b)\includegraphics[width=0.19\textwidth]{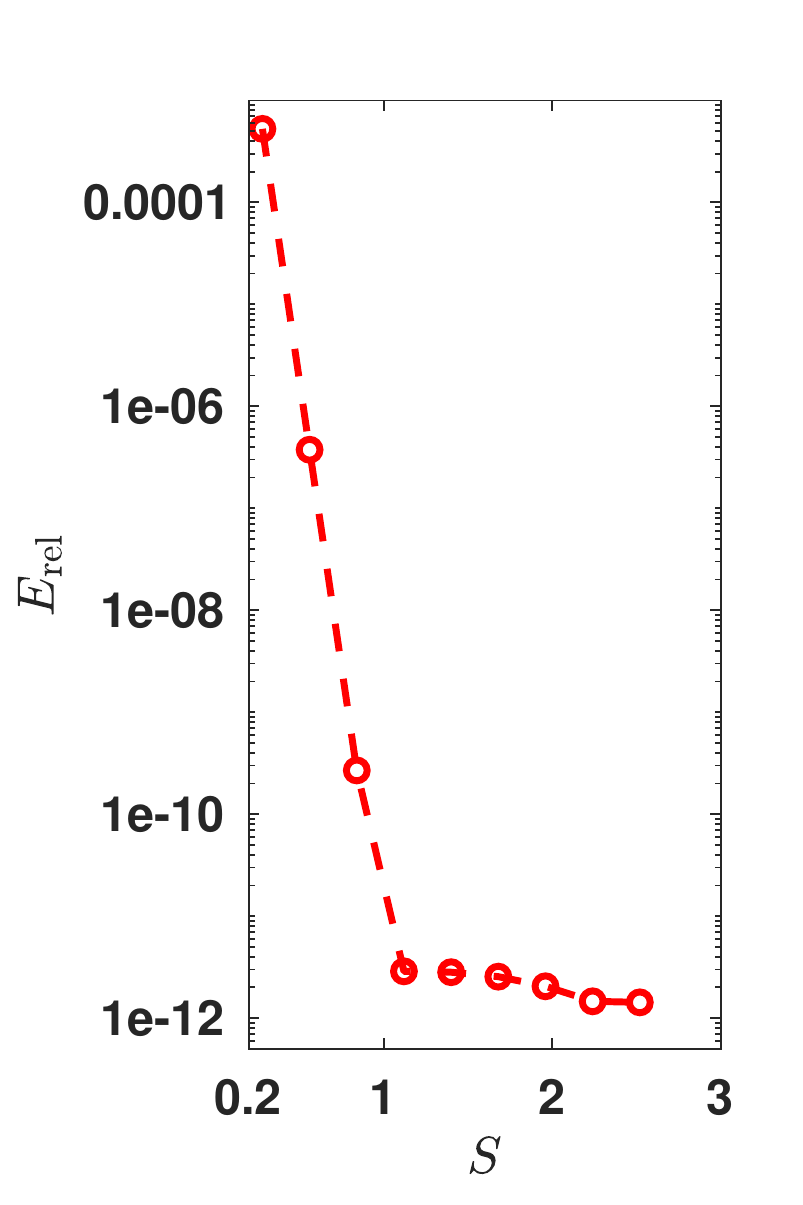}
  (c)\includegraphics[width=0.19\textwidth]{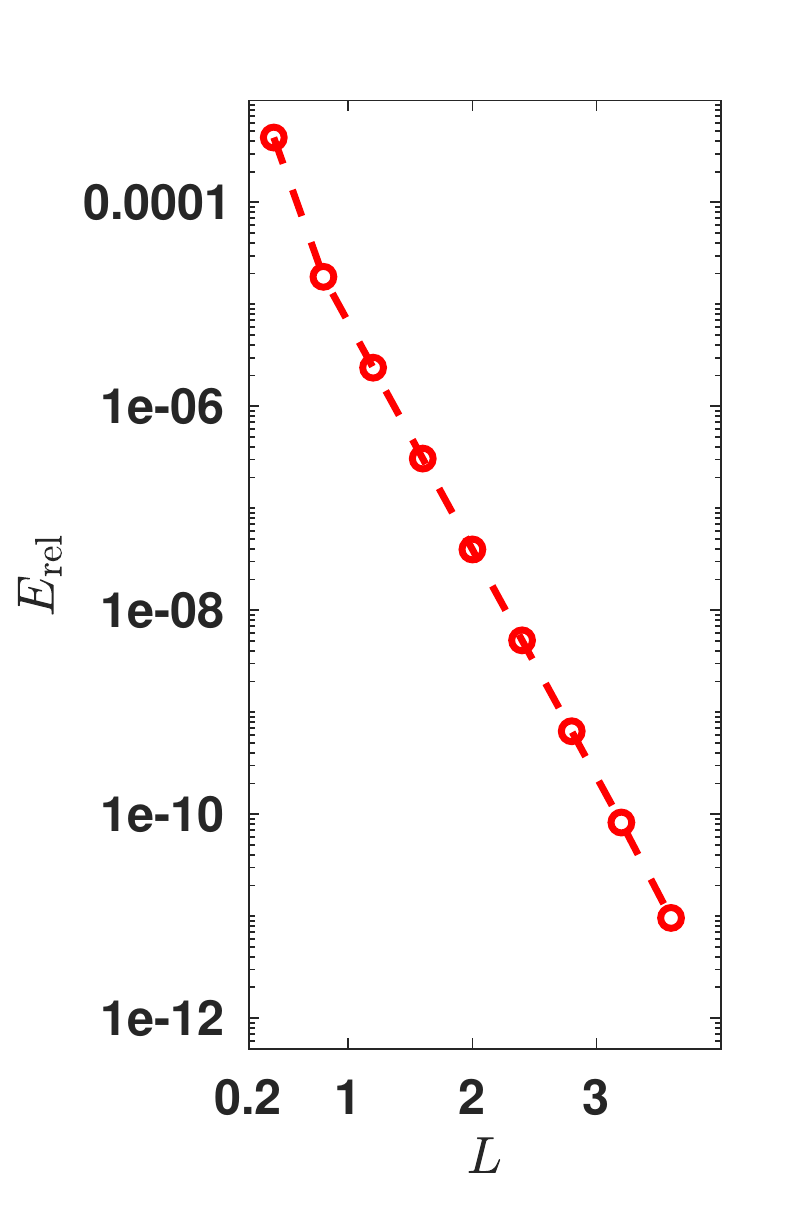}
	\caption{Example 2: (a) Numerical solution of {real} part of the
          total wave field $u$ in $[-5.5,5.5]\!\times \![-2.0,3.0]$
    excited by a plane incident wave of angle $\theta=\frac{\pi}{3}$. Convergence history of relative error {$E_{\rm rel}$} versus: (b)
    PML absorbing constant $S$ for fixed PML Thickness $\tcr{L}\!=\!4$, (c) PML
    Thickness $\tcr{L}$ for fixed PML absorbing constant $S\!=\!2.8$; vertical
    axes are logarithmically scaled.}
  \label{fig:ex2:pw}
\end{figure}

For the \tcr{plane-wave} incidence, we take $\theta=\frac{\pi}{3}$. Employing the method in
section 6.3, we discretize each smooth segment of any unit cell by 700 grid
points. Taking $S\!=\!2.8$ and $L\!=\!4$, we evaluate the wave field in
$[-5.5,5.5]\!\times\! [-2.0,3.0]$ and use this as the reference solution. In
{Figure}~\ref{fig:ex2:pw}, (a) shows the field pattern, and (b) and (c) show
the convergence history of relative error {$E_{\rm rel}$} versus one of the two
PML parameters $S$ and $L$, respectively. For both incidences, the convergence
curves in Figures~\ref{fig:ex2:ps} and \ref{fig:ex2:pw} decay exponentially,
indicating that nearly $12$ significant digits are revealed by the proposed
PML-based BIE method.

\begin{figure}[!ht]
  \centering
  (a)\includegraphics[width=0.23\textwidth]{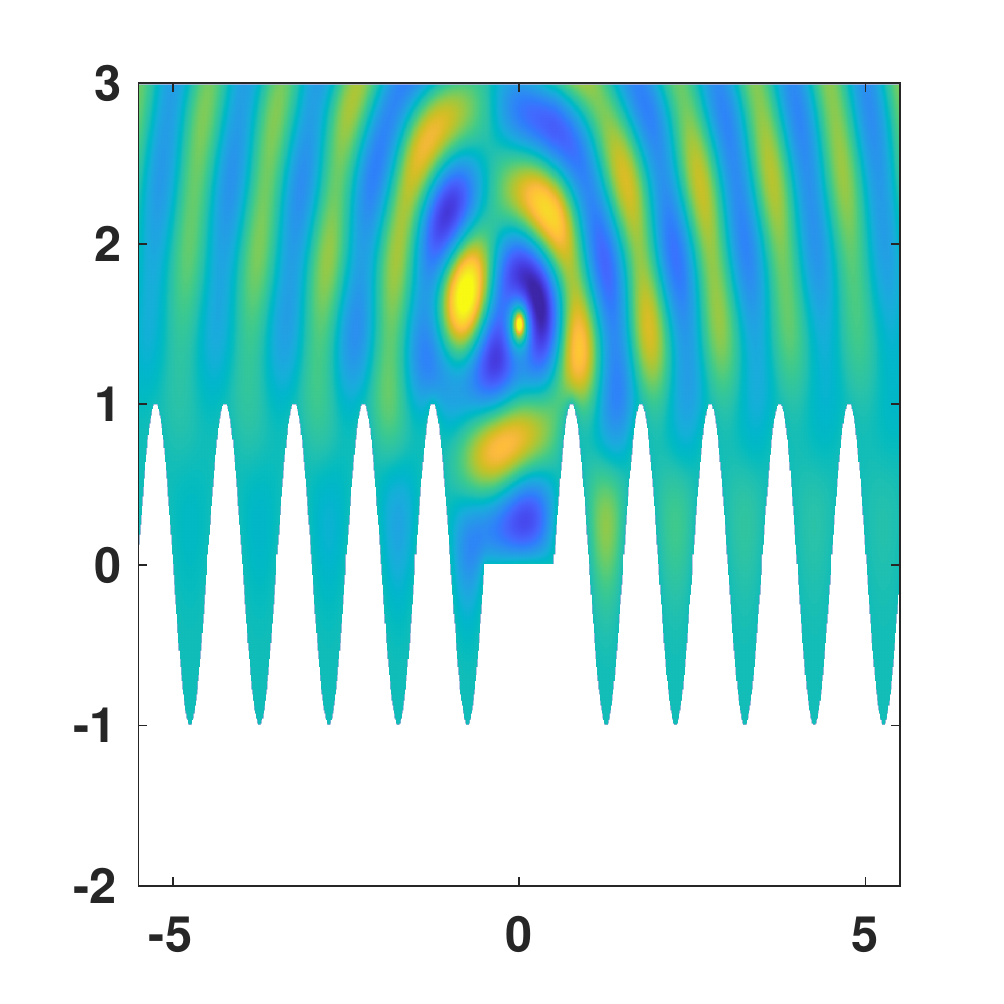}
  (b)\includegraphics[width=0.23\textwidth]{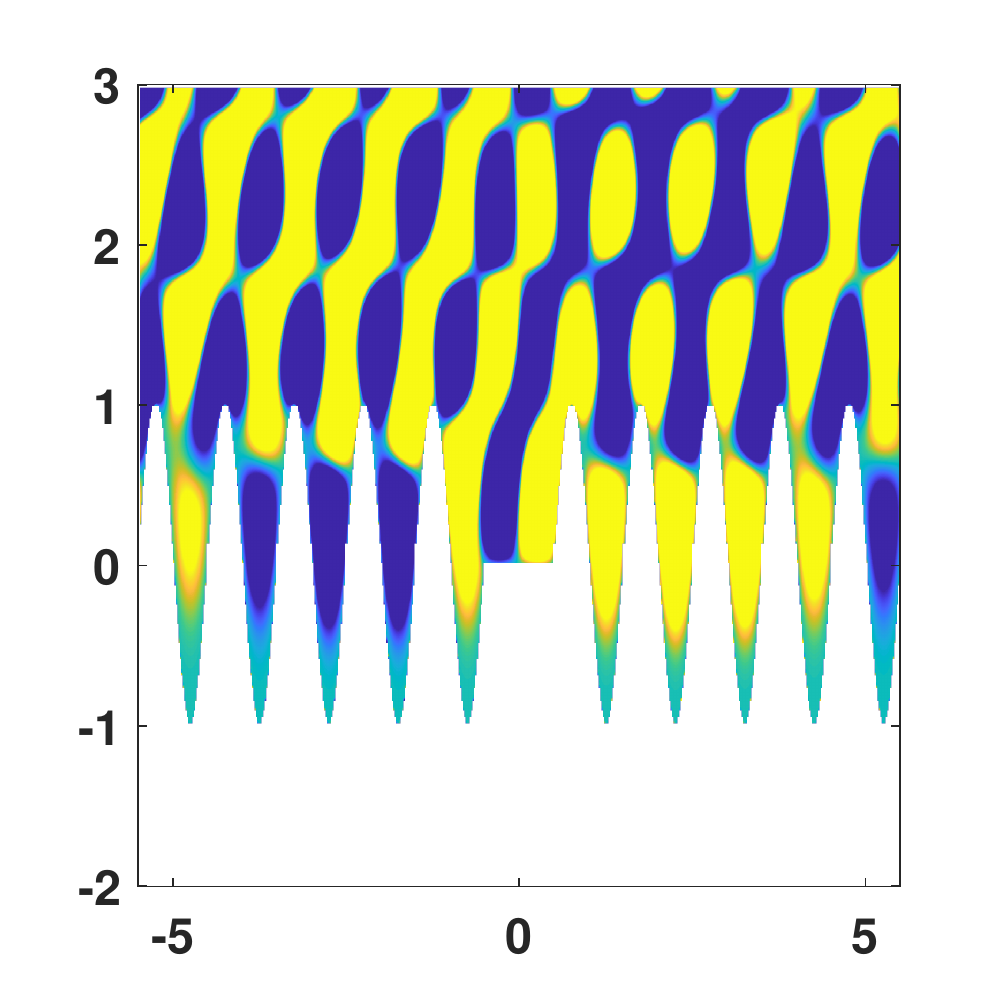}
  (c)\includegraphics[width=0.19\textwidth]{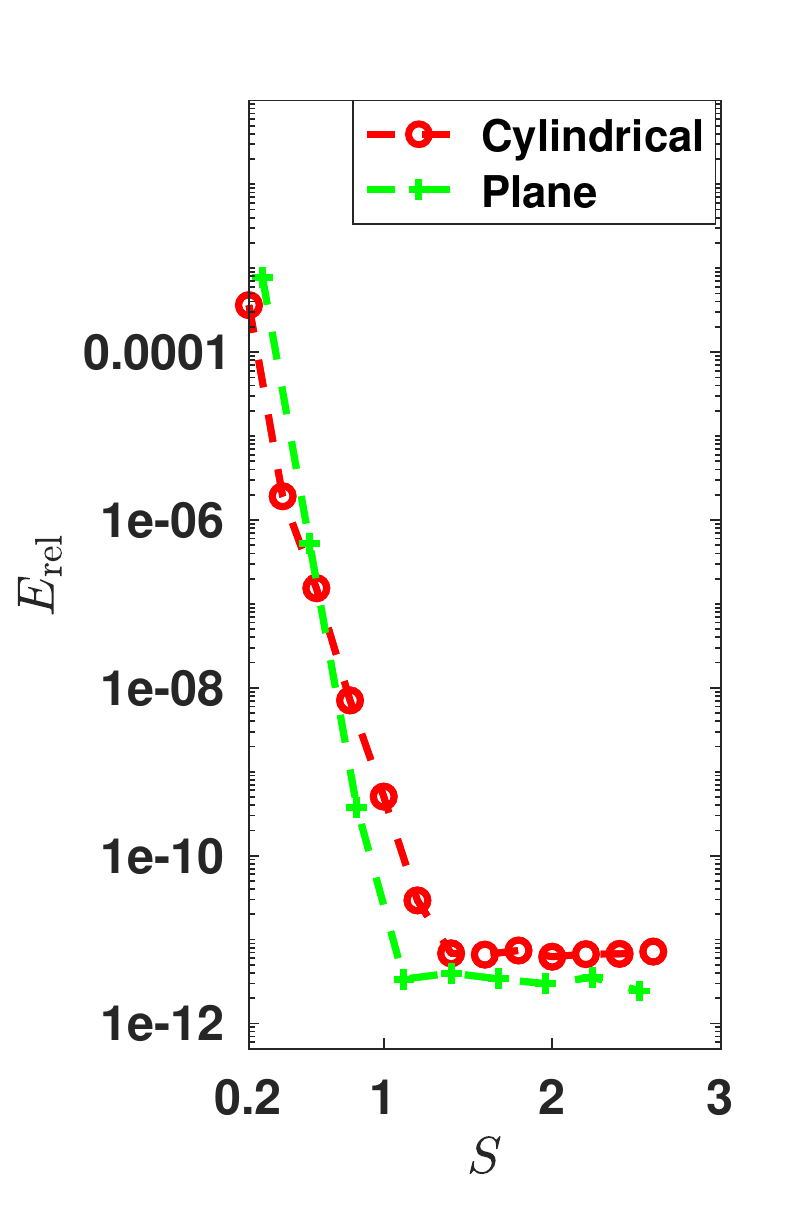}
  (d)\includegraphics[width=0.19\textwidth]{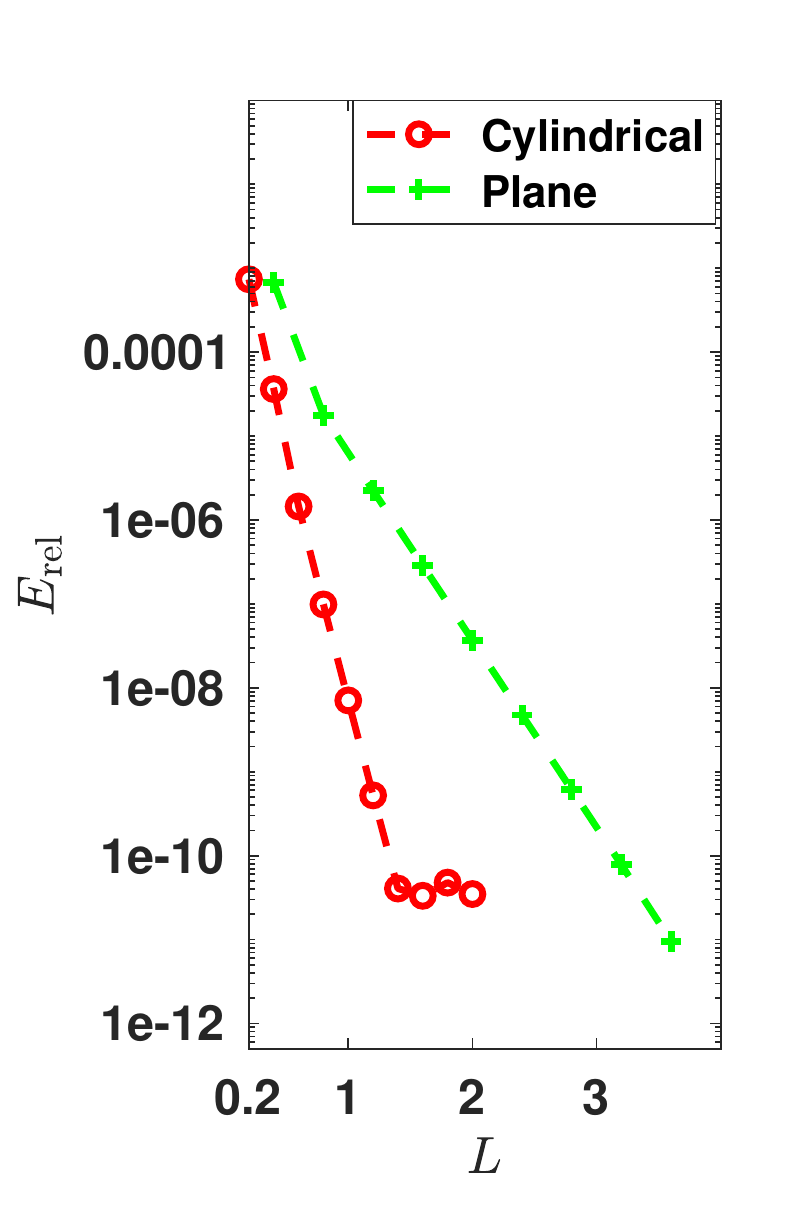}
	\caption{Example 3: Numerical solution of {real} part of the
          total wave field $u$ in $[-5.5,5.5]\!\times\![-2.0,3.0]$
    excited by: (a) a cylindrical wave by source ${y}\!=\!(0,1.5)$; (b) a plane
    wave of incident angle $\theta=\frac{\pi}{3}$. Convergence history of relative error {$E_{\rm rel}$} versus: (c) PML
    Thickness $L$ for fixed PML absorbing constant $S\!=\!2.8$ for both incidences; (d)
    PML absorbing constant $S$ for fixed PML Thickness $L\!=\!2.2$ ($4.0$) for
    cylindrical (\tcr{plane-wave}) incidence.}
  \label{fig:ex3:ps}
\end{figure}

{\bf Example 3: a locally perturbed sine curve.} In the third example, we assume
that the sine curve $\Gamma: x_2\!=\!\sin(2\pi x_1+\pi)$ is locally perturbed
with the part between $x_1\!=\!-0.5$ and $x_1\!=\!0.5$
replaced by {the} line segment $\{(x_1,0)\!: x_1\!\in\![-0.5,0.5]\}$, as shown in
{Figure}~\ref{fig:ex3:ps} (a).
For the cylindrical incidence, we discretize each smooth segment of any unit
cell by 600 grid points. Taking $S\!=\!2.8$ and $L\!=\!2.2$, we
evaluate the wave field in $[-5.5,5.5]\!\times\! [-2.0,3.0]$ and use this as the
reference solution, the field pattern of which is shown in 
{Figure}~\ref{fig:ex3:ps} (a).  The 'x' lines in
Figure~\ref{fig:ex1:2} (b) show the curve of $||\phi^{\rm
  og,s}|_{\Gamma_{j}^+}||_{\infty}$ against $j$.

For the plane incidence, we take $\theta=\frac{\pi}{3}$ and discretize each smooth segment of any unit cell by
700 grid points. Taking $S\!=\!2.8$ and $L\!=\!4$, we evaluate the wave field
in $[-5.5,5.5]\!\times\! [-2.0,3.0]$ and use this as the reference solution, the
field pattern of which is shown in {Figure}~\ref{fig:ex3:ps} (b).

For both incidences, {Figure}~\ref{fig:ex3:ps} (c) and (d) 
show the convergence history of relative error {$E_{\rm rel}$} versus one of the
two PML parameters $S$ and $L$, respectively. The convergence curves decay
exponentially and indicate that nearly $11$ significant digits are revealed by the
proposed PML-based BIE method.

\begin{figure}[!ht]
  \centering
  (a)\includegraphics[width=0.23\textwidth]{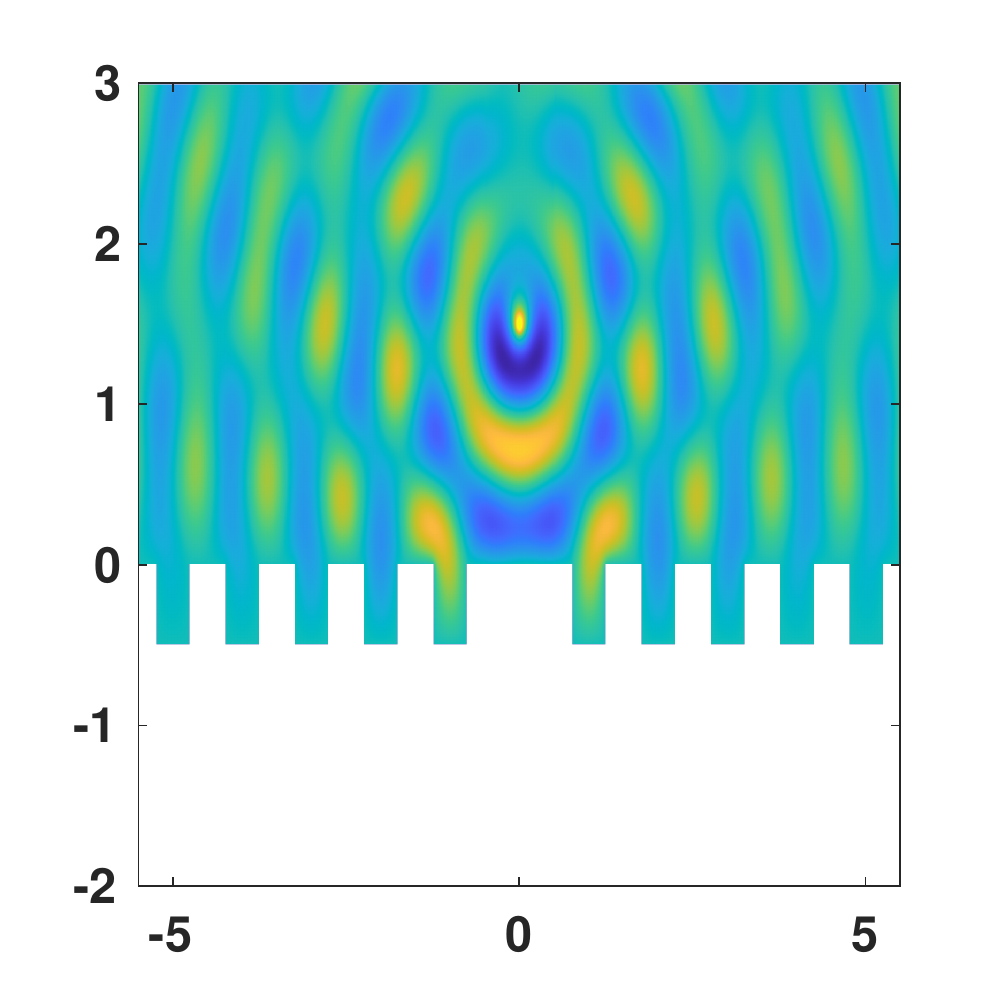}
  (b)\includegraphics[width=0.23\textwidth]{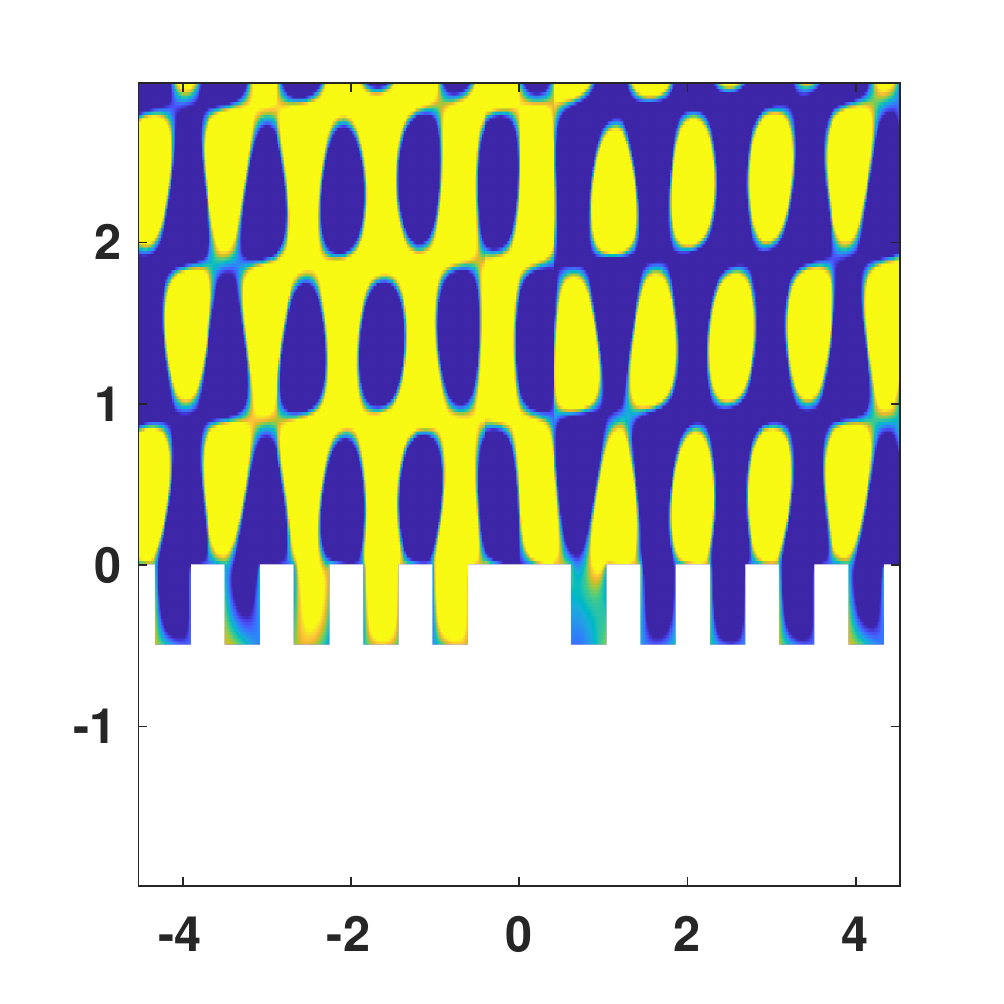}
  (c)\includegraphics[width=0.19\textwidth]{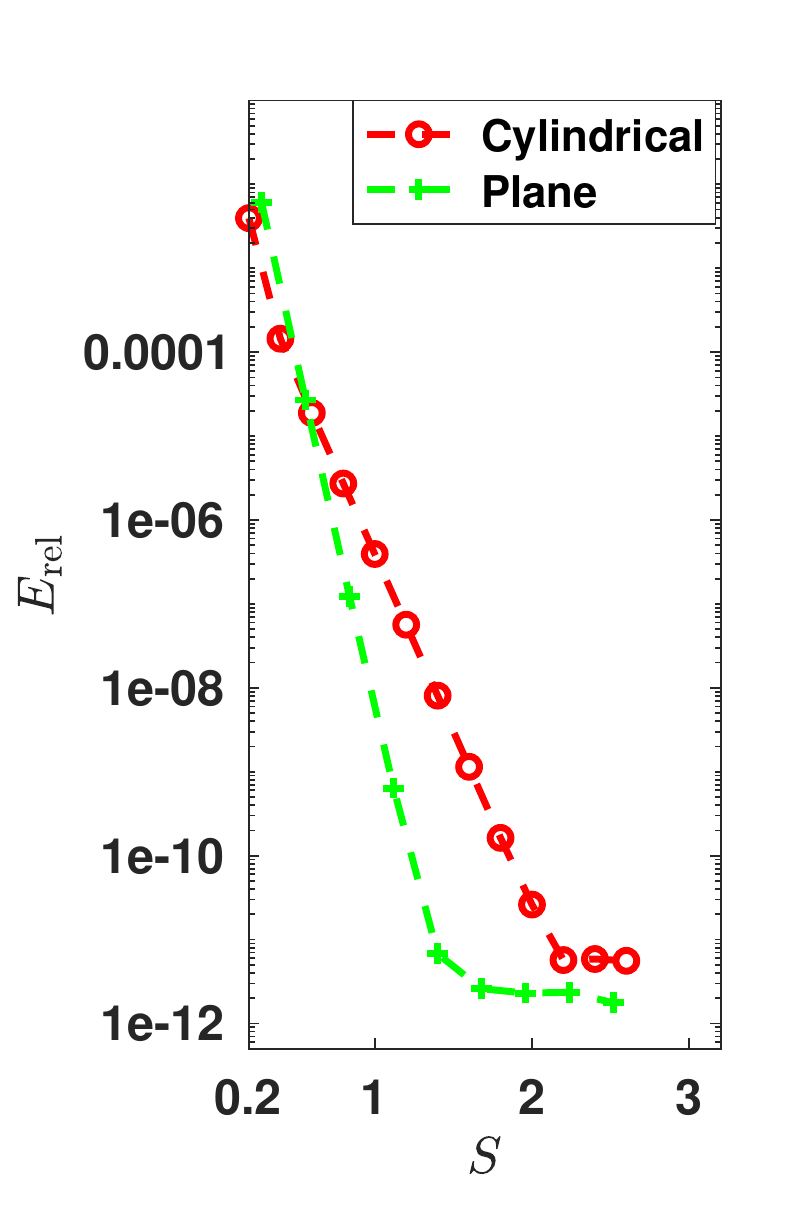}
  (d)\includegraphics[width=0.19\textwidth]{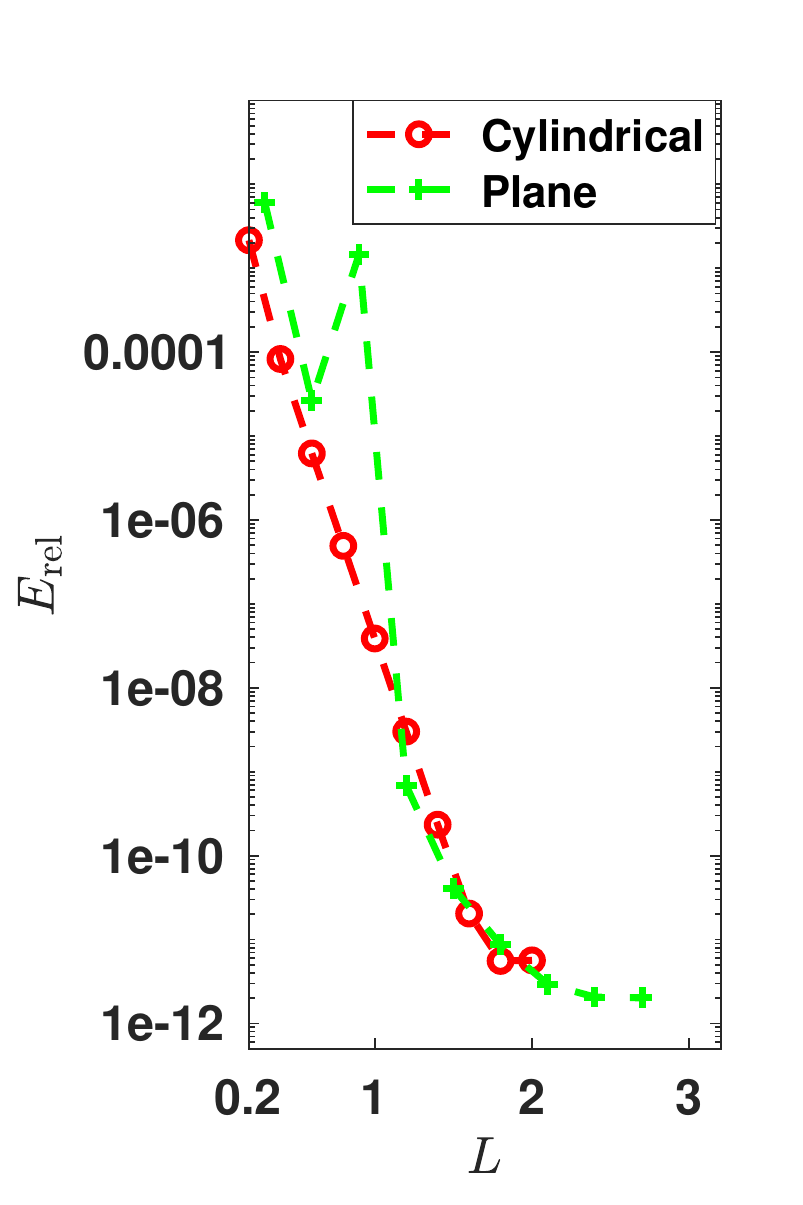}
	\caption{Example 4: Numerical solution of {real} part of the
          total wave field $u$ in $[-5.5,5.5]\!\times\![-2.0,3.0]$
    excited by: (a) a cylindrical wave by source ${y}\!=\!(0,1.5)$; (b) a plane
    wave of incident angle $\theta=\frac{\pi}{6}$. Convergence history of relative error {$E_{\rm rel}$} versus: (c) PML
    Thickness $L$ for fixed PML absorbing constant $S\!=\!2.8$ for both incidences; (d)
    PML absorbing constant $S$ for fixed PML Thickness $L\!=\!2.2$ ($3.0$) for
    cylindrical (\tcr{plane-wave}) incidence.}
  \label{fig:ex4:ps}
\end{figure}

{\bf Example 4: \tcr{a locally perturbed binary grating.}} In the last example, we assume that $\Gamma$ consists of
periodic rectangular grooves of depth $0.5$ and width $0.25$, with the
part between $x_1\!=\!-0.5$ and $x_1\!=\!0.5$
replaced by {the} line segment $\{(x_1,0)\!: x_1\!\in\![-0.5,0.5]\}$, as shown in
{Figure}~\ref{fig:ex4:ps}(a).
For the cylindrical incidence, we discretize each smooth segment of any unit
cell by 600 grid points. Taking $S\!=\!2.8$ and $L\!=\!2.2$, we
evaluate the wave field in $[-5.5,5.5]\!\times\! [-2.0,3.0]$ and use this as the
reference solution, the field pattern of which is shown in 
{Figure}~\ref{fig:ex4:ps} (a). The  '$\Diamond$' lines in
Figure~\ref{fig:ex1:2} (a) show the convergence curve of $E_{\rm Ric}$ in
(\ref{eq:eric}) against the number of iterations $l$. The '$\Diamond$' lines in
Figure~\ref{fig:ex1:2} (b) show the curve of $||\phi^{\rm
  og,s}|_{\Gamma_{j}^+}||_{\infty}$ against $j$.

For the \tcr{plane-wave} incidence, we take $\theta=\frac{\pi}6$ and discretize each smooth segment of any unit cell by
600 grid points. Taking $S\!=\!2.8$ and $L\!=\!3$, we evaluate the wave field
in $[-5.5,5.5]\!\times\! [-2.0,3.0]$ and use this as the reference solution, the
field pattern of which is shown 
in {Figure}~\ref{fig:ex4:ps} (b).

For both incidences, {Figure}~\ref{fig:ex4:ps} (c) and (d)
show the convergence history of relative error {$E_{\rm rel}$} versus one of the
two PML parameters $S$ and $L$, respectively.  The
convergence curves decay exponentially and 
indicate that nearly $12$ significant digits are revealed by the proposed PML-based
BIE method.

\section{Conclusion}
This paper studied the perfectly-matched-layer (PML) theory for wave scattering
in a half space of homogeneous medium bounded by a two-dimensional, perfectly
conducting, and locally defected periodic surface, and developed a high-accuracy
boundary-integral-equation (BIE) solver. By placing a PML in the vertical
direction to truncate the unbounded domain to a strip, we proved that the PML
solution converges to the true solution in the physical subregion of the strip
at an algebraic order of the PML thickness. Laterally, the unbounded strip is
divided into three regions: a region containing the defect and two
semi-waveguide regions of periodic subsurfaces, separated by two vertical line
segments. We proved the well-posedness of an associated scattering problem in
\tcr{both} semi-waveguide so as to well define a Neumann-to-Dirichlet (NtD) operator
on the associated vertical segment. The two NtD operators, serving as exact
lateral boundary conditions, reformulate the unbounded strip problem as a
boundary value problem \tcr{over} the defected region. \tcr{Each} NtD operator is closely
related to a Neumann-marching operator, governed by a nonlinear \tcr{Riccati}
equation, which was efficiently solved by \tcr{an} RDP method and a
high-accuracy PML-based BIE method so that the boundary value problem on the defected
region can be solved finally. Our future research plan shall focus on the following two aspects:
\begin{itemize}
  \item[(1).] Extend the current work to study locally defected periodic
    structures of stratified media. In such case, propagating Bloch modes may
    exist so that the related Neumann marching operators ${\cal R}_p^{\pm}$
    may not be contracting. 
  \item[(2).] Rigorously justify that the PML solution
    converges exponentially to the true solution in any compact subset of the
    strip, \tcr{as has been demonstrated by numerical experiments.}
\end{itemize}

\bibliographystyle{plain}
\bibliography{wt}
\end{document}